\documentclass[11pt]{elsarticle}
\usepackage{amsmath,amssymb,amsthm}
\usepackage{hyperref}
\usepackage{mathrsfs}
\usepackage{graphicx}
\usepackage{tikz}
\usepackage{enumerate}
\usepackage{float}
\usepackage{caption}
\usepackage{ltxtable}

\theoremstyle{plain}
\newtheorem{theorem}{Theorem}[section]
\newtheorem{lemma}[theorem]{Lemma}

\newtheorem{corollary}[theorem]{Corollary}

\theoremstyle{definition}
\newtheorem{definition}[theorem]{Definition}
\newtheorem{remark}[theorem]{Remark}
\newtheorem{example}[theorem]{Example}
\newtheorem{problem}[theorem]{Problem}

\newcommand{\ua}{\mathord{\uparrow}}
\newcommand{\da}{\mathord{\downarrow}}
\newcommand{\rom}[1]{\rm{\uppercase\expandafter{\romannumeral #1}}}

\makeatletter
\def\ps@pprintTitle{%
  \let\@oddhead\@empty
  \let\@evenhead\@empty
  \def\@oddfoot{\reset@font\hfil\thepage\hfil}
  \let\@evenfoot\@oddfoot
}
\makeatother
\makeatletter
\def\hlinew#1{%
  \noalign{\ifnum0=`}\fi\hrule \@height #1 \futurelet
   \reserved@a\@xhline}
\makeatother%

\begin{document}

\begin{frontmatter}

\title{Not every countable complete lattice is sober \tnoteref{t1}}
\tnotetext[t1]{This work is supported by the National Natural Science Foundation of China (No.11771134) and by Hunan Provincial Innovation Foundation For Postgraduate (CX20200419)}

\author[1]{Hualin Miao}
\ead{miaohualinmiao@163.com}
\author[2]{Xiaoyong Xi\corref{a1}}
\cortext[a1]{Corresponding author.}
\ead{littlebrook@jsnu.edu.cn}
\author[1]{Qingguo Li}
\ead{liqingguoli@aliyun.com}
\author[3]{Dongsheng Zhao}
\ead{dongsheng.zhao@nie.edu.sg}
\address[1]{School of Mathematics, Hunan University, Changsha, Hunan, 410082, China}
 \address[2]{School of mathematics and Statistics, Jiangsu Normal University, Xuzhou 221116, China}
\address[3]{School of Mathematics, Nanyang Technological University, 1 Nanyang Walk, Singapore 637616}

\begin{abstract}
The study of the sobriety of Scott spaces has got an relative long history in domain theory. Lawson and Hoffmann independently  proved that the Scott space of every continuous directed complete poset (usually called domain) is sober.  Johnstone constructed the first directed complete poset whose Scott space is non-sober. Not long after, Isbell gave a complete lattice with non-sober Scott space. Based on Isbell's example, Xu, Xi and Zhao showed that there is even a complete Heyting algebra whose Scott space is non-sober. Achim Jung then asked whether every countable complete lattice has a sober Scott space.

Let $\Sigma P$ be the Scott space of poset $P$. In this paper, we first prove that the topology of the product space $\Sigma P\times \Sigma Q$   coincides with the Scott topology on the product poset  $P\times Q$  if the set $Id(P)$ and $Id(Q)$ of all non-trivial ideals of posets $P$ and $Q$ are both countable. Based on this result, we deduce that a directed complete poset $P$ has a sober Scott space,  if $Id(P)$ is countable and the space $\Sigma P$  is coherent and well-filtered. Thus a complete lattice $L$ with $Id(L)$ countable has a  sober Scott space.
Making use the obtained results, we then construct  a countable complete lattice whose Scott space is non-sober and thus give a negative answer to Jung's problem.
\end{abstract}

\end{frontmatter}

\section{Introduction}

Sobriety is one of the earliest  studied major properties of $T_0$ topological spaces. It has been used in the characterization of spectra spaces  of commutative rings. In the recent years, this property and some of its weaker forms have been extensively investigated from various different perspectives. The Scott topology is the most important topology in domain theory which bridges a strong link between topological  and order structures. Lawson \cite{Lawson-1979} and Hoffmann\cite{Hoffmann-1981b} proved independently that the Scott space of every domain (continuous directed complete poset)
is sober. At the early time, it was an open problem whether the Scott space of every directed complete poset (dcpo, for short)  is sober. Johnstone constructed the first counterexample  to give a negative answer\cite{Johnstone-81}. Soon, Isbell\cite{Isbell} came up with a complete lattice whose Scott space is non sober. However, Isbell's complete lattice is neither distributive nor countable.

A poset $P$ will be called sober if its Scott space $\Sigma P$ is sober.
 In \cite{Jung}, Achim Jung posed  two problems. One of the them  is  whether every  distributive complete lattice is sober.

By using  Isbell's complete lattice,  Xu, Xi and Zhao gave  a negative answer to  this problem \cite{Xiaoquan22}.

The second  problem by Jung (also mentioned by Xu and Zhao in \cite{Xiaoquan})  is the following one:

\begin{problem}\label{a}
  Is there a non-sober countable complete lattice?
\end{problem}

In the current paper, we will give an  answer to this  problem.

Here is the outline of our paper.

For each poset $P$, let $\sigma(P)$ be the Scott topology on $P$ and $\Sigma P=(P, \sigma(P))$ be the Scott space of $P$ (see Section 1 for definitions).  It is well-known that for two posets $P$ and $Q$, the product topology on $\Sigma P\times \Sigma Q$ is usually  strictly coarser than the Scott topology $\sigma(P\times Q)$ on the product poset $P\times Q$. In Section 2,  we propose  a condition that guarantee the coincides of  $\sigma(P\times Q)$ and the topology of the product space $\Sigma(P)\times \Sigma(Q)$.
We then  prove that for a  dcpo $P$, $\Sigma P$  is sober if  $\Sigma P$ is coherent and well-filtered, and the set $Id(P)$ of all non-trivial ideals of posets $P$ is countable.

By \cite{xiaodong2} and \cite{Xi03}, it is known that the Scott space of a complete lattice is both well-filtered and coherent.
It then follows that if a complete lattice $L$ has only a countable number of non-trivial ideals, then $\Sigma L$ is sober.

 In Section 3, we construct a non-sober countable complete lattice, thus answering the second problem by Jung. In addition,  a countable distributive non-sober complete lattice is obtained  by modifying the constructed countable non-sober complete lattice.

\section{Preliminaries}
In this section we collect some basic definitions and results to be used later.
 For more details on them, we refer the reader to \cite{clad3} and \cite{nht2}.

Let $P$ be a poset. A nonempty subset $D$ of $P$ is \emph{directed} if every two elements of $D$ have an upper bound in $D$. If $D$ is also a lower set ($D=\downarrow D=\{x\in P: x\le d \mbox{ for some } d\in D\}$), then $D$ is called an \emph{ideal}. A poset is called a \emph{directed complete poset} (dcpo, for short) if its every directed subset has a  supremum. A \emph{complete lattice} is a poset in which every subset has a supremum and an infimum. A subset $U$ of a poset $P$ is \emph{Scott open} if (i) it is an upper set ($U=\ua U=\{x\in P: u\le x \mbox{ for some } u\in U\}$) and (ii) for every directed subset $D$ of $P$ with $\sup D$ existing and $\sup D\in U$, it follows that $D\cap U\not=\emptyset$. The complements of Scott open sets are called \emph{Scott closed} sets. The collection of all Scott open subsets of $P$ form a topology on $P$, called the \emph{Scott topology} of $P$ and is denoted by $\sigma(P)$. The collection of all Scott closed subsets of $P$ is denoted by $\Gamma (P)$. The space $(P,\sigma(P))$ called the Scott space of $P$ is written as $\Sigma P$.

A subset $K$ of a topological space $X$ is \emph{compact} if every open cover of $K$ has a finite subcover. A set $K$ of a topological space is called \emph{saturated} if it is the intersection of its open neighborhood ($K=\ua K$ in its specialization order).  The \emph{saturation} sat$A$ of a set $A$ is the intersection of all its open neighborhood.

\begin{definition}(\cite{clad3})
 (1) A topological space $X$ is \emph{sober} if it is $T_{0}$ and every irreducible closed subset of $X$ is the closure of a (unique) point.

 (2) We shall say that a $T_{0}$ space $X$ is \emph{well-filtered} if for each filter basis $\mathcal{C}$ of compact saturated sets and each open set $U$ with $\bigcap \mathcal{C} \subseteq U$, there is a $K\in \mathcal{C}$ with $K\subseteq U$.

 (3) A space $X$ is \emph{coherent} if the intersection of any two compact saturated sets is again compact.
\end{definition}
\begin{definition}(\cite{nht2})
  An ideal $I$ of a poset $P$ is non-trivial if $\max I\notin I$. We use $Id(P)$ to denote  the set of all non-trivial ideals of a poset $P$.
\end{definition}
\begin{corollary}(\cite{clad3})\label{b}
If $L$ is a sup semilattice such that the sup operation is jointly Scott-continuous, then $\Sigma L$ is  sober.
\end{corollary}

\section{Countable ideals}

The following  lemma is  critical for our  later discussions.
\begin{lemma}\label{c}
 Let $P,Q$ be two posets. If ${\mid}Id(P){\mid},{\mid}Id(Q){\mid}$ are both countable, then $\Sigma(P\times Q)=\Sigma P\times \Sigma Q$.
 \begin{proof}
 Obviously, $\sigma(P)\times \sigma(Q)\subseteq \sigma(P\times Q)$. It remains to prove that $\sigma(P\times Q)\subseteq \sigma(P)\times \sigma(Q)$. Let $U$ be a nonempty Scott open set  and $(a_{1},b_{1})\in U$. We denote $Id(P)$ and $Id(Q)$ by $\{I_{n}^{P}\mid n\in \mathbb{N}\}$ and $\{I_{n}^{Q}\mid n\in \mathbb{N}\}$, respectively.

 For $n=1$, $A_{1}=\{a_{1}\},B_{1}=\{b_{1}\}$.

 For $n=2$, we  define $A_{2}$ and $B_{2}$ below:

  If $\sup I_{1}^{P}\in \ua A_{1}$, then $(\sup I_{1}^{P},b_{1})\in U$. It follows that there exists $d_{1}^{P}\in I_{1}^{P}$ such that $(d_{1}^{P},b_{1})\in U$ by the Scott openness of $U$. Let $A_{2}=\{d_{1}^{P}\}$ in this case, and  $A_{2}=\emptyset$ otherwise. Note that $(A_{1}\cup A_{2})\times B_{1}\subseteq U$.

  If $\sup I_{1}^{Q}\in \ua B_{1}$, we have $(A_{1}\cup A_{2})\times \{\sup I_{1}^{Q}\}\subseteq U$.
   For each $a\in A_{1}\cup A_{2}$, we can choose a $d_a$ form $I_1^Q$ satisftying $(a,d_a)\in U$ by the Scott openness of $U$.
   Since $A_1\cup A_2$ is finite and $I_1^Q$ is directed, there exists $d_{1}^{Q}\in I_{1}^{Q}$ such that $(A_{1}\cup A_{2})\times \{d_{1}^{Q}\}\subseteq U$. Let $B_{2}=\{d_{1}^{Q}\}$ in this case, and  $B_{2}=\emptyset$ otherwise. We conclude that $(A_{1}\cup A_{2})\times (B_{1}\cup B_{2})\subseteq U$.

  For $n=3$, we first consider the two index sets:
  \begin{align*}
    E_{1} & =\Big\{i\in\{1\}\mid \sup I_{i}^{P}\notin \ua A_{1} ~ {\rm and}~ \sup I_{i}^{P}\in \ua A_{2}\Big\}\cup\Big\{i\in \{2\}\mid \sup I_{i}^{P}\in \ua (A_{1}\cup A_{2})\Big\}~ {\rm and} \\
    F_{1} & =\Big\{i\in\{1\}\mid \sup I_{i}^{Q} \notin \ua B_{1} ~ {\rm and}~ \sup I_{i}^{Q}\in \ua B_{2}\Big\}\cup\Big\{i\in \{2\}\mid \sup I_{i}^{Q}\in \ua (B_{1}\cup B_{2})\Big\}.
  \end{align*}
However, if $\sup I_{1}^{P} \notin \ua A_{1}$ and $\sup I_{1}^{Q} \notin \ua B_{1}$, then $A_2=\emptyset$ and $B_2=\emptyset$ from the above step.
In this way, $\Big\{i\in\{1\}\mid \sup I_{i}^{P} \notin \ua A_{1} $ and $ \sup I_{i}^{P}\in \ua A_{2}\Big\}$ and $\Big\{i\in\{1\}\mid \sup I_{i}^{Q} \notin \ua B_{1} $ and $\sup I_{i}^{Q}\in \ua B_{2}\Big\}$ must be empty.
Next, we define $A_{3}$ and $B_{3}$ in the similar way as before.

 If $E_1\neq \emptyset$, then $E_1=\{2\}$. We have $\{\sup I_{2}^{P}\}\times (B_{1}\cup B_{2})\subseteq U$. Through the similar discussion process, we can deduce that  there exists $d_{2}^{P}\in I_{2}^{P}$ such that $\{d_{2}^{P}\}\times (B_{1}\cup B_{2})\subseteq U$ because $B_{1}\cup B_{2}$ is  finite and $I_{2}^{P}$ is directed. Let $A_{3}=\{d_{2}^{P}\}$ in this case, and $A_{3}=\emptyset$ otherwise. Note that $(A_{1}\cup A_{2}\cup A_{3})\times (B_{1}\cup B_{2})\subseteq U$.

 If $F_{1}\neq \emptyset$, then $F_{1}=\{2\}$. Thus $(A_{1}\cup A_{2}\cup A_{3})\times \{\sup I_{2}^{Q}\}\subseteq U$. Note that $A_{1}\cup A_{2}\cup A_{3}$ is a finite set. It follows that there exists $d_{2}^{Q}\in I_{2}^{Q}$ such that $(A_{1}\cup A_{2}\cup A_{3})\times \{d_{2}^{Q}\}\subseteq U$. Let $B_{3}=\{d_{2}^{Q}\}$ in this case, and $B_{3}=\emptyset$ otherwise. We conclude that $(A_{1}\cup A_{2}\cup A_{3})\times (B_{1}\cup B_{2}\cup B_{3})\subseteq U$.

 For $n=4$, we also consider the two index sets:
  \begin{align*}
    E_{2} & =\Big\{i\in\{1,2\}\mid \sup I_{i}^{P}\notin \bigcup_{k=1}^{i}\ua A_{k} ~ {\rm and}~ \sup I_{i}^{P}\in \bigcup_{k=i+1}^{3}\ua A_{k}\Big\}\cup\Big\{i\in \{3\}\mid \sup I_{i}^{P}\in \bigcup_{k=1}^{3}\ua A_{k}\Big\}, \\
    F_{2} & =\Big\{i\in\{1,2\}\mid \sup I_{i}^{Q} \notin \bigcup_{k=1}^{i}\ua B_{k}~ {\rm and}~ \sup I_{i}^{Q}\in \bigcup_{k=i+1}^{3}\ua B_{k}\Big\}\cup\Big\{i\in \{3\}\mid \sup I_{i}^{Q}\in \bigcup_{k=1}^{3}\ua B_{k}\Big\}.
  \end{align*}
Next, we define $A_{4}$ and $B_{4}$ in the following:

 If $E_2\neq \emptyset$, then $i\in \{1,2\}$ implies $\sup I_{i}^{P}\in \bigcup_{k=i+1}^{3}\ua A_{k}\subseteq \bigcup_{k=1}^{3}\ua A_{k}$, and $i=3$ implies $\sup I_{i}^{P}\in \bigcup_{k=1}^{3}\ua A_{k}$.
 Thus $\sup I_{i}^{P}\in \bigcup_{k=1}^{3}\ua A_{k}$ for all $i\in E_2$.
 So for each $i\in E_2$, $\{\sup I_{i}^{P}\}\times (B_{1}\cup B_{2}\cup B_3)\subseteq U$ implies that there exists $d_{i}^{P}\in I_{i}^{P}$ such that $\{d_{i}^{P}\}\times (B_{1}\cup B_{2}\cup B_3)\subseteq U$ because $B_{1}\cup B_{2}\cup B_3$ is  finite and $I_{i}^{P}$ is directed. Let $A_{4}=\{d_{i}^{P}\mid i\in E_2\}$ in this case, and $A_{4}=\emptyset$ otherwise. Note that $$\Big(\bigcup_{k=1}^{4}\ua A_k\Big)\times \Big(\bigcup_{k=1}^{3}\ua B_k\Big)\subseteq U.$$

 If $F_{2}\neq \emptyset$, then $i\in \{1,2\}$ implies $\sup I_{i}^{Q}\in \bigcup_{k=i+1}^{3}\ua B_{k}\subseteq \bigcup_{k=1}^{3}\ua B_{k}$, and $i=3$ implies $\sup I_{i}^{Q}\in \bigcup_{k=1}^{3}\ua B_{k}$.
 Thus $\sup I_{i}^{Q}\in \bigcup_{k=1}^{3}\ua B_{k}$ for all $i\in F_2$.
 So for each $i\in F_2$, $(\bigcup_{k=1}^{4}A_{k})\times \{\sup I_{i}^{Q}\}\subseteq U$ implies that there exists $d_{i}^{Q}\in I_{i}^{Q}$ such that $(\bigcup_{k=1}^{4}A_{k})\times \{d_{i}^{Q}\}\subseteq U$ since $\bigcup_{k=1}^{4}A_{k}$ is a finite set and $d_{i}^{Q}$ is directed.
  Let $B_{4}=\{d_{i}^{Q}\mid i\in F_2\}$ in this case, and $B_{4}=\emptyset$ otherwise. We conclude that $$\Big(\bigcup_{k=1}^{4}A_{k}\Big)\times \Big(\bigcup_{k=1}^{4}B_{k}\Big)\subseteq U.$$

  For $n\geq 4$,  we assume that $$\Big(\bigcup_{k=1}^{n-1}A_k \Big)\times \Big(\bigcup_{k=1}^{n-1} B_k  \Big)\subseteq U.$$
  Then we define $A_{n}$ and $B_{n}$ inductively.

We  first consider the following two index sets:
\begin{align*}
  E_{n-2}  = &\Big\{i\in\{1,\ldots, n-2\}\mid \sup I_{i}^{P} \notin \bigcup_{k=1}^{i}\ua A_{k} ~{\rm and}~\sup I_{i}^{P}\in \bigcup_{k=i+1}^{n-1}\ua A_{k}\Big\}\\
  &\cup\Big\{i\in \{n-1\}\mid \sup I_{i}^{P}\in \bigcup_{k=1}^{n-1}\ua A_{k}\Big\},\\
  F_{n-2}  =& \Big\{i \in \{1,\ldots, n-2\}\mid \sup I_{i}^{Q} \notin \bigcup_{k=1}^{i} \ua B_{k}~{\rm and}~\sup I_{i}^{Q}\in \bigcup_{k=i+1}^{n-1} \ua B_{k}\Big\}\\
  &\cup\Big\{i\in \{n-1\}\mid \sup I_{i}^{Q}\in \bigcup_{k=1}^{n-1} \ua B_{k}\Big\}.
\end{align*}
Note that $\Big\{i\in\{1,\ldots, n-2\}\mid \sup I_{i}^{P} \notin \bigcup_{k=1}^{i}\ua A_{k} ~{\rm and}~\sup I_{i}^{P}\in \bigcup_{k=i+1}^{n-1}\ua A_{k}\Big\}$ and $\Big\{i \in \{1,\ldots, n-2\}\mid \sup I_{i}^{Q} \notin \bigcup_{k=1}^{i} \ua B_{k}~{\rm and}~\sup I_{i}^{Q}\in \bigcup_{k=i+1}^{n-1} \ua B_{k}\Big\}$ may not be empty.

If $E_{n-2}\neq \emptyset$, similarly, we can deduce \{$\sup I_{i}^{P}\}\times \big(\bigcup_{k=1}^{n-1} B_{k}\big)\subseteq U$ for any $i\in E_{n-1}$. Note that $\bigcup_{k=1}^{n-1} B_{k}$ is a finite set and each $I_{i}^{P}$ is directed. Thus there exists $d_{i}^{P}\in I_{i}^{P}$ such that $\{d_{i}^{P}\}\times \big(\bigcup_{k=1}^{n-1} B_{k}\big)\subseteq U$ for any $i\in E_{n-2}$.
Let $A_{n}=\{d_{i}^{P}{\mid} i\in E_{n-1}\}$ in this case, and  $A_{n}=\emptyset$ otherwise. It follows that $$\Big(\bigcup_{k=1}^{n} A_{k}\Big)\times \Big(\bigcup_{k=1}^{n-1} B_{k}\Big)\subseteq U.$$

If $F_{n-2}\neq \emptyset$, then $\big(\bigcup_{k=1}^{n} A_{k}\big)\times \{\sup I_{i}^{Q}\}\subseteq U$ for any $i\in F_{n-2}$. Note that $\bigcup_{k=1}^{n} A_{k}$ is a finite set. This means that there exists $d_{i}^{Q}\in I_{i}^{Q}$ such that $\big(\bigcup_{k=1}^{n} A_{k}\big)\times \{d_{i}^{Q}\}\subseteq U$ for any $i\in F_{n-2}$. Let $B_{n}=\{d_{i}^{Q}\mid i\in F_{n-2}\}$ in this case, and $B_{n}=\emptyset$ otherwise. We conclude that $$\Big(\bigcup_{k=1}^{n} A_{k}\Big)\times \Big(\bigcup_{k=1}^{n} B_{k}\Big)\subseteq U.$$

  Let $A=\bigcup_{n\in \mathbb{N}}A_{n}$ and $B=\bigcup_{n\in \mathbb{N}}B_{n}$. It is easy to see that $(a_{1},b_{1})\in A_1\times B_1 \subseteq\ua A\times \ua B \subseteq U$. It suffices to prove that $\ua A,\ua B$ are both Scott open.

  Let $D$ be a directed subset of $P$. If $\sup D \in D$, then $D\cap \ua A\neq \emptyset$. If $\sup D \notin D$, i.e., $D$ contains no maximal element,
  then $\da D \in Id(P)$. Thus, there exists $n_0\in \mathbb{N}$ such that $\da D=I_{n_{0}}^{P}$.
  Therefore, $\sup D\in \ua A$ can imply that $\sup I_{n_{0}}^{P}\in \ua A$. Let $n_{1}=\inf\{n\in \mathbb{N}\mid \sup I_{n_{0}}^{P}\in \ua A_{n}\}$. Then $\sup I_{n_{0}}^{P}\in \ua A_{n_{1}}$. Now we need to distinguish between the following two cases for $n_{0},n_{1}$.

  Case 1, $n_{0}<n_{1}$. If $n_0=1,n_1=2$, then $\sup I_1^P \notin \ua A_1$ implies $A_2= \emptyset$, which contradicts the condition $\sup I_1^P \in \ua A_2$. So $n_1 \geqslant 3$. The fact that $\sup I_{n_{0}}^{P}\notin \bigcup_{k=1}^{n_{0}}\ua A_{k}$ and $\sup I_{n_{0}}^{P}\in \bigcup_{k=n_{0}+1}^{n_{1}}\ua A_{k}$ can imply $n_0 \in E_{n_1-1}$. This means that $I_{n_{0}}^{P}\cap A_{n_{1}+1}\neq \emptyset$. Hence, $D\cap \ua A\neq \emptyset$.

  Case 2, $n_{0}\geq n_{1}$.
  If $n_0=n_1=1$, then $\sup I_1^P \in \ua A_1$ implies $I_1^P \cap A_2 \neq \emptyset$. If $n_0 \geqslant 2$,   then $\sup I_{n_{0}}^{P}\in \ua A_{n_1}\subseteq \bigcup_{k=1}^{n_{0}}\ua A_{k}$, which implies $n_0 \in E_{n_0-1}$. It follows that $I_{n_{0}}^{P}\cap A_{n_{0}+1}\neq \emptyset$. Therefore, $D\cap \ua A\neq \emptyset$.

  Hence, $\ua A$ is Scott open, and  $\ua B$ is Scott open by the similar proof.
 \end{proof}
\end{lemma}
The following example reveals that the above lemma does not hold on the contrary.
\begin{example}
  Let $L=\mathbb{R}\times \mathbb{N}$, where $\mathbb{R}$ denotes the set of all real numbers and  $\mathbb{N}$  all natural numbers. We define an order $\leq$ on $L$ as follows:

  $(r,m)\leq (s,n)$ if and only if  $r=s$ and $m\leq n$.

Obviously, $L$ is continuous. Then $\sigma(L\times L)=\sigma(L) \times \sigma(L)$. But it is easy to see that ${\mid}Id(L){\mid}$ is uncountable.
\end{example}

By the above lemma, we can get the following corollary.
\begin{corollary}\label{g}
  Let $L$ be a dcpo with ${\mid}Id(L){\mid}$ countable. If $\Sigma L$ is coherent and well-filtered, then $\Sigma L$ is sober.
  \begin{proof}
    Let $A$ be an irreducible closed subset of $\Sigma L$. It suffices to prove that $A$ is directed, which means that $\ua x\cap \ua y\cap A\neq \emptyset$ for any $x,y\in A$.

   Write $B=\{(a,b)\in L\times L\mid \ua a\cap \ua b\subseteq L\backslash A\}$.
   We claim that $B$ is Scott open in $L\times L$.
   Obviously, $B$ is an upper set.
   Let $(x_{i},y_{i})_{i\in I}$ be a directed subset of $L\times L$ with $\sup_{i\in I}(x_{i},y_{i})\in B$.
   Then $(\sup_{i\in I}x_{i},\sup_{i\in I}y_{i})\in B$, which is equivalent to saying that $\ua \sup_{i\in I}x_{i}\cap \ua \sup_{i\in I}y_{i}\subseteq L\backslash A$. It follows that $\bigcap_{i\in I} (\ua x_{i}\cap \ua y_{i})\subseteq L\backslash A$.
   Since $\Sigma L$ is coherent and well-filtered, we can find some index  $i\in I$ such that $\ua x_{i}\cap \ua y_{i}\subseteq L\backslash A$.
   This implies that $(x_{i},y_{i})\in B$. Thus $B$ is Scott open.

   It is worth noting that ${\mid}Id(L){\mid}$ is countable.
    $\Sigma(L\times L)=\Sigma L\times \Sigma L$ from Lemma \ref{c}.
    For the sake of contradiction, we assume that there is $x,y\in A$ such that  $\ua x\cap \ua y\cap A=\emptyset$.
   The fact that $(x,y)\in B\subseteq \sigma(L\times L)$ implies that we can find $U_{x},U_{y}\in \sigma (L)$ such that $(x,y)\in U_{x}\times  U_{y}\subseteq B$.
   Note that $x\in U_{x}\cap A$ and $y\in U_{y}\cap A$.
   By the irreducibility of $A$, we have  $A\cap U_{x}\cap U_{y}\neq \emptyset$.
   Pick $a\in A\cap U_{x}\cap U_{y}$. Then $(a,a)\in U_{x}\times U_{y}\subseteq B$, that is, $a\in \ua a \cap \ua a\subseteq L\backslash A$. It contradicts  the assumption that $a\in A$.
   Hence, $A$ is directed and $\sup A\in A$. So $A=\da \sup A$.
  \end{proof}
\end{corollary}
\begin{example}(\cite{meet5})
  Jia constructs a dcpo $P=\mathbb{N}\times\mathbb{N}\times(\mathbb{N}\cup \{\infty\})$.
  The order $\leqslant$ on $P$ is defined as follows:

  $(i_1,j_1,m_1)\leqslant (i_2,j_2,m_2)$ if and only if:

  $\bullet$ $i_1=i_2,j_1=j_2, m_1\leqslant m_2\leqslant \infty$;

  $\bullet$ $i_2=i_1+1,m_1\leqslant j_2$, $m_2=\infty$.

  \begin{tikzpicture}[line width=0.5pt,scale=1.1]
\fill[black] (0,0) circle (1.5pt);
\fill[black] (0,1) circle (1.5pt);
\fill[black] (0,2) circle (1.5pt);
\fill[black] (0,4) circle (1.5pt);
\draw (0,0)--(0,2);
\draw [dashed](0,2)--(0,4);
\fill[black] (1,0) circle (1.5pt);
\fill[black] (1,1) circle (1.5pt);
\fill[black] (1,2) circle (1.5pt);
\fill[black] (1,4) circle (1.5pt);
\draw (1,0)--(1,2);
\draw [dashed](1,2)--(1,4);
\fill[black] (2,0) circle (1.5pt);
\fill[black] (2,1) circle (1.5pt);
\fill[black] (2,2) circle (1.5pt);
\fill[black] (2,4) circle (1.5pt);
\draw (2,0)--(2,2);
\draw [dashed](2,2)--(2,4);
\draw [dashed](-1,4)--(0,4);
\draw (0,0.4)--(2,0.4);
\draw (0,-0.4)--(2,-0.4);
\draw (2,-0.4) arc (-90:90:0.4);
\draw [dashed](-1,0.4)--(0,0.4);
\draw [dashed](-1,-0.4)--(0,-0.4);

\draw (0,0.6)--(2,0.6);
\draw (0,1.4)--(2,1.4);
\draw (2,0.6) arc (-90:90:0.4);
\draw [dashed](-1,1.4)--(0,1.4);
\draw [dashed](-1,0.6)--(0,0.6);

\draw (0,1.6)--(2,1.6);
\draw (0,2.4)--(2,2.4);
\draw (2,1.6) arc (-90:90:0.4);
\draw [dashed](-1,2.4)--(0,2.4);
\draw [dashed](-1,1.6)--(0,1.6);


\fill[black] (6,0) circle (1.5pt);
\fill[black] (6,1) circle (1.5pt);
\fill[black] (6,2) circle (1.5pt);
\fill[black] (6,4) circle (1.5pt);
\draw (6,0)--(6,2);
\draw [dashed](6,2)--(6,4);
\fill[black] (4,0) circle (1.5pt);
\fill[black] (4,1) circle (1.5pt);
\fill[black] (4,2) circle (1.5pt);
\fill[black] (4,4) circle (1.5pt);
\draw (4,0)--(4,2);
\draw [dashed](4,2)--(4,4);
\fill[black] (5,0) circle (1.5pt);
\fill[black] (5,1) circle (1.5pt);
\fill[black] (5,2) circle (1.5pt);
\fill[black] (5,4) circle (1.5pt);
\draw (5,0)--(5,2);
\draw [dashed](5,2)--(5,4);
\draw [dashed](4,4)--(3,4);
\draw (4,0.4)--(6,0.4);
\draw (4,-0.4)--(6,-0.4);
\draw (6,-0.4) arc (-90:90:0.4);
\draw [dashed](4,0.4)--(3,0.4);
\draw [dashed](4,-0.4)--(3,-0.4);

\draw (4,0.6)--(6,0.6);
\draw (4,1.4)--(6,1.4);
\draw (6,0.6) arc (-90:90:0.4);
\draw [dashed](4,1.4)--(3,1.4);
\draw [dashed](4,0.6)--(3,0.6);

\draw (4,1.6)--(6,1.6);
\draw (4,2.4)--(6,2.4);
\draw (6,1.6) arc (-90:90:0.4);
\draw [dashed](4,2.4)--(3,2.4);
\draw [dashed](4,1.6)--(3,1.6);


\fill[black] (8,0) circle (1.5pt);
\fill[black] (8,1) circle (1.5pt);
\fill[black] (8,2) circle (1.5pt);
\fill[black] (8,4) circle (1.5pt);
\draw (8,0)--(8,2);
\draw [dashed](8,2)--(8,4);
\fill[black] (9,0) circle (1.5pt);
\fill[black] (9,1) circle (1.5pt);
\fill[black] (9,2) circle (1.5pt);
\fill[black] (9,4) circle (1.5pt);
\draw (9,0)--(9,2);
\draw [dashed](9,2)--(9,4);
\fill[black] (10,0) circle (1.5pt);
\fill[black] (10,1) circle (1.5pt);
\fill[black] (10,2) circle (1.5pt);
\fill[black] (10,4) circle (1.5pt);
\draw (10,0)--(10,2);
\draw [dashed](10,2)--(10,4);
\draw [dashed](8,4)--(7,4);
\draw (8,0.4)--(10,0.4);
\draw (8,-0.4)--(10,-0.4);
\draw (10,-0.4) arc (-90:90:0.4);
\draw [dashed](8,0.4)--(7,0.4);
\draw [dashed](8,-0.4)--(7,-0.4);

\draw (8,0.6)--(10,0.6);
\draw (8,1.4)--(10,1.4);
\draw (10,0.6) arc (-90:90:0.4);
\draw [dashed](8,1.4)--(7,1.4);
\draw [dashed](8,0.6)--(7,0.6);

\draw (8,1.6)--(10,1.6);
\draw (8,2.4)--(10,2.4);
\draw (10,1.6) arc (-90:90:0.4);
\draw [dashed](8,2.4)--(7,2.4);
\draw [dashed](8,1.6)--(7,1.6);

\draw (2.4,0)--(6,4);
\draw (2.4,1)--(5,4);
\draw (2.4,2)--(4,4);
\draw (6.4,0)--(10,4);
\draw (6.4,1)--(9,4);
\draw (6.4,2)--(8,4);

\draw [dashed](11,0)--(12,0);
\draw [dashed](11,4)--(12,4);
\draw [dashed](11,3)--(12,3);

\node[font=\footnotesize] (1) at(1,-0.25) {$(1,2,1)$};
\node[font=\footnotesize] (1) at(5,-0.25) {$(2,2,1)$};
\node[font=\footnotesize] (2) at(6,4.25) {$(2,1,\infty)$};
\node[font=\footnotesize] (3) at(10,4.25) {$(3,1,\infty)$};
\end{tikzpicture}

  In \cite{meet5}, it shows that $\Sigma P$ is well-filtered and not sober. We also find two facts:
  \begin{itemize}
    \item $Id(P)$ is countable.
    \item $P$ is not coherent
  \end{itemize}

Let $D$ be a directed subset with no maximal element. It is easy to verify that $D$ is infinite and is contained in $\da (i,j,\infty)$ for some $i,j\in \mathbb{N}$ with its supremum being $(i,j,\infty)$.
If $(i_1,j_1,m_1)$, $(i_2,j_2,m_2)\in D$ with $i_1\neq i_2$, then $D$ must have a greatest element, which contradicts the hypothesis on $D$.
After the similar discussion, we have $D \subseteq \{(i,j,m)\mid m\in \mathbb{N}\}$ for fixed $i,j \in\mathbb{N}$.
Thus, $Id(P)=\big\{\{i\}\times\{j\}\times \mathbb{N}\mid i,j \in\mathbb{N}\big\}$, which is countable, obviously.

As for the coherence, we only need to find that the intersection of two principle filters is not compact.
We claim that $\ua (1,2,1)\cap \ua (1,3,1)=\{(2,j,\infty)\mid j\in\mathbb{N}\}$ is not compact.
Let $C_j=\bigcup_{k\geqslant j}\da (2,k,\infty)\cup \{(1,n,\infty)\mid n\in \mathbb{N}\}$.
Obviously, $\{C_j\mid j\in\mathbb{N}\}$ is a filtered family of Scott closed subsets and $\{(2,j,\infty)\mid j\in\mathbb{N}\}$ meets all $C_j$. But the intersection $\bigcap_{j\in \mathbb{N}}C_j\cap \{(2,j,\infty)\mid j\in\mathbb{N}\}=\emptyset$. So $\ua (1,2,1)\cap \ua (1,3,1)$ is not compact.

This example indicates that the condition in the above corollary is essential.
\end{example}
The following theorem  gives a partial answer to Problem \ref{a} based on the above corollary.
\begin{theorem}
  Let $L$ be a complete lattice. If ${\mid}Id(L){\mid}$ is countable, then $\Sigma L$ is sober.
  \begin{proof}
    From \cite{xiaodong2} and \cite{Xi03}, we deduce that $\Sigma L$ is well-filtered and coherent. The result is evident by Corollary \ref{g}.
  \end{proof}
\end{theorem}

\section{A countable non-sober complete lattice}
In this section, we present a counterexample in order to solve the open problem mentioned in the introduction posed by Jung.
\begin{example}
  Let $L=\mathbb{N}\cup \mathbb{N}^{<\mathbb{N}}\cup \{\top\}$, where $\mathbb{N}$ is the set of positive natural numbers and $\mathbb{N}^{<\mathbb{N}}$ the set of non-empty finite sequences of natural numbers. We define an order $\leq $ on $L$ as follows:

  $x\leq y$ if and only if:

  $\bullet$ $x\leq y$ in $\mathbb{N}$;

  $\bullet$ $x,y\in \mathbb{N}^{<\mathbb{N}}$, $y=x.t$, $t\in \mathbb{N}^{<\mathbb{N}}$ or $y=x$;

  $\bullet$ $x\in L$, $y=\top$.

  Then $L$ can be easily depicted as Fig.\ref{F.1}.
\begin{figure}[H]
 \centering
  \includegraphics[height=5cm]{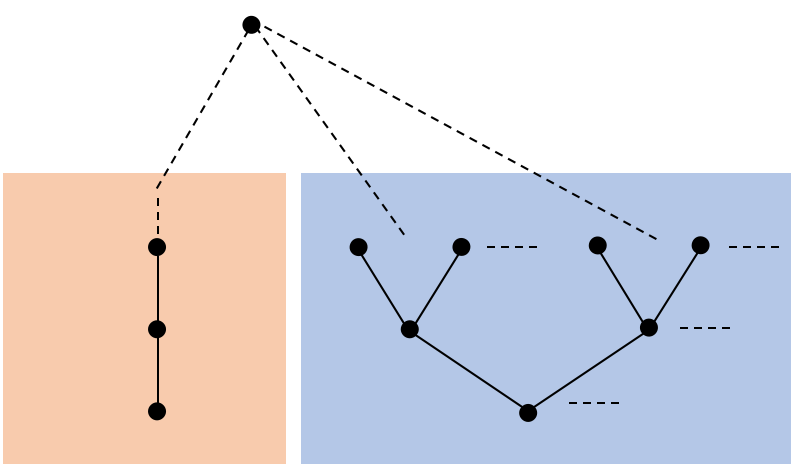}\\
  \caption{The basic gadget of $P$}\label{F.1}
\end{figure}
\end{example}
To construct the final counterexample, we need to give a monotone injection $f_{m,n}:\mathbb{N}^{<\mathbb{N}}\rightarrow \mathbb{N}$ for any $(m,n)\in \mathbb{N}\times \mathbb{N}$ with $m<n$. In the following, we provide a specific construction for them.
\begin{remark}\label{mon}
  By induction way, there exists a monotone injection $f:\mathbb{N}^{<\mathbb{N}}\rightarrow \mathbb{N}$.

\end{remark}

Let $P=\mathbb{N}\times L$. Fix a map $i:\mathbb{N}\times \mathbb{N}\rightarrow \mathcal{P}(\mathbb{N})$, where $\mathcal{P}(\mathbb{N})$ is the powerset of $P$, satisfying the following properties:

   $\bullet$ $i(m,n)$ is an infinite subset of $\mathbb{N}$ for each $(m,n)\in \mathbb{N}\times \mathbb{N}$ with $m<n$.

   $\bullet$ $n$ is a strict lower bound of $i(m,n)$ for each $(m,n)\in \mathbb{N}\times \mathbb{N}$ with $m<n$.

   $\bullet$ $i(m_{1},n_{1})\cap i(m_{2},n_{2})=\emptyset$ for any two distinct elements $(m_{1},n_{1}),(m_{2},n_{2})$ on $\mathbb{N}\times \mathbb{N}$ with $m_{1}<n_{1},m_{2}<n_{2}$.

By Remark \ref{mon}, we can fixed a monotone bijection $f_{m,n}:\mathbb{N}^{<\mathbb{N}} \rightarrow i(m,n)$ for any $(m,n)\in \mathbb{N}\times \mathbb{N}$ with $m<n$.
  Let $L_{n}=\{(n,x)\in P\mid x\in L\}$. In this section, $s\in \mathbb{N}^{<\mathbb{N}}$ with ${\mid}s{\mid}=1$ sometimes is considered as a natural number $s$. We define the following relations on $P$:

$\bullet$ $(n,x)<_{1}(m,y)$ if $n=m$, $x< y$ in $L_{n}$;

$\bullet$ $(n,x)<_{2}(m,y)$ if $y=\top$, $x\in \mathbb{N}^{<\mathbb{N}}$ and there exists $k\in \mathbb{N}$ with $k>n$ such that $m\in i(n,k)$ and $m=f_{n,k}(x)$.

$\bullet$ $(n,x)<_{3}(m,y)$ if $y=\top$, $x\in \mathbb{N}$ and there exists $d\in \mathbb{N}$ with $d<n$ such that $m\in i(d,n)$ and $m=f_{d,n}(x)$.

$\bullet$ $(n,x)<_{4}(m,y)$ if $y=\top$, $x\in \mathbb{N}$ and there exists $a,b\in \mathbb{N}$, $s\in\mathbb{N}$ with $a<b$ such that $f_{a,b}(s)=n$ and $f_{a,b}(s.x)=m$.

Then $<:=<_{1}\cup<_{2}\cup<_{3}\cup <_{4}\cup<_{1};<_{2}\cup<_{1};<_{3}\cup <_{1};<_{4}$
 (we use $;$ for relation composition) is transitive and irreflexive. So $\leq := (< \cup =)$ is an order relation.

If $n=m$, then the strict order is depicted as in Fig.\ref{F.1}. Otherwise, $n\neq m$, then the other strict order of $P$ can be easily depicted as in the following figures, respectively.
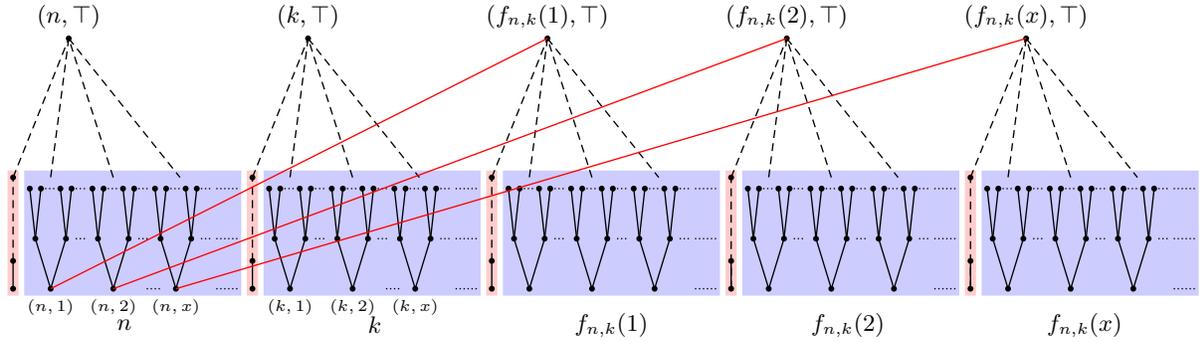
\begin{figure}[H]
  \centering
  \begin{tikzpicture} [line width=0.5pt,scale=0.37]
\fill[red,opacity=0.2] (-0.2,-0.25) rectangle (0.2,4.25)
[xshift=8.6cm] (-0.2,-0.25) rectangle (0.2,4.25)
[xshift=8.6cm] (-0.2,-0.25) rectangle (0.2,4.25)
[xshift=8.6cm] (-0.2,-0.25) rectangle (0.2,4.25)
[xshift=8.6cm] (-0.2,-0.25) rectangle (0.2,4.25);
\fill[blue,opacity=0.2] (0.4,-0.25) rectangle (8.2,4.25)
[xshift=8.6cm] (0.4,-0.25) rectangle (8.2,4.25)
[xshift=8.6cm] (0.4,-0.25) rectangle (8.2,4.25)
[xshift=8.6cm] (0.4,-0.25) rectangle (8.2,4.25)
[xshift=8.6cm] (0.4,-0.25) rectangle (8.2,4.25);
\fill[black] (0,0) circle (3pt)
[xshift=8.6cm] (0,0) circle (3pt)
[xshift=8.6cm] (0,0) circle (3pt)
[xshift=8.6cm] (0,0) circle (3pt)
[xshift=8.6cm] (0,0) circle (3pt);
\fill[black] (0,1) circle (3pt)
[xshift=8.6cm] (0,1) circle (3pt)
[xshift=8.6cm] (0,1) circle (3pt)
[xshift=8.6cm] (0,1) circle (3pt)
[xshift=8.6cm] (0,1) circle (3pt);
\fill[black] (0,4) circle (3pt)
[xshift=8.6cm] (0,4) circle (3pt)
[xshift=8.6cm] (0,4) circle (3pt)
[xshift=8.6cm] (0,4) circle (3pt)
[xshift=8.6cm] (0,4) circle (3pt);
\draw (0,0)--(0,1)
[xshift=8.6cm] (0,0)--(0,1)
[xshift=8.6cm] (0,0)--(0,1)
[xshift=8.6cm] (0,0)--(0,1)
[xshift=8.6cm] (0,0)--(0,1);
\draw [densely dashed](0,1)--(0,4)
[xshift=8.6cm] (0,1)--(0,4)
[xshift=8.6cm] (0,1)--(0,4)
[xshift=8.6cm] (0,0)--(0,4)
[xshift=8.6cm] (0,0)--(0,4);

\fill[black] (0.6,3.6) circle (3pt)
[xshift=1.1cm] (0.6,3.6) circle (3pt)
[xshift=1.15cm] (0.6,3.6) circle (3pt)
[xshift=1.1cm] (0.6,3.6) circle (3pt)
[xshift=1.15cm] (0.6,3.6) circle (3pt)
[xshift=1.1cm] (0.6,3.6) circle (3pt)
[xshift=3cm] (0.6,3.6) circle (3pt)
[xshift=1.1cm] (0.6,3.6) circle (3pt)
[xshift=1.15cm] (0.6,3.6) circle (3pt)
[xshift=1.1cm] (0.6,3.6) circle (3pt)
[xshift=1.15cm] (0.6,3.6) circle (3pt)
[xshift=1.1cm] (0.6,3.6) circle (3pt)
[xshift=3cm] (0.6,3.6) circle (3pt)
[xshift=1.1cm] (0.6,3.6) circle (3pt)
[xshift=1.15cm] (0.6,3.6) circle (3pt)
[xshift=1.1cm] (0.6,3.6) circle (3pt)
[xshift=1.15cm] (0.6,3.6) circle (3pt)
[xshift=1.1cm] (0.6,3.6) circle (3pt)
[xshift=3cm] (0.6,3.6) circle (3pt)
[xshift=1.1cm] (0.6,3.6) circle (3pt)
[xshift=1.15cm] (0.6,3.6) circle (3pt)
[xshift=1.1cm] (0.6,3.6) circle (3pt)
[xshift=1.15cm] (0.6,3.6) circle (3pt)
[xshift=1.1cm] (0.6,3.6) circle (3pt)
[xshift=3cm] (0.6,3.6) circle (3pt)
[xshift=1.1cm] (0.6,3.6) circle (3pt)
[xshift=1.15cm] (0.6,3.6) circle (3pt)
[xshift=1.1cm] (0.6,3.6) circle (3pt)
[xshift=1.15cm] (0.6,3.6) circle (3pt)
[xshift=1.1cm] (0.6,3.6) circle (3pt);

\fill[black] (0.8,1.8) circle (3pt)
[xshift=1.1cm](0.8,1.8) circle (3pt)
[xshift=1.15cm](0.8,1.8) circle (3pt)
[xshift=1.1cm](0.8,1.8) circle (3pt)
[xshift=1.15cm](0.8,1.8) circle (3pt)
[xshift=1.1cm](0.8,1.8) circle (3pt)
[xshift=3cm] (0.8,1.8) circle (3pt)
[xshift=1.1cm] (0.8,1.8) circle (3pt)
[xshift=1.15cm](0.8,1.8) circle (3pt)
[xshift=1.1cm](0.8,1.8) circle (3pt)
[xshift=1.15cm] (0.8,1.8) circle (3pt)
[xshift=1.1cm] (0.8,1.8) circle (3pt)
[xshift=3cm] (0.8,1.8) circle (3pt)
[xshift=1.1cm] (0.8,1.8) circle (3pt)
[xshift=1.15cm](0.8,1.8) circle (3pt)
[xshift=1.1cm](0.8,1.8) circle (3pt)
[xshift=1.15cm] (0.8,1.8) circle (3pt)
[xshift=1.1cm] (0.8,1.8) circle (3pt)
[xshift=3cm] (0.8,1.8) circle (3pt)
[xshift=1.1cm] (0.8,1.8) circle (3pt)
[xshift=1.15cm](0.8,1.8) circle (3pt)
[xshift=1.1cm](0.8,1.8) circle (3pt)
[xshift=1.15cm] (0.8,1.8) circle (3pt)
[xshift=1.1cm] (0.8,1.8) circle (3pt)
[xshift=3cm] (0.8,1.8) circle (3pt)
[xshift=1.1cm] (0.8,1.8) circle (3pt)
[xshift=1.15cm](0.8,1.8) circle (3pt)
[xshift=1.1cm](0.8,1.8) circle (3pt)
[xshift=1.15cm] (0.8,1.8) circle (3pt)
[xshift=1.1cm] (0.8,1.8) circle (3pt);

\fill[black] (1,3.6) circle (3pt)
[xshift=1.1cm](1,3.6) circle (3pt)
[xshift=1.15cm](1,3.6) circle (3pt)
[xshift=1.1cm](1,3.6) circle (3pt)
[xshift=1.15cm](1,3.6) circle (3pt)
[xshift=1.1cm] (1,3.6) circle (3pt)
[xshift=3cm](1,3.6) circle (3pt)
[xshift=1.1cm](1,3.6) circle (3pt)
[xshift=1.15cm](1,3.6) circle (3pt)
[xshift=1.1cm](1,3.6) circle (3pt)
[xshift=1.15cm] (1,3.6) circle (3pt)
[xshift=1.1cm] (1,3.6) circle (3pt)
[xshift=3cm](1,3.6) circle (3pt)
[xshift=1.1cm](1,3.6) circle (3pt)
[xshift=1.15cm](1,3.6) circle (3pt)
[xshift=1.1cm](1,3.6) circle (3pt)
[xshift=1.15cm] (1,3.6) circle (3pt)
[xshift=1.1cm] (1,3.6) circle (3pt)
[xshift=3cm](1,3.6) circle (3pt)
[xshift=1.1cm](1,3.6) circle (3pt)
[xshift=1.15cm](1,3.6) circle (3pt)
[xshift=1.1cm](1,3.6) circle (3pt)
[xshift=1.15cm] (1,3.6) circle (3pt)
[xshift=1.1cm] (1,3.6) circle (3pt)
[xshift=3cm](1,3.6) circle (3pt)
[xshift=1.1cm](1,3.6) circle (3pt)
[xshift=1.15cm](1,3.6) circle (3pt)
[xshift=1.1cm](1,3.6) circle (3pt)
[xshift=1.15cm] (1,3.6) circle (3pt)
[xshift=1.1cm] (1,3.6) circle (3pt);
\draw (0.6,3.6)--(0.8,1.8)
[xshift=1.1cm](0.6,3.6)--(0.8,1.8)
[xshift=1.15cm](0.6,3.6)--(0.8,1.8)
[xshift=1.1cm](0.6,3.6)--(0.8,1.8)
[xshift=1.15cm](0.6,3.6)--(0.8,1.8)
[xshift=1.1cm](0.6,3.6)--(0.8,1.8)
[xshift=3cm] (0.6,3.6)--(0.8,1.8)
[xshift=1.1cm](0.6,3.6)--(0.8,1.8)
[xshift=1.15cm](0.6,3.6)--(0.8,1.8)
[xshift=1.1cm](0.6,3.6)--(0.8,1.8)
[xshift=1.15cm] (0.6,3.6)--(0.8,1.8)
[xshift=1.1cm] (0.6,3.6)--(0.8,1.8)
[xshift=3cm] (0.6,3.6)--(0.8,1.8)
[xshift=1.1cm](0.6,3.6)--(0.8,1.8)
[xshift=1.15cm](0.6,3.6)--(0.8,1.8)
[xshift=1.1cm](0.6,3.6)--(0.8,1.8)
[xshift=1.15cm] (0.6,3.6)--(0.8,1.8)
[xshift=1.1cm] (0.6,3.6)--(0.8,1.8)
[xshift=3cm] (0.6,3.6)--(0.8,1.8)
[xshift=1.1cm](0.6,3.6)--(0.8,1.8)
[xshift=1.15cm](0.6,3.6)--(0.8,1.8)
[xshift=1.1cm](0.6,3.6)--(0.8,1.8)
[xshift=1.15cm] (0.6,3.6)--(0.8,1.8)
[xshift=1.1cm] (0.6,3.6)--(0.8,1.8)
[xshift=3cm] (0.6,3.6)--(0.8,1.8)
[xshift=1.1cm](0.6,3.6)--(0.8,1.8)
[xshift=1.15cm](0.6,3.6)--(0.8,1.8)
[xshift=1.1cm](0.6,3.6)--(0.8,1.8)
[xshift=1.15cm] (0.6,3.6)--(0.8,1.8)
[xshift=1.1cm] (0.6,3.6)--(0.8,1.8);

\draw (0.8,1.8)--(1,3.6)
[xshift=1.1cm](0.8,1.8)--(1,3.6)
[xshift=1.15cm](0.8,1.8)--(1,3.6)
[xshift=1.1cm](0.8,1.8)--(1,3.6)
[xshift=1.15cm](0.8,1.8)--(1,3.6)
[xshift=1.1cm](0.8,1.8)--(1,3.6)
[xshift=3cm] (0.8,1.8)--(1,3.6)
[xshift=1.1cm](0.8,1.8)--(1,3.6)
[xshift=1.15cm](0.8,1.8)--(1,3.6)
[xshift=1.1cm](0.8,1.8)--(1,3.6)
[xshift=1.15cm] (0.8,1.8)--(1,3.6)
[xshift=1.1cm] (0.8,1.8)--(1,3.6)
[xshift=3cm] (0.8,1.8)--(1,3.6)
[xshift=1.1cm](0.8,1.8)--(1,3.6)
[xshift=1.15cm](0.8,1.8)--(1,3.6)
[xshift=1.1cm](0.8,1.8)--(1,3.6)
[xshift=1.15cm] (0.8,1.8)--(1,3.6)
[xshift=1.1cm] (0.8,1.8)--(1,3.6)
[xshift=3cm] (0.8,1.8)--(1,3.6)
[xshift=1.1cm](0.8,1.8)--(1,3.6)
[xshift=1.15cm](0.8,1.8)--(1,3.6)
[xshift=1.1cm](0.8,1.8)--(1,3.6)
[xshift=1.15cm] (0.8,1.8)--(1,3.6)
[xshift=1.1cm] (0.8,1.8)--(1,3.6)
[xshift=3cm] (0.8,1.8)--(1,3.6)
[xshift=1.1cm](0.8,1.8)--(1,3.6)
[xshift=1.15cm](0.8,1.8)--(1,3.6)
[xshift=1.1cm](0.8,1.8)--(1,3.6)
[xshift=1.15cm] (0.8,1.8)--(1,3.6)
[xshift=1.1cm] (0.8,1.8)--(1,3.6);
\draw [densely dotted](1.15,3.6)--(1.5,3.6)
[xshift=1.1cm] (1.15,3.6)--(1.5,3.6)
[yshift=-1.8cm] (1.15,3.6)--(1.5,3.6)
[shift={(1.15cm, 1.8cm)}](1.15,3.6)--(1.5,3.6)
[xshift=1.1cm] (1.15,3.6)--(1.5,3.6)
[yshift=-1.8cm] (1.15,3.6)--(1.5,3.6)
[shift={(1.15cm,1.8cm)}] (1.15,3.6)--(1.5,3.6)
[xshift=1.1cm](1.15,3.6)--(1.5,3.6)
[yshift=-1.8cm] (1.15,3.6)--(1.5,3.6)
[shift={(3cm,1.8cm)}] (1.15,3.6)--(1.5,3.6)
[xshift=1.1cm] (1.15,3.6)--(1.5,3.6)
[yshift=-1.8cm] (1.15,3.6)--(1.5,3.6)
[shift={(1.15cm, 1.8cm)}](1.15,3.6)--(1.5,3.6)
[xshift=1.1cm] (1.15,3.6)--(1.5,3.6)
[yshift=-1.8cm] (1.15,3.6)--(1.5,3.6)
[shift={(1.15cm,1.8cm)}] (1.15,3.6)--(1.5,3.6)
[xshift=1.1cm](1.15,3.6)--(1.5,3.6)
[yshift=-1.8cm] (1.15,3.6)--(1.5,3.6)
[shift={(3cm,1.8cm)}] (1.15,3.6)--(1.5,3.6)
[xshift=1.1cm] (1.15,3.6)--(1.5,3.6)
[yshift=-1.8cm] (1.15,3.6)--(1.5,3.6)
[shift={(1.15cm, 1.8cm)}](1.15,3.6)--(1.5,3.6)
[xshift=1.1cm] (1.15,3.6)--(1.5,3.6)
[yshift=-1.8cm] (1.15,3.6)--(1.5,3.6)
[shift={(1.15cm,1.8cm)}] (1.15,3.6)--(1.5,3.6)
[xshift=1.1cm](1.15,3.6)--(1.5,3.6)
[yshift=-1.8cm] (1.15,3.6)--(1.5,3.6)
[shift={(3cm,1.8cm)}] (1.15,3.6)--(1.5,3.6)
[xshift=1.1cm] (1.15,3.6)--(1.5,3.6)
[yshift=-1.8cm] (1.15,3.6)--(1.5,3.6)
[shift={(1.15cm, 1.8cm)}](1.15,3.6)--(1.5,3.6)
[xshift=1.1cm] (1.15,3.6)--(1.5,3.6)
[yshift=-1.8cm] (1.15,3.6)--(1.5,3.6)
[shift={(1.15cm,1.8cm)}] (1.15,3.6)--(1.5,3.6)
[xshift=1.1cm](1.15,3.6)--(1.5,3.6)
[yshift=-1.8cm] (1.15,3.6)--(1.5,3.6)
[shift={(3cm,1.8cm)}] (1.15,3.6)--(1.5,3.6)
[xshift=1.1cm] (1.15,3.6)--(1.5,3.6)
[yshift=-1.8cm] (1.15,3.6)--(1.5,3.6)
[shift={(1.15cm, 1.8cm)}](1.15,3.6)--(1.5,3.6)
[xshift=1.1cm] (1.15,3.6)--(1.5,3.6)
[yshift=-1.8cm] (1.15,3.6)--(1.5,3.6)
[shift={(1.15cm,1.8cm)}] (1.15,3.6)--(1.5,3.6)
[xshift=1.1cm](1.15,3.6)--(1.5,3.6)
[yshift=-1.8cm] (1.15,3.6)--(1.5,3.6);
\fill[black] (1.35,0) circle (3pt)
[xshift=2.25cm](1.35,0) circle (3pt)
[xshift=2.25cm] (1.35,0) circle (3pt)
[xshift=4.1cm](1.35,0) circle (3pt)
[xshift=2.25cm] (1.35,0) circle (3pt)
[xshift=2.25cm](1.35,0) circle (3pt)
[xshift=4.1cm] (1.35,0) circle (3pt)
[xshift=2.25cm] (1.35,0) circle (3pt)
[xshift=2.25cm](1.35,0) circle (3pt)
[xshift=4.1cm](1.35,0) circle (3pt)
[xshift=2.25cm] (1.35,0) circle (3pt)
[xshift=2.25cm](1.35,0) circle (3pt)
[xshift=4.1cm](1.35,0) circle (3pt)
[xshift=2.25cm] (1.35,0) circle (3pt)
[xshift=2.25cm](1.35,0) circle (3pt);

\draw (1.35,0)--(1.9,1.8)
[xshift=2.25cm] (1.35,0)--(1.9,1.8)
[xshift=2.25cm] (1.35,0)--(1.9,1.8)
[xshift=4.1cm] (1.35,0)--(1.9,1.8)
[xshift=2.25cm] (1.35,0)--(1.9,1.8)
[xshift=2.25cm] (1.35,0)--(1.9,1.8)
[xshift=4.1cm] (1.35,0)--(1.9,1.8)
[xshift=2.25cm] (1.35,0)--(1.9,1.8)
[xshift=2.25cm] (1.35,0)--(1.9,1.8)
[xshift=4.1cm] (1.35,0)--(1.9,1.8)
[xshift=2.25cm] (1.35,0)--(1.9,1.8)
[xshift=2.25cm] (1.35,0)--(1.9,1.8)
[xshift=4.1cm] (1.35,0)--(1.9,1.8)
[xshift=2.25cm] (1.35,0)--(1.9,1.8)
[xshift=2.25cm] (1.35,0)--(1.9,1.8);

\draw (0.8,1.8)--(1.35,0)
[xshift=2.25cm]  (0.8,1.8)--(1.35,0)
[xshift=2.25cm]  (0.8,1.8)--(1.35,0)
[xshift=4.1cm]  (0.8,1.8)--(1.35,0)
[xshift=2.25cm]  (0.8,1.8)--(1.35,0)
[xshift=2.25cm]  (0.8,1.8)--(1.35,0)
[xshift=4.1cm]  (0.8,1.8)--(1.35,0)
[xshift=2.25cm]  (0.8,1.8)--(1.35,0)
[xshift=2.25cm]  (0.8,1.8)--(1.35,0)
[xshift=4.1cm]  (0.8,1.8)--(1.35,0)
[xshift=2.25cm]  (0.8,1.8)--(1.35,0)
[xshift=2.25cm]  (0.8,1.8)--(1.35,0)
[xshift=4.1cm]  (0.8,1.8)--(1.35,0)
[xshift=2.25cm]  (0.8,1.8)--(1.35,0)
[xshift=2.25cm]  (0.8,1.8)--(1.35,0);
\draw [densely dotted](7.3,3.6)--(8.1,3.6)
[yshift=-1.8cm] (7.3,3.6)--(8.1,3.6)
[yshift=-1.8cm] (7.3,3.6)--(8.1,3.6)
[xshift=8.6cm] (7.3,3.6)--(8.1,3.6)
[yshift=1.8cm] (7.3,3.6)--(8.1,3.6)
[yshift=1.8cm] (7.3,3.6)--(8.1,3.6)
[xshift=8.6cm] (7.3,3.6)--(8.1,3.6)
[yshift=-1.8cm] (7.3,3.6)--(8.1,3.6)
[yshift=-1.8cm] (7.3,3.6)--(8.1,3.6)
[xshift=8.6cm] (7.3,3.6)--(8.1,3.6)
[yshift=1.8cm] (7.3,3.6)--(8.1,3.6)
[yshift=1.8cm] (7.3,3.6)--(8.1,3.6)
[xshift=8.6cm] (7.3,3.6)--(8.1,3.6)
[yshift=-1.8cm] (7.3,3.6)--(8.1,3.6)
[yshift=-1.8cm] (7.3,3.6)--(8.1,3.6);
%
\fill[black] (2,9) circle (3pt)
[xshift=8.6cm] (2,9) circle (3pt)
[xshift=8.6cm] (2,9) circle (3pt)
[xshift=8.6cm] (2,9) circle (3pt)
[xshift=8.6cm] (2,9) circle (3pt);
\draw[densely dashed] (2,9)--(1.35,4)
[xshift=8.6cm] (2,9)--(1.35,4)
[xshift=8.6cm] (2,9)--(1.35,4)
[xshift=8.6cm] ((2,9)--(1.35,4)
[xshift=8.6cm] (2,9)--(1.35,4);
\draw[densely dashed] (0,4)--(2,9)
[xshift=8.6cm] (0,4)--(2,9)
[xshift=8.6cm] (0,4)--(2,9)
[xshift=8.6cm] (0,4)--(2,9)
[xshift=8.6cm] (0,4)--(2,9);
\draw[densely dashed] (2,9)--(6,4)
[xshift=8.6cm] (2,9)--(6,4)
[xshift=8.6cm] (2,9)--(6,4)
[xshift=8.6cm] (2,9)--(6,4)
[xshift=8.6cm] (2,9)--(6,4);

\draw[densely dashed] (3.6,4)--(2,9)
[xshift=8.6cm] (3.6,4)--(2,9)
[xshift=8.6cm] (3.6,4)--(2,9)
[xshift=8.6cm] (3.6,4)--(2,9)
[xshift=8.6cm] (3.6,4)--(2,9);
\draw [densely dotted](4.8,0)--(5.3,0)
[xshift=8.6cm] (4.8,0)--(5.3,0);
\draw[red] (1.35,0)--(19.2,9);
\draw[red] (3.6,0)--(27.8,9);
\draw[red] (5.85,0)--(36.4,9);
\node (top1)[above] at(2,9){\footnotesize$(n,\top)$};
\node (top2)[above] at(10.6,9){\footnotesize$(k,\top)$};
\node (top3)[above] at(19.2,9){\footnotesize$(f_{n,k}(1),\top)$};
\node (top4)[above] at(27.8,9){\footnotesize$(f_{n,k}(2),\top)$};
\node (top5)[above] at(36.4,9){\footnotesize$(f_{n,k}(x),\top)$};
\node (d1) at(4,-1.3){\footnotesize$n$};
\node (d2) at(13,-1.3){\footnotesize$k$};
\node (d3) at(21.5,-1.3){\footnotesize$f_{n,k}(1)$};
\node (d3) at(30,-1.3){\footnotesize$f_{n,k}(2)$};
\node (d3) at(38.5,-1.3){\footnotesize$f_{n,k}(x)$};
\node (l1) [below] at(1.35,0){\tiny$(n,1)$};
\node (l2) [below] at(3.6,0){\tiny$(n,2)$};
\node (l3) [below] at(5.85,0){\tiny$(n,x)$};
\node (l1) [below] at(9.95,0){\tiny$(k,1)$};
\node (l2) [below] at(12.2,0){\tiny$(k,2)$};
\node (l3) [below] at(14.45,0){\tiny$(k,x)$};
\end{tikzpicture}

  \caption{ The strict order $<_{2}$}\label{F.2}
\end{figure}
The red lines in Fig.\ref{F.2} illustrate three specific cases:
$(n,1)<_{2}(f_{n,k}(1),\top)$, $(n,2)<_{2}(f_{n,k}(2),\top)$ and $(n,x)<_{2}(f_{n,k}(x),\top)$.
\begin{figure}[H]
  \centering
  \begin{tikzpicture} [line width=0.5pt,scale=0.37]
\fill[red,opacity=0.2] (-0.2,-0.25) rectangle (0.2,4.25)
[xshift=8.6cm] (-0.2,-0.25) rectangle (0.2,4.25)
[xshift=8.6cm] (-0.2,-0.25) rectangle (0.2,4.25)
[xshift=8.6cm] (-0.2,-0.25) rectangle (0.2,4.25)
[xshift=8.6cm] (-0.2,-0.25) rectangle (0.2,4.25);
\fill[blue,opacity=0.2] (0.4,-0.25) rectangle (8.2,4.25)
[xshift=8.6cm] (0.4,-0.25) rectangle (8.2,4.25)
[xshift=8.6cm] (0.4,-0.25) rectangle (8.2,4.25)
[xshift=8.6cm] (0.4,-0.25) rectangle (8.2,4.25)
[xshift=8.6cm] (0.4,-0.25) rectangle (8.2,4.25);
\fill[black] (0,0) circle (3pt)
[xshift=8.6cm] (0,0) circle (3pt)
[xshift=8.6cm] (0,0) circle (3pt)
[xshift=8.6cm] (0,0) circle (3pt)
[xshift=8.6cm] (0,0) circle (3pt);
\fill[black] (0,1) circle (3pt)
[xshift=8.6cm] (0,1) circle (3pt)
[xshift=8.6cm] (0,1) circle (3pt)
[xshift=8.6cm] (0,1) circle (3pt)
[xshift=8.6cm] (0,1) circle (3pt);
\fill[black] (0,4) circle (3pt)
[xshift=8.6cm] (0,4) circle (3pt)
[xshift=8.6cm] (0,4) circle (3pt)
[xshift=8.6cm] (0,4) circle (3pt)
[xshift=8.6cm] (0,4) circle (3pt);
\draw (0,0)--(0,1)
[xshift=8.6cm] (0,0)--(0,1)
[xshift=8.6cm] (0,0)--(0,1)
[xshift=8.6cm] (0,0)--(0,1)
[xshift=8.6cm] (0,0)--(0,1);
\draw [densely dashed](0,1)--(0,4)
[xshift=8.6cm] (0,1)--(0,4)
[xshift=8.6cm] (0,1)--(0,4)
[xshift=8.6cm] (0,0)--(0,4)
[xshift=8.6cm] (0,0)--(0,4);

\fill[black] (0.6,3.6) circle (3pt)
[xshift=1.1cm] (0.6,3.6) circle (3pt)
[xshift=1.15cm] (0.6,3.6) circle (3pt)
[xshift=1.1cm] (0.6,3.6) circle (3pt)
[xshift=1.15cm] (0.6,3.6) circle (3pt)
[xshift=1.1cm] (0.6,3.6) circle (3pt)
[xshift=3cm] (0.6,3.6) circle (3pt)
[xshift=1.1cm] (0.6,3.6) circle (3pt)
[xshift=1.15cm] (0.6,3.6) circle (3pt)
[xshift=1.1cm] (0.6,3.6) circle (3pt)
[xshift=1.15cm] (0.6,3.6) circle (3pt)
[xshift=1.1cm] (0.6,3.6) circle (3pt)
[xshift=3cm] (0.6,3.6) circle (3pt)
[xshift=1.1cm] (0.6,3.6) circle (3pt)
[xshift=1.15cm] (0.6,3.6) circle (3pt)
[xshift=1.1cm] (0.6,3.6) circle (3pt)
[xshift=1.15cm] (0.6,3.6) circle (3pt)
[xshift=1.1cm] (0.6,3.6) circle (3pt)
[xshift=3cm] (0.6,3.6) circle (3pt)
[xshift=1.1cm] (0.6,3.6) circle (3pt)
[xshift=1.15cm] (0.6,3.6) circle (3pt)
[xshift=1.1cm] (0.6,3.6) circle (3pt)
[xshift=1.15cm] (0.6,3.6) circle (3pt)
[xshift=1.1cm] (0.6,3.6) circle (3pt)
[xshift=3cm] (0.6,3.6) circle (3pt)
[xshift=1.1cm] (0.6,3.6) circle (3pt)
[xshift=1.15cm] (0.6,3.6) circle (3pt)
[xshift=1.1cm] (0.6,3.6) circle (3pt)
[xshift=1.15cm] (0.6,3.6) circle (3pt)
[xshift=1.1cm] (0.6,3.6) circle (3pt);

\fill[black] (0.8,1.8) circle (3pt)
[xshift=1.1cm](0.8,1.8) circle (3pt)
[xshift=1.15cm](0.8,1.8) circle (3pt)
[xshift=1.1cm](0.8,1.8) circle (3pt)
[xshift=1.15cm](0.8,1.8) circle (3pt)
[xshift=1.1cm](0.8,1.8) circle (3pt)
[xshift=3cm] (0.8,1.8) circle (3pt)
[xshift=1.1cm] (0.8,1.8) circle (3pt)
[xshift=1.15cm](0.8,1.8) circle (3pt)
[xshift=1.1cm](0.8,1.8) circle (3pt)
[xshift=1.15cm] (0.8,1.8) circle (3pt)
[xshift=1.1cm] (0.8,1.8) circle (3pt)
[xshift=3cm] (0.8,1.8) circle (3pt)
[xshift=1.1cm] (0.8,1.8) circle (3pt)
[xshift=1.15cm](0.8,1.8) circle (3pt)
[xshift=1.1cm](0.8,1.8) circle (3pt)
[xshift=1.15cm] (0.8,1.8) circle (3pt)
[xshift=1.1cm] (0.8,1.8) circle (3pt)
[xshift=3cm] (0.8,1.8) circle (3pt)
[xshift=1.1cm] (0.8,1.8) circle (3pt)
[xshift=1.15cm](0.8,1.8) circle (3pt)
[xshift=1.1cm](0.8,1.8) circle (3pt)
[xshift=1.15cm] (0.8,1.8) circle (3pt)
[xshift=1.1cm] (0.8,1.8) circle (3pt)
[xshift=3cm] (0.8,1.8) circle (3pt)
[xshift=1.1cm] (0.8,1.8) circle (3pt)
[xshift=1.15cm](0.8,1.8) circle (3pt)
[xshift=1.1cm](0.8,1.8) circle (3pt)
[xshift=1.15cm] (0.8,1.8) circle (3pt)
[xshift=1.1cm] (0.8,1.8) circle (3pt);

\fill[black] (1,3.6) circle (3pt)
[xshift=1.1cm](1,3.6) circle (3pt)
[xshift=1.15cm](1,3.6) circle (3pt)
[xshift=1.1cm](1,3.6) circle (3pt)
[xshift=1.15cm](1,3.6) circle (3pt)
[xshift=1.1cm] (1,3.6) circle (3pt)
[xshift=3cm](1,3.6) circle (3pt)
[xshift=1.1cm](1,3.6) circle (3pt)
[xshift=1.15cm](1,3.6) circle (3pt)
[xshift=1.1cm](1,3.6) circle (3pt)
[xshift=1.15cm] (1,3.6) circle (3pt)
[xshift=1.1cm] (1,3.6) circle (3pt)
[xshift=3cm](1,3.6) circle (3pt)
[xshift=1.1cm](1,3.6) circle (3pt)
[xshift=1.15cm](1,3.6) circle (3pt)
[xshift=1.1cm](1,3.6) circle (3pt)
[xshift=1.15cm] (1,3.6) circle (3pt)
[xshift=1.1cm] (1,3.6) circle (3pt)
[xshift=3cm](1,3.6) circle (3pt)
[xshift=1.1cm](1,3.6) circle (3pt)
[xshift=1.15cm](1,3.6) circle (3pt)
[xshift=1.1cm](1,3.6) circle (3pt)
[xshift=1.15cm] (1,3.6) circle (3pt)
[xshift=1.1cm] (1,3.6) circle (3pt)
[xshift=3cm](1,3.6) circle (3pt)
[xshift=1.1cm](1,3.6) circle (3pt)
[xshift=1.15cm](1,3.6) circle (3pt)
[xshift=1.1cm](1,3.6) circle (3pt)
[xshift=1.15cm] (1,3.6) circle (3pt)
[xshift=1.1cm] (1,3.6) circle (3pt);
\draw (0.6,3.6)--(0.8,1.8)
[xshift=1.1cm](0.6,3.6)--(0.8,1.8)
[xshift=1.15cm](0.6,3.6)--(0.8,1.8)
[xshift=1.1cm](0.6,3.6)--(0.8,1.8)
[xshift=1.15cm](0.6,3.6)--(0.8,1.8)
[xshift=1.1cm](0.6,3.6)--(0.8,1.8)
[xshift=3cm] (0.6,3.6)--(0.8,1.8)
[xshift=1.1cm](0.6,3.6)--(0.8,1.8)
[xshift=1.15cm](0.6,3.6)--(0.8,1.8)
[xshift=1.1cm](0.6,3.6)--(0.8,1.8)
[xshift=1.15cm] (0.6,3.6)--(0.8,1.8)
[xshift=1.1cm] (0.6,3.6)--(0.8,1.8)
[xshift=3cm] (0.6,3.6)--(0.8,1.8)
[xshift=1.1cm](0.6,3.6)--(0.8,1.8)
[xshift=1.15cm](0.6,3.6)--(0.8,1.8)
[xshift=1.1cm](0.6,3.6)--(0.8,1.8)
[xshift=1.15cm] (0.6,3.6)--(0.8,1.8)
[xshift=1.1cm] (0.6,3.6)--(0.8,1.8)
[xshift=3cm] (0.6,3.6)--(0.8,1.8)
[xshift=1.1cm](0.6,3.6)--(0.8,1.8)
[xshift=1.15cm](0.6,3.6)--(0.8,1.8)
[xshift=1.1cm](0.6,3.6)--(0.8,1.8)
[xshift=1.15cm] (0.6,3.6)--(0.8,1.8)
[xshift=1.1cm] (0.6,3.6)--(0.8,1.8)
[xshift=3cm] (0.6,3.6)--(0.8,1.8)
[xshift=1.1cm](0.6,3.6)--(0.8,1.8)
[xshift=1.15cm](0.6,3.6)--(0.8,1.8)
[xshift=1.1cm](0.6,3.6)--(0.8,1.8)
[xshift=1.15cm] (0.6,3.6)--(0.8,1.8)
[xshift=1.1cm] (0.6,3.6)--(0.8,1.8);

\draw (0.8,1.8)--(1,3.6)
[xshift=1.1cm](0.8,1.8)--(1,3.6)
[xshift=1.15cm](0.8,1.8)--(1,3.6)
[xshift=1.1cm](0.8,1.8)--(1,3.6)
[xshift=1.15cm](0.8,1.8)--(1,3.6)
[xshift=1.1cm](0.8,1.8)--(1,3.6)
[xshift=3cm] (0.8,1.8)--(1,3.6)
[xshift=1.1cm](0.8,1.8)--(1,3.6)
[xshift=1.15cm](0.8,1.8)--(1,3.6)
[xshift=1.1cm](0.8,1.8)--(1,3.6)
[xshift=1.15cm] (0.8,1.8)--(1,3.6)
[xshift=1.1cm] (0.8,1.8)--(1,3.6)
[xshift=3cm] (0.8,1.8)--(1,3.6)
[xshift=1.1cm](0.8,1.8)--(1,3.6)
[xshift=1.15cm](0.8,1.8)--(1,3.6)
[xshift=1.1cm](0.8,1.8)--(1,3.6)
[xshift=1.15cm] (0.8,1.8)--(1,3.6)
[xshift=1.1cm] (0.8,1.8)--(1,3.6)
[xshift=3cm] (0.8,1.8)--(1,3.6)
[xshift=1.1cm](0.8,1.8)--(1,3.6)
[xshift=1.15cm](0.8,1.8)--(1,3.6)
[xshift=1.1cm](0.8,1.8)--(1,3.6)
[xshift=1.15cm] (0.8,1.8)--(1,3.6)
[xshift=1.1cm] (0.8,1.8)--(1,3.6)
[xshift=3cm] (0.8,1.8)--(1,3.6)
[xshift=1.1cm](0.8,1.8)--(1,3.6)
[xshift=1.15cm](0.8,1.8)--(1,3.6)
[xshift=1.1cm](0.8,1.8)--(1,3.6)
[xshift=1.15cm] (0.8,1.8)--(1,3.6)
[xshift=1.1cm] (0.8,1.8)--(1,3.6);
\draw [densely dotted](1.15,3.6)--(1.5,3.6)
[xshift=1.1cm] (1.15,3.6)--(1.5,3.6)
[yshift=-1.8cm] (1.15,3.6)--(1.5,3.6)
[shift={(1.15cm, 1.8cm)}](1.15,3.6)--(1.5,3.6)
[xshift=1.1cm] (1.15,3.6)--(1.5,3.6)
[yshift=-1.8cm] (1.15,3.6)--(1.5,3.6)
[shift={(1.15cm,1.8cm)}] (1.15,3.6)--(1.5,3.6)
[xshift=1.1cm](1.15,3.6)--(1.5,3.6)
[yshift=-1.8cm] (1.15,3.6)--(1.5,3.6)
[shift={(3cm,1.8cm)}] (1.15,3.6)--(1.5,3.6)
[xshift=1.1cm] (1.15,3.6)--(1.5,3.6)
[yshift=-1.8cm] (1.15,3.6)--(1.5,3.6)
[shift={(1.15cm, 1.8cm)}](1.15,3.6)--(1.5,3.6)
[xshift=1.1cm] (1.15,3.6)--(1.5,3.6)
[yshift=-1.8cm] (1.15,3.6)--(1.5,3.6)
[shift={(1.15cm,1.8cm)}] (1.15,3.6)--(1.5,3.6)
[xshift=1.1cm](1.15,3.6)--(1.5,3.6)
[yshift=-1.8cm] (1.15,3.6)--(1.5,3.6)
[shift={(3cm,1.8cm)}] (1.15,3.6)--(1.5,3.6)
[xshift=1.1cm] (1.15,3.6)--(1.5,3.6)
[yshift=-1.8cm] (1.15,3.6)--(1.5,3.6)
[shift={(1.15cm, 1.8cm)}](1.15,3.6)--(1.5,3.6)
[xshift=1.1cm] (1.15,3.6)--(1.5,3.6)
[yshift=-1.8cm] (1.15,3.6)--(1.5,3.6)
[shift={(1.15cm,1.8cm)}] (1.15,3.6)--(1.5,3.6)
[xshift=1.1cm](1.15,3.6)--(1.5,3.6)
[yshift=-1.8cm] (1.15,3.6)--(1.5,3.6)
[shift={(3cm,1.8cm)}] (1.15,3.6)--(1.5,3.6)
[xshift=1.1cm] (1.15,3.6)--(1.5,3.6)
[yshift=-1.8cm] (1.15,3.6)--(1.5,3.6)
[shift={(1.15cm, 1.8cm)}](1.15,3.6)--(1.5,3.6)
[xshift=1.1cm] (1.15,3.6)--(1.5,3.6)
[yshift=-1.8cm] (1.15,3.6)--(1.5,3.6)
[shift={(1.15cm,1.8cm)}] (1.15,3.6)--(1.5,3.6)
[xshift=1.1cm](1.15,3.6)--(1.5,3.6)
[yshift=-1.8cm] (1.15,3.6)--(1.5,3.6)
[shift={(3cm,1.8cm)}] (1.15,3.6)--(1.5,3.6)
[xshift=1.1cm] (1.15,3.6)--(1.5,3.6)
[yshift=-1.8cm] (1.15,3.6)--(1.5,3.6)
[shift={(1.15cm, 1.8cm)}](1.15,3.6)--(1.5,3.6)
[xshift=1.1cm] (1.15,3.6)--(1.5,3.6)
[yshift=-1.8cm] (1.15,3.6)--(1.5,3.6)
[shift={(1.15cm,1.8cm)}] (1.15,3.6)--(1.5,3.6)
[xshift=1.1cm](1.15,3.6)--(1.5,3.6)
[yshift=-1.8cm] (1.15,3.6)--(1.5,3.6);
\fill[black] (1.35,0) circle (3pt)
[xshift=2.25cm](1.35,0) circle (3pt)
[xshift=2.25cm] (1.35,0) circle (3pt)
[xshift=4.1cm](1.35,0) circle (3pt)
[xshift=2.25cm] (1.35,0) circle (3pt)
[xshift=2.25cm](1.35,0) circle (3pt)
[xshift=4.1cm] (1.35,0) circle (3pt)
[xshift=2.25cm] (1.35,0) circle (3pt)
[xshift=2.25cm](1.35,0) circle (3pt)
[xshift=4.1cm](1.35,0) circle (3pt)
[xshift=2.25cm] (1.35,0) circle (3pt)
[xshift=2.25cm](1.35,0) circle (3pt)
[xshift=4.1cm](1.35,0) circle (3pt)
[xshift=2.25cm] (1.35,0) circle (3pt)
[xshift=2.25cm](1.35,0) circle (3pt);

\draw (1.35,0)--(1.9,1.8)
[xshift=2.25cm] (1.35,0)--(1.9,1.8)
[xshift=2.25cm] (1.35,0)--(1.9,1.8)
[xshift=4.1cm] (1.35,0)--(1.9,1.8)
[xshift=2.25cm] (1.35,0)--(1.9,1.8)
[xshift=2.25cm] (1.35,0)--(1.9,1.8)
[xshift=4.1cm] (1.35,0)--(1.9,1.8)
[xshift=2.25cm] (1.35,0)--(1.9,1.8)
[xshift=2.25cm] (1.35,0)--(1.9,1.8)
[xshift=4.1cm] (1.35,0)--(1.9,1.8)
[xshift=2.25cm] (1.35,0)--(1.9,1.8)
[xshift=2.25cm] (1.35,0)--(1.9,1.8)
[xshift=4.1cm] (1.35,0)--(1.9,1.8)
[xshift=2.25cm] (1.35,0)--(1.9,1.8)
[xshift=2.25cm] (1.35,0)--(1.9,1.8);

\draw (0.8,1.8)--(1.35,0)
[xshift=2.25cm]  (0.8,1.8)--(1.35,0)
[xshift=2.25cm]  (0.8,1.8)--(1.35,0)
[xshift=4.1cm]  (0.8,1.8)--(1.35,0)
[xshift=2.25cm]  (0.8,1.8)--(1.35,0)
[xshift=2.25cm]  (0.8,1.8)--(1.35,0)
[xshift=4.1cm]  (0.8,1.8)--(1.35,0)
[xshift=2.25cm]  (0.8,1.8)--(1.35,0)
[xshift=2.25cm]  (0.8,1.8)--(1.35,0)
[xshift=4.1cm]  (0.8,1.8)--(1.35,0)
[xshift=2.25cm]  (0.8,1.8)--(1.35,0)
[xshift=2.25cm]  (0.8,1.8)--(1.35,0)
[xshift=4.1cm]  (0.8,1.8)--(1.35,0)
[xshift=2.25cm]  (0.8,1.8)--(1.35,0)
[xshift=2.25cm]  (0.8,1.8)--(1.35,0);
\draw [densely dotted](7.3,3.6)--(8.1,3.6)
[yshift=-1.8cm] (7.3,3.6)--(8.1,3.6)
[yshift=-1.8cm] (7.3,3.6)--(8.1,3.6)
[xshift=8.6cm] (7.3,3.6)--(8.1,3.6)
[yshift=1.8cm] (7.3,3.6)--(8.1,3.6)
[yshift=1.8cm] (7.3,3.6)--(8.1,3.6)
[xshift=8.6cm] (7.3,3.6)--(8.1,3.6)
[yshift=-1.8cm] (7.3,3.6)--(8.1,3.6)
[yshift=-1.8cm] (7.3,3.6)--(8.1,3.6)
[xshift=8.6cm] (7.3,3.6)--(8.1,3.6)
[yshift=1.8cm] (7.3,3.6)--(8.1,3.6)
[yshift=1.8cm] (7.3,3.6)--(8.1,3.6)
[xshift=8.6cm] (7.3,3.6)--(8.1,3.6)
[yshift=-1.8cm] (7.3,3.6)--(8.1,3.6)
[yshift=-1.8cm] (7.3,3.6)--(8.1,3.6);
%
\fill[black] (2,9) circle (3pt)
[xshift=8.6cm] (2,9) circle (3pt)
[xshift=8.6cm] (2,9) circle (3pt)
[xshift=8.6cm] (2,9) circle (3pt)
[xshift=8.6cm] (2,9) circle (3pt);
\draw[densely dashed] (2,9)--(1.35,4)
[xshift=8.6cm] (2,9)--(1.35,4)
[xshift=8.6cm] (2,9)--(1.35,4)
[xshift=8.6cm] ((2,9)--(1.35,4)
[xshift=8.6cm] (2,9)--(1.35,4);
\draw[densely dashed] (0,4)--(2,9)
[xshift=8.6cm] (0,4)--(2,9)
[xshift=8.6cm] (0,4)--(2,9)
[xshift=8.6cm] (0,4)--(2,9)
[xshift=8.6cm] (0,4)--(2,9);
\draw[densely dashed] (2,9)--(6,4)
[xshift=8.6cm] (2,9)--(6,4)
[xshift=8.6cm] (2,9)--(6,4)
[xshift=8.6cm] (2,9)--(6,4)
[xshift=8.6cm] (2,9)--(6,4);

\draw[densely dashed] (3.6,4)--(2,9)
[xshift=8.6cm] (3.6,4)--(2,9)
[xshift=8.6cm] (3.6,4)--(2,9)
[xshift=8.6cm] (3.6,4)--(2,9)
[xshift=8.6cm] (3.6,4)--(2,9);
\draw [densely dotted](4.8,0)--(5.3,0)
[xshift=8.6cm] (4.8,0)--(5.3,0);
\draw[red] (8.6,0)--(19.2,9);
\draw[red] (8.6,1)--(27.8,9);
\draw[red] (8.6,4)--(36.4,9);
\node (top1)[above] at(2,9){\footnotesize$(d,\top)$};
\node (top2)[above] at(10.6,9){\footnotesize$(n,\top)$};
\node (top3)[above] at(19.2,9){\footnotesize$(f_{d,n}(1),\top)$};
\node (top4)[above] at(27.8,9){\footnotesize$(f_{d,n}(2),\top)$};
\node (top5)[above] at(36.4,9){\footnotesize$(f_{d,n}(x),\top)$};
\node (d1) at(4,-1.3){\footnotesize$d$};
\node (d2) at(13,-1.3){\footnotesize$n$};
\node (d3) at(21.5,-1.3){\footnotesize$f_{d,n}(1)$};
\node (d3) at(30,-1.3){\footnotesize$f_{d,n}(2)$};
\node (d3) at(38.5,-1.3){\footnotesize$f_{d,n}(x)$};
\node (l1) [below] at(8.6,0){\tiny$(n,1)$};
\node (l2) [below] at(8.6,1.4){\tiny$(n,2)$};
\node (l3) [below] at(8.6,4.4){\tiny$(n,x)$};
\node (l4) [below] at(1.35,0){\tiny$(d,1)$};
\node (l5) [below] at(3.6,0){\tiny$(d,2)$};
\node (l6) [below] at(5.85,0){\tiny$(d,x)$};
\end{tikzpicture}

  \caption{ The strict order $<_{3}$}\label{F.3}
\end{figure}
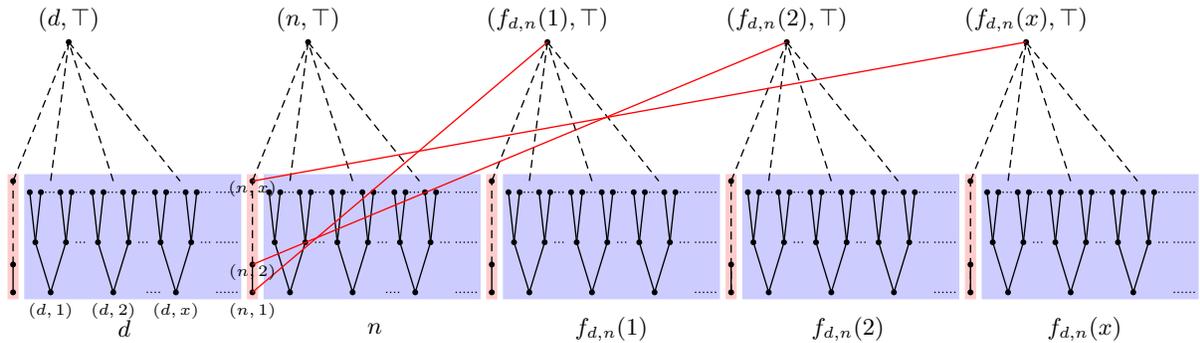
The red lines in Fig.\ref{F.3} illustrate the cases:
$(n,1)<_{3}(f_{d,n}(1),\top)$, $(n,2)<_{3}(f_{d,n}(2),\top)$ and $(n,x)<_{3}(f_{d,n}(x),\top)$.
\begin{figure}[H]
  \centering
  \begin{tikzpicture} [line width=0.5pt,scale=0.37]
\fill[red,opacity=0.2] (-0.2,-0.25) rectangle (0.2,4.25)
[xshift=8.6cm] (-0.2,-0.25) rectangle (0.2,4.25)
[xshift=8.6cm] (-0.2,-0.25) rectangle (0.2,4.25)
[xshift=8.6cm] (-0.2,-0.25) rectangle (0.2,4.25)
[xshift=8.6cm] (-0.2,-0.25) rectangle (0.2,4.25);
\fill[blue,opacity=0.2] (0.4,-0.25) rectangle (8.2,4.25)
[xshift=8.6cm] (0.4,-0.25) rectangle (8.2,4.25)
[xshift=8.6cm] (0.4,-0.25) rectangle (8.2,4.25)
[xshift=8.6cm] (0.4,-0.25) rectangle (8.2,4.25)
[xshift=8.6cm] (0.4,-0.25) rectangle (8.2,4.25);
\fill[black] (0,0) circle (3pt)
[xshift=8.6cm] (0,0) circle (3pt)
[xshift=8.6cm] (0,0) circle (3pt)
[xshift=8.6cm] (0,0) circle (3pt)
[xshift=8.6cm] (0,0) circle (3pt);
\fill[black] (0,1) circle (3pt)
[xshift=8.6cm] (0,1) circle (3pt)
[xshift=8.6cm] (0,1) circle (3pt)
[xshift=8.6cm] (0,1) circle (3pt)
[xshift=8.6cm] (0,1) circle (3pt);
\fill[black] (0,4) circle (3pt)
[xshift=8.6cm] (0,4) circle (3pt)
[xshift=8.6cm] (0,4) circle (3pt)
[xshift=8.6cm] (0,4) circle (3pt)
[xshift=8.6cm] (0,4) circle (3pt);
\draw (0,0)--(0,1)
[xshift=8.6cm] (0,0)--(0,1)
[xshift=8.6cm] (0,0)--(0,1)
[xshift=8.6cm] (0,0)--(0,1)
[xshift=8.6cm] (0,0)--(0,1);
\draw [densely dashed](0,1)--(0,4)
[xshift=8.6cm] (0,1)--(0,4)
[xshift=8.6cm] (0,1)--(0,4)
[xshift=8.6cm] (0,0)--(0,4)
[xshift=8.6cm] (0,0)--(0,4);

\fill[black] (0.6,3.6) circle (3pt)
[xshift=1.1cm] (0.6,3.6) circle (3pt)
[xshift=1.15cm] (0.6,3.6) circle (3pt)
[xshift=1.1cm] (0.6,3.6) circle (3pt)
[xshift=1.15cm] (0.6,3.6) circle (3pt)
[xshift=1.1cm] (0.6,3.6) circle (3pt)
[xshift=3cm] (0.6,3.6) circle (3pt)
[xshift=1.1cm] (0.6,3.6) circle (3pt)
[xshift=1.15cm] (0.6,3.6) circle (3pt)
[xshift=1.1cm] (0.6,3.6) circle (3pt)
[xshift=1.15cm] (0.6,3.6) circle (3pt)
[xshift=1.1cm] (0.6,3.6) circle (3pt)
[xshift=3cm] (0.6,3.6) circle (3pt)
[xshift=1.1cm] (0.6,3.6) circle (3pt)
[xshift=1.15cm] (0.6,3.6) circle (3pt)
[xshift=1.1cm] (0.6,3.6) circle (3pt)
[xshift=1.15cm] (0.6,3.6) circle (3pt)
[xshift=1.1cm] (0.6,3.6) circle (3pt)
[xshift=3cm] (0.6,3.6) circle (3pt)
[xshift=1.1cm] (0.6,3.6) circle (3pt)
[xshift=1.15cm] (0.6,3.6) circle (3pt)
[xshift=1.1cm] (0.6,3.6) circle (3pt)
[xshift=1.15cm] (0.6,3.6) circle (3pt)
[xshift=1.1cm] (0.6,3.6) circle (3pt)
[xshift=3cm] (0.6,3.6) circle (3pt)
[xshift=1.1cm] (0.6,3.6) circle (3pt)
[xshift=1.15cm] (0.6,3.6) circle (3pt)
[xshift=1.1cm] (0.6,3.6) circle (3pt)
[xshift=1.15cm] (0.6,3.6) circle (3pt)
[xshift=1.1cm] (0.6,3.6) circle (3pt);

\fill[black] (0.8,1.8) circle (3pt)
[xshift=1.1cm](0.8,1.8) circle (3pt)
[xshift=1.15cm](0.8,1.8) circle (3pt)
[xshift=1.1cm](0.8,1.8) circle (3pt)
[xshift=1.15cm](0.8,1.8) circle (3pt)
[xshift=1.1cm](0.8,1.8) circle (3pt)
[xshift=3cm] (0.8,1.8) circle (3pt)
[xshift=1.1cm] (0.8,1.8) circle (3pt)
[xshift=1.15cm](0.8,1.8) circle (3pt)
[xshift=1.1cm](0.8,1.8) circle (3pt)
[xshift=1.15cm] (0.8,1.8) circle (3pt)
[xshift=1.1cm] (0.8,1.8) circle (3pt)
[xshift=3cm] (0.8,1.8) circle (3pt)
[xshift=1.1cm] (0.8,1.8) circle (3pt)
[xshift=1.15cm](0.8,1.8) circle (3pt)
[xshift=1.1cm](0.8,1.8) circle (3pt)
[xshift=1.15cm] (0.8,1.8) circle (3pt)
[xshift=1.1cm] (0.8,1.8) circle (3pt)
[xshift=3cm] (0.8,1.8) circle (3pt)
[xshift=1.1cm] (0.8,1.8) circle (3pt)
[xshift=1.15cm](0.8,1.8) circle (3pt)
[xshift=1.1cm](0.8,1.8) circle (3pt)
[xshift=1.15cm] (0.8,1.8) circle (3pt)
[xshift=1.1cm] (0.8,1.8) circle (3pt)
[xshift=3cm] (0.8,1.8) circle (3pt)
[xshift=1.1cm] (0.8,1.8) circle (3pt)
[xshift=1.15cm](0.8,1.8) circle (3pt)
[xshift=1.1cm](0.8,1.8) circle (3pt)
[xshift=1.15cm] (0.8,1.8) circle (3pt)
[xshift=1.1cm] (0.8,1.8) circle (3pt);

\fill[black] (1,3.6) circle (3pt)
[xshift=1.1cm](1,3.6) circle (3pt)
[xshift=1.15cm](1,3.6) circle (3pt)
[xshift=1.1cm](1,3.6) circle (3pt)
[xshift=1.15cm](1,3.6) circle (3pt)
[xshift=1.1cm] (1,3.6) circle (3pt)
[xshift=3cm](1,3.6) circle (3pt)
[xshift=1.1cm](1,3.6) circle (3pt)
[xshift=1.15cm](1,3.6) circle (3pt)
[xshift=1.1cm](1,3.6) circle (3pt)
[xshift=1.15cm] (1,3.6) circle (3pt)
[xshift=1.1cm] (1,3.6) circle (3pt)
[xshift=3cm](1,3.6) circle (3pt)
[xshift=1.1cm](1,3.6) circle (3pt)
[xshift=1.15cm](1,3.6) circle (3pt)
[xshift=1.1cm](1,3.6) circle (3pt)
[xshift=1.15cm] (1,3.6) circle (3pt)
[xshift=1.1cm] (1,3.6) circle (3pt)
[xshift=3cm](1,3.6) circle (3pt)
[xshift=1.1cm](1,3.6) circle (3pt)
[xshift=1.15cm](1,3.6) circle (3pt)
[xshift=1.1cm](1,3.6) circle (3pt)
[xshift=1.15cm] (1,3.6) circle (3pt)
[xshift=1.1cm] (1,3.6) circle (3pt)
[xshift=3cm](1,3.6) circle (3pt)
[xshift=1.1cm](1,3.6) circle (3pt)
[xshift=1.15cm](1,3.6) circle (3pt)
[xshift=1.1cm](1,3.6) circle (3pt)
[xshift=1.15cm] (1,3.6) circle (3pt)
[xshift=1.1cm] (1,3.6) circle (3pt);
\draw (0.6,3.6)--(0.8,1.8)
[xshift=1.1cm](0.6,3.6)--(0.8,1.8)
[xshift=1.15cm](0.6,3.6)--(0.8,1.8)
[xshift=1.1cm](0.6,3.6)--(0.8,1.8)
[xshift=1.15cm](0.6,3.6)--(0.8,1.8)
[xshift=1.1cm](0.6,3.6)--(0.8,1.8)
[xshift=3cm] (0.6,3.6)--(0.8,1.8)
[xshift=1.1cm](0.6,3.6)--(0.8,1.8)
[xshift=1.15cm](0.6,3.6)--(0.8,1.8)
[xshift=1.1cm](0.6,3.6)--(0.8,1.8)
[xshift=1.15cm] (0.6,3.6)--(0.8,1.8)
[xshift=1.1cm] (0.6,3.6)--(0.8,1.8)
[xshift=3cm] (0.6,3.6)--(0.8,1.8)
[xshift=1.1cm](0.6,3.6)--(0.8,1.8)
[xshift=1.15cm](0.6,3.6)--(0.8,1.8)
[xshift=1.1cm](0.6,3.6)--(0.8,1.8)
[xshift=1.15cm] (0.6,3.6)--(0.8,1.8)
[xshift=1.1cm] (0.6,3.6)--(0.8,1.8)
[xshift=3cm] (0.6,3.6)--(0.8,1.8)
[xshift=1.1cm](0.6,3.6)--(0.8,1.8)
[xshift=1.15cm](0.6,3.6)--(0.8,1.8)
[xshift=1.1cm](0.6,3.6)--(0.8,1.8)
[xshift=1.15cm] (0.6,3.6)--(0.8,1.8)
[xshift=1.1cm] (0.6,3.6)--(0.8,1.8)
[xshift=3cm] (0.6,3.6)--(0.8,1.8)
[xshift=1.1cm](0.6,3.6)--(0.8,1.8)
[xshift=1.15cm](0.6,3.6)--(0.8,1.8)
[xshift=1.1cm](0.6,3.6)--(0.8,1.8)
[xshift=1.15cm] (0.6,3.6)--(0.8,1.8)
[xshift=1.1cm] (0.6,3.6)--(0.8,1.8);

\draw (0.8,1.8)--(1,3.6)
[xshift=1.1cm](0.8,1.8)--(1,3.6)
[xshift=1.15cm](0.8,1.8)--(1,3.6)
[xshift=1.1cm](0.8,1.8)--(1,3.6)
[xshift=1.15cm](0.8,1.8)--(1,3.6)
[xshift=1.1cm](0.8,1.8)--(1,3.6)
[xshift=3cm] (0.8,1.8)--(1,3.6)
[xshift=1.1cm](0.8,1.8)--(1,3.6)
[xshift=1.15cm](0.8,1.8)--(1,3.6)
[xshift=1.1cm](0.8,1.8)--(1,3.6)
[xshift=1.15cm] (0.8,1.8)--(1,3.6)
[xshift=1.1cm] (0.8,1.8)--(1,3.6)
[xshift=3cm] (0.8,1.8)--(1,3.6)
[xshift=1.1cm](0.8,1.8)--(1,3.6)
[xshift=1.15cm](0.8,1.8)--(1,3.6)
[xshift=1.1cm](0.8,1.8)--(1,3.6)
[xshift=1.15cm] (0.8,1.8)--(1,3.6)
[xshift=1.1cm] (0.8,1.8)--(1,3.6)
[xshift=3cm] (0.8,1.8)--(1,3.6)
[xshift=1.1cm](0.8,1.8)--(1,3.6)
[xshift=1.15cm](0.8,1.8)--(1,3.6)
[xshift=1.1cm](0.8,1.8)--(1,3.6)
[xshift=1.15cm] (0.8,1.8)--(1,3.6)
[xshift=1.1cm] (0.8,1.8)--(1,3.6)
[xshift=3cm] (0.8,1.8)--(1,3.6)
[xshift=1.1cm](0.8,1.8)--(1,3.6)
[xshift=1.15cm](0.8,1.8)--(1,3.6)
[xshift=1.1cm](0.8,1.8)--(1,3.6)
[xshift=1.15cm] (0.8,1.8)--(1,3.6)
[xshift=1.1cm] (0.8,1.8)--(1,3.6);
\draw [densely dotted](1.15,3.6)--(1.5,3.6)
[xshift=1.1cm] (1.15,3.6)--(1.5,3.6)
[yshift=-1.8cm] (1.15,3.6)--(1.5,3.6)
[shift={(1.15cm, 1.8cm)}](1.15,3.6)--(1.5,3.6)
[xshift=1.1cm] (1.15,3.6)--(1.5,3.6)
[yshift=-1.8cm] (1.15,3.6)--(1.5,3.6)
[shift={(1.15cm,1.8cm)}] (1.15,3.6)--(1.5,3.6)
[xshift=1.1cm](1.15,3.6)--(1.5,3.6)
[yshift=-1.8cm] (1.15,3.6)--(1.5,3.6)
[shift={(3cm,1.8cm)}] (1.15,3.6)--(1.5,3.6)
[xshift=1.1cm] (1.15,3.6)--(1.5,3.6)
[yshift=-1.8cm] (1.15,3.6)--(1.5,3.6)
[shift={(1.15cm, 1.8cm)}](1.15,3.6)--(1.5,3.6)
[xshift=1.1cm] (1.15,3.6)--(1.5,3.6)
[yshift=-1.8cm] (1.15,3.6)--(1.5,3.6)
[shift={(1.15cm,1.8cm)}] (1.15,3.6)--(1.5,3.6)
[xshift=1.1cm](1.15,3.6)--(1.5,3.6)
[yshift=-1.8cm] (1.15,3.6)--(1.5,3.6)
[shift={(3cm,1.8cm)}] (1.15,3.6)--(1.5,3.6)
[xshift=1.1cm] (1.15,3.6)--(1.5,3.6)
[yshift=-1.8cm] (1.15,3.6)--(1.5,3.6)
[shift={(1.15cm, 1.8cm)}](1.15,3.6)--(1.5,3.6)
[xshift=1.1cm] (1.15,3.6)--(1.5,3.6)
[yshift=-1.8cm] (1.15,3.6)--(1.5,3.6)
[shift={(1.15cm,1.8cm)}] (1.15,3.6)--(1.5,3.6)
[xshift=1.1cm](1.15,3.6)--(1.5,3.6)
[yshift=-1.8cm] (1.15,3.6)--(1.5,3.6)
[shift={(3cm,1.8cm)}] (1.15,3.6)--(1.5,3.6)
[xshift=1.1cm] (1.15,3.6)--(1.5,3.6)
[yshift=-1.8cm] (1.15,3.6)--(1.5,3.6)
[shift={(1.15cm, 1.8cm)}](1.15,3.6)--(1.5,3.6)
[xshift=1.1cm] (1.15,3.6)--(1.5,3.6)
[yshift=-1.8cm] (1.15,3.6)--(1.5,3.6)
[shift={(1.15cm,1.8cm)}] (1.15,3.6)--(1.5,3.6)
[xshift=1.1cm](1.15,3.6)--(1.5,3.6)
[yshift=-1.8cm] (1.15,3.6)--(1.5,3.6)
[shift={(3cm,1.8cm)}] (1.15,3.6)--(1.5,3.6)
[xshift=1.1cm] (1.15,3.6)--(1.5,3.6)
[yshift=-1.8cm] (1.15,3.6)--(1.5,3.6)
[shift={(1.15cm, 1.8cm)}](1.15,3.6)--(1.5,3.6)
[xshift=1.1cm] (1.15,3.6)--(1.5,3.6)
[yshift=-1.8cm] (1.15,3.6)--(1.5,3.6)
[shift={(1.15cm,1.8cm)}] (1.15,3.6)--(1.5,3.6)
[xshift=1.1cm](1.15,3.6)--(1.5,3.6)
[yshift=-1.8cm] (1.15,3.6)--(1.5,3.6);
\fill[black] (1.35,0) circle (3pt)
[xshift=2.25cm](1.35,0) circle (3pt)
[xshift=2.25cm] (1.35,0) circle (3pt)
[xshift=4.1cm](1.35,0) circle (3pt)
[xshift=2.25cm] (1.35,0) circle (3pt)
[xshift=2.25cm](1.35,0) circle (3pt)
[xshift=4.1cm] (1.35,0) circle (3pt)
[xshift=2.25cm] (1.35,0) circle (3pt)
[xshift=2.25cm](1.35,0) circle (3pt)
[xshift=4.1cm](1.35,0) circle (3pt)
[xshift=2.25cm] (1.35,0) circle (3pt)
[xshift=2.25cm](1.35,0) circle (3pt)
[xshift=4.1cm](1.35,0) circle (3pt)
[xshift=2.25cm] (1.35,0) circle (3pt)
[xshift=2.25cm](1.35,0) circle (3pt);

\draw (1.35,0)--(1.9,1.8)
[xshift=2.25cm] (1.35,0)--(1.9,1.8)
[xshift=2.25cm] (1.35,0)--(1.9,1.8)
[xshift=4.1cm] (1.35,0)--(1.9,1.8)
[xshift=2.25cm] (1.35,0)--(1.9,1.8)
[xshift=2.25cm] (1.35,0)--(1.9,1.8)
[xshift=4.1cm] (1.35,0)--(1.9,1.8)
[xshift=2.25cm] (1.35,0)--(1.9,1.8)
[xshift=2.25cm] (1.35,0)--(1.9,1.8)
[xshift=4.1cm] (1.35,0)--(1.9,1.8)
[xshift=2.25cm] (1.35,0)--(1.9,1.8)
[xshift=2.25cm] (1.35,0)--(1.9,1.8)
[xshift=4.1cm] (1.35,0)--(1.9,1.8)
[xshift=2.25cm] (1.35,0)--(1.9,1.8)
[xshift=2.25cm] (1.35,0)--(1.9,1.8);

\draw (0.8,1.8)--(1.35,0)
[xshift=2.25cm]  (0.8,1.8)--(1.35,0)
[xshift=2.25cm]  (0.8,1.8)--(1.35,0)
[xshift=4.1cm]  (0.8,1.8)--(1.35,0)
[xshift=2.25cm]  (0.8,1.8)--(1.35,0)
[xshift=2.25cm]  (0.8,1.8)--(1.35,0)
[xshift=4.1cm]  (0.8,1.8)--(1.35,0)
[xshift=2.25cm]  (0.8,1.8)--(1.35,0)
[xshift=2.25cm]  (0.8,1.8)--(1.35,0)
[xshift=4.1cm]  (0.8,1.8)--(1.35,0)
[xshift=2.25cm]  (0.8,1.8)--(1.35,0)
[xshift=2.25cm]  (0.8,1.8)--(1.35,0)
[xshift=4.1cm]  (0.8,1.8)--(1.35,0)
[xshift=2.25cm]  (0.8,1.8)--(1.35,0)
[xshift=2.25cm]  (0.8,1.8)--(1.35,0);
\draw [densely dotted](7.3,3.6)--(8.1,3.6)
[yshift=-1.8cm] (7.3,3.6)--(8.1,3.6)
[yshift=-1.8cm] (7.3,3.6)--(8.1,3.6)
[xshift=8.6cm] (7.3,3.6)--(8.1,3.6)
[yshift=1.8cm] (7.3,3.6)--(8.1,3.6)
[yshift=1.8cm] (7.3,3.6)--(8.1,3.6)
[xshift=8.6cm] (7.3,3.6)--(8.1,3.6)
[yshift=-1.8cm] (7.3,3.6)--(8.1,3.6)
[yshift=-1.8cm] (7.3,3.6)--(8.1,3.6)
[xshift=8.6cm] (7.3,3.6)--(8.1,3.6)
[yshift=1.8cm] (7.3,3.6)--(8.1,3.6)
[yshift=1.8cm] (7.3,3.6)--(8.1,3.6)
[xshift=8.6cm] (7.3,3.6)--(8.1,3.6)
[yshift=-1.8cm] (7.3,3.6)--(8.1,3.6)
[yshift=-1.8cm] (7.3,3.6)--(8.1,3.6);
%
\fill[black] (2,9) circle (3pt)
[xshift=8.6cm] (2,9) circle (3pt)
[xshift=8.6cm] (2,9) circle (3pt)
[xshift=8.6cm] (2,9) circle (3pt)
[xshift=8.6cm] (2,9) circle (3pt);
\draw[densely dashed] (2,9)--(1.35,4)
[xshift=8.6cm] (2,9)--(1.35,4)
[xshift=8.6cm] (2,9)--(1.35,4)
[xshift=8.6cm] ((2,9)--(1.35,4)
[xshift=8.6cm] (2,9)--(1.35,4);
\draw[densely dashed] (0,4)--(2,9)
[xshift=8.6cm] (0,4)--(2,9)
[xshift=8.6cm] (0,4)--(2,9)
[xshift=8.6cm] (0,4)--(2,9)
[xshift=8.6cm] (0,4)--(2,9);
\draw[densely dashed] (2,9)--(6,4)
[xshift=8.6cm] (2,9)--(6,4)
[xshift=8.6cm] (2,9)--(6,4)
[xshift=8.6cm] (2,9)--(6,4)
[xshift=8.6cm] (2,9)--(6,4);

\draw[densely dashed] (3.6,4)--(2,9)
[xshift=8.6cm] (3.6,4)--(2,9)
[xshift=8.6cm] (3.6,4)--(2,9)
[xshift=8.6cm] (3.6,4)--(2,9)
[xshift=8.6cm] (3.6,4)--(2,9);
\draw[red] (17.2,0)--(27.8,9);
\draw[red] (17.2,4)--(36.4,9);
\node (top1)[above] at(2,9){\footnotesize$(a,\top)$};
\node (top2)[above] at(10.6,9){\footnotesize$(b,\top)$};
\node (top3)[above] at(19.2,9){\footnotesize$(f_{a,b}(1),\top)$};
\node (top4)[above] at(27.8,9){\footnotesize$(f_{a,b}(1.1),\top)$};
\node (top5)[above] at(36.4,9){\footnotesize$(f_{a,b}(1.x),\top)$};
\node (d1) at(4,-1.3){\footnotesize$a$};
\node (d2) at(13,-1.3){\footnotesize$b$};
\node (d3) at(21.5,-1.3){\footnotesize$f_{a,b}(1)$};
\node (d3) at(30,-1.3){\footnotesize$f_{a,b}(1.1)$};
\node (d3) at(38.5,-1.3){\footnotesize$f_{a,b}(1.x)$};
\node (l2) [below] at(17.2,0){\tiny$(f_{a,b}(1),1)$};
\node (l3) [below] at(17.2,4){\tiny$(f_{a,b}(1),x)$};
\end{tikzpicture}

  \caption{ The strict order $<_{4}$}\label{F.4}
\end{figure}
The red lines in Fig.\ref{F.4} illustrate the cases: $(f_{a,b}(1),1)<_{4}(f_{a,b}(1.1),\top)$ and $(f_{a,b}(1),x)<_{4}(f_{a,b}(1.x),\top)$.

\begin{figure}[H]
  \centering
  \begin{tikzpicture} [line width=0.5pt,scale=0.37]
\fill[red,opacity=0.2] (-0.2,-0.25) rectangle (0.2,4.25)
[xshift=8.6cm] (-0.2,-0.25) rectangle (0.2,4.25)
[xshift=8.6cm] (-0.2,-0.25) rectangle (0.2,4.25)
[xshift=8.6cm] (-0.2,-0.25) rectangle (0.2,4.25)
[xshift=8.6cm] (-0.2,-0.25) rectangle (0.2,4.25);
\fill[blue,opacity=0.2] (0.4,-0.25) rectangle (8.2,4.25)
[xshift=8.6cm] (0.4,-0.25) rectangle (8.2,4.25)
[xshift=8.6cm] (0.4,-0.25) rectangle (8.2,4.25)
[xshift=8.6cm] (0.4,-0.25) rectangle (8.2,4.25)
[xshift=8.6cm] (0.4,-0.25) rectangle (8.2,4.25);
\fill[black] (0,0) circle (3pt)
[xshift=8.6cm] (0,0) circle (3pt)
[xshift=8.6cm] (0,0) circle (3pt)
[xshift=8.6cm] (0,0) circle (3pt)
[xshift=8.6cm] (0,0) circle (3pt);
\fill[black] (0,1) circle (3pt)
[xshift=8.6cm] (0,1) circle (3pt)
[xshift=8.6cm] (0,1) circle (3pt)
[xshift=8.6cm] (0,1) circle (3pt)
[xshift=8.6cm] (0,1) circle (3pt);
\fill[black] (0,4) circle (3pt)
[xshift=8.6cm] (0,4) circle (3pt)
[xshift=8.6cm] (0,4) circle (3pt)
[xshift=8.6cm] (0,4) circle (3pt)
[xshift=8.6cm] (0,4) circle (3pt);
\draw (0,0)--(0,1)
[xshift=8.6cm] (0,0)--(0,1)
[xshift=8.6cm] (0,0)--(0,1)
[xshift=8.6cm] (0,0)--(0,1)
[xshift=8.6cm] (0,0)--(0,1);
\draw [densely dashed](0,1)--(0,4)
[xshift=8.6cm] (0,1)--(0,4)
[xshift=8.6cm] (0,1)--(0,4)
[xshift=8.6cm] (0,0)--(0,4)
[xshift=8.6cm] (0,0)--(0,4);

\fill[black] (0.6,3.6) circle (3pt)
[xshift=1.1cm] (0.6,3.6) circle (3pt)
[xshift=1.15cm] (0.6,3.6) circle (3pt)
[xshift=1.1cm] (0.6,3.6) circle (3pt)
[xshift=1.15cm] (0.6,3.6) circle (3pt)
[xshift=1.1cm] (0.6,3.6) circle (3pt)
[xshift=3cm] (0.6,3.6) circle (3pt)
[xshift=1.1cm] (0.6,3.6) circle (3pt)
[xshift=1.15cm] (0.6,3.6) circle (3pt)
[xshift=1.1cm] (0.6,3.6) circle (3pt)
[xshift=1.15cm] (0.6,3.6) circle (3pt)
[xshift=1.1cm] (0.6,3.6) circle (3pt)
[xshift=3cm] (0.6,3.6) circle (3pt)
[xshift=1.1cm] (0.6,3.6) circle (3pt)
[xshift=1.15cm] (0.6,3.6) circle (3pt)
[xshift=1.1cm] (0.6,3.6) circle (3pt)
[xshift=1.15cm] (0.6,3.6) circle (3pt)
[xshift=1.1cm] (0.6,3.6) circle (3pt)
[xshift=3cm] (0.6,3.6) circle (3pt)
[xshift=1.1cm] (0.6,3.6) circle (3pt)
[xshift=1.15cm] (0.6,3.6) circle (3pt)
[xshift=1.1cm] (0.6,3.6) circle (3pt)
[xshift=1.15cm] (0.6,3.6) circle (3pt)
[xshift=1.1cm] (0.6,3.6) circle (3pt)
[xshift=3cm] (0.6,3.6) circle (3pt)
[xshift=1.1cm] (0.6,3.6) circle (3pt)
[xshift=1.15cm] (0.6,3.6) circle (3pt)
[xshift=1.1cm] (0.6,3.6) circle (3pt)
[xshift=1.15cm] (0.6,3.6) circle (3pt)
[xshift=1.1cm] (0.6,3.6) circle (3pt);

\fill[black] (0.8,1.8) circle (3pt)
[xshift=1.1cm](0.8,1.8) circle (3pt)
[xshift=1.15cm](0.8,1.8) circle (3pt)
[xshift=1.1cm](0.8,1.8) circle (3pt)
[xshift=1.15cm](0.8,1.8) circle (3pt)
[xshift=1.1cm](0.8,1.8) circle (3pt)
[xshift=3cm] (0.8,1.8) circle (3pt)
[xshift=1.1cm] (0.8,1.8) circle (3pt)
[xshift=1.15cm](0.8,1.8) circle (3pt)
[xshift=1.1cm](0.8,1.8) circle (3pt)
[xshift=1.15cm] (0.8,1.8) circle (3pt)
[xshift=1.1cm] (0.8,1.8) circle (3pt)
[xshift=3cm] (0.8,1.8) circle (3pt)
[xshift=1.1cm] (0.8,1.8) circle (3pt)
[xshift=1.15cm](0.8,1.8) circle (3pt)
[xshift=1.1cm](0.8,1.8) circle (3pt)
[xshift=1.15cm] (0.8,1.8) circle (3pt)
[xshift=1.1cm] (0.8,1.8) circle (3pt)
[xshift=3cm] (0.8,1.8) circle (3pt)
[xshift=1.1cm] (0.8,1.8) circle (3pt)
[xshift=1.15cm](0.8,1.8) circle (3pt)
[xshift=1.1cm](0.8,1.8) circle (3pt)
[xshift=1.15cm] (0.8,1.8) circle (3pt)
[xshift=1.1cm] (0.8,1.8) circle (3pt)
[xshift=3cm] (0.8,1.8) circle (3pt)
[xshift=1.1cm] (0.8,1.8) circle (3pt)
[xshift=1.15cm](0.8,1.8) circle (3pt)
[xshift=1.1cm](0.8,1.8) circle (3pt)
[xshift=1.15cm] (0.8,1.8) circle (3pt)
[xshift=1.1cm] (0.8,1.8) circle (3pt);

\fill[black] (1,3.6) circle (3pt)
[xshift=1.1cm](1,3.6) circle (3pt)
[xshift=1.15cm](1,3.6) circle (3pt)
[xshift=1.1cm](1,3.6) circle (3pt)
[xshift=1.15cm](1,3.6) circle (3pt)
[xshift=1.1cm] (1,3.6) circle (3pt)
[xshift=3cm](1,3.6) circle (3pt)
[xshift=1.1cm](1,3.6) circle (3pt)
[xshift=1.15cm](1,3.6) circle (3pt)
[xshift=1.1cm](1,3.6) circle (3pt)
[xshift=1.15cm] (1,3.6) circle (3pt)
[xshift=1.1cm] (1,3.6) circle (3pt)
[xshift=3cm](1,3.6) circle (3pt)
[xshift=1.1cm](1,3.6) circle (3pt)
[xshift=1.15cm](1,3.6) circle (3pt)
[xshift=1.1cm](1,3.6) circle (3pt)
[xshift=1.15cm] (1,3.6) circle (3pt)
[xshift=1.1cm] (1,3.6) circle (3pt)
[xshift=3cm](1,3.6) circle (3pt)
[xshift=1.1cm](1,3.6) circle (3pt)
[xshift=1.15cm](1,3.6) circle (3pt)
[xshift=1.1cm](1,3.6) circle (3pt)
[xshift=1.15cm] (1,3.6) circle (3pt)
[xshift=1.1cm] (1,3.6) circle (3pt)
[xshift=3cm](1,3.6) circle (3pt)
[xshift=1.1cm](1,3.6) circle (3pt)
[xshift=1.15cm](1,3.6) circle (3pt)
[xshift=1.1cm](1,3.6) circle (3pt)
[xshift=1.15cm] (1,3.6) circle (3pt)
[xshift=1.1cm] (1,3.6) circle (3pt);
\draw (0.6,3.6)--(0.8,1.8)
[xshift=1.1cm](0.6,3.6)--(0.8,1.8)
[xshift=1.15cm](0.6,3.6)--(0.8,1.8)
[xshift=1.1cm](0.6,3.6)--(0.8,1.8)
[xshift=1.15cm](0.6,3.6)--(0.8,1.8)
[xshift=1.1cm](0.6,3.6)--(0.8,1.8)
[xshift=3cm] (0.6,3.6)--(0.8,1.8)
[xshift=1.1cm](0.6,3.6)--(0.8,1.8)
[xshift=1.15cm](0.6,3.6)--(0.8,1.8)
[xshift=1.1cm](0.6,3.6)--(0.8,1.8)
[xshift=1.15cm] (0.6,3.6)--(0.8,1.8)
[xshift=1.1cm] (0.6,3.6)--(0.8,1.8)
[xshift=3cm] (0.6,3.6)--(0.8,1.8)
[xshift=1.1cm](0.6,3.6)--(0.8,1.8)
[xshift=1.15cm](0.6,3.6)--(0.8,1.8)
[xshift=1.1cm](0.6,3.6)--(0.8,1.8)
[xshift=1.15cm] (0.6,3.6)--(0.8,1.8)
[xshift=1.1cm] (0.6,3.6)--(0.8,1.8)
[xshift=3cm] (0.6,3.6)--(0.8,1.8)
[xshift=1.1cm](0.6,3.6)--(0.8,1.8)
[xshift=1.15cm](0.6,3.6)--(0.8,1.8)
[xshift=1.1cm](0.6,3.6)--(0.8,1.8)
[xshift=1.15cm] (0.6,3.6)--(0.8,1.8)
[xshift=1.1cm] (0.6,3.6)--(0.8,1.8)
[xshift=3cm] (0.6,3.6)--(0.8,1.8)
[xshift=1.1cm](0.6,3.6)--(0.8,1.8)
[xshift=1.15cm](0.6,3.6)--(0.8,1.8)
[xshift=1.1cm](0.6,3.6)--(0.8,1.8)
[xshift=1.15cm] (0.6,3.6)--(0.8,1.8)
[xshift=1.1cm] (0.6,3.6)--(0.8,1.8);

\draw (0.8,1.8)--(1,3.6)
[xshift=1.1cm](0.8,1.8)--(1,3.6)
[xshift=1.15cm](0.8,1.8)--(1,3.6)
[xshift=1.1cm](0.8,1.8)--(1,3.6)
[xshift=1.15cm](0.8,1.8)--(1,3.6)
[xshift=1.1cm](0.8,1.8)--(1,3.6)
[xshift=3cm] (0.8,1.8)--(1,3.6)
[xshift=1.1cm](0.8,1.8)--(1,3.6)
[xshift=1.15cm](0.8,1.8)--(1,3.6)
[xshift=1.1cm](0.8,1.8)--(1,3.6)
[xshift=1.15cm] (0.8,1.8)--(1,3.6)
[xshift=1.1cm] (0.8,1.8)--(1,3.6)
[xshift=3cm] (0.8,1.8)--(1,3.6)
[xshift=1.1cm](0.8,1.8)--(1,3.6)
[xshift=1.15cm](0.8,1.8)--(1,3.6)
[xshift=1.1cm](0.8,1.8)--(1,3.6)
[xshift=1.15cm] (0.8,1.8)--(1,3.6)
[xshift=1.1cm] (0.8,1.8)--(1,3.6)
[xshift=3cm] (0.8,1.8)--(1,3.6)
[xshift=1.1cm](0.8,1.8)--(1,3.6)
[xshift=1.15cm](0.8,1.8)--(1,3.6)
[xshift=1.1cm](0.8,1.8)--(1,3.6)
[xshift=1.15cm] (0.8,1.8)--(1,3.6)
[xshift=1.1cm] (0.8,1.8)--(1,3.6)
[xshift=3cm] (0.8,1.8)--(1,3.6)
[xshift=1.1cm](0.8,1.8)--(1,3.6)
[xshift=1.15cm](0.8,1.8)--(1,3.6)
[xshift=1.1cm](0.8,1.8)--(1,3.6)
[xshift=1.15cm] (0.8,1.8)--(1,3.6)
[xshift=1.1cm] (0.8,1.8)--(1,3.6);
\draw [densely dotted](1.15,3.6)--(1.5,3.6)
[xshift=1.1cm] (1.15,3.6)--(1.5,3.6)
[yshift=-1.8cm] (1.15,3.6)--(1.5,3.6)
[shift={(1.15cm, 1.8cm)}](1.15,3.6)--(1.5,3.6)
[xshift=1.1cm] (1.15,3.6)--(1.5,3.6)
[yshift=-1.8cm] (1.15,3.6)--(1.5,3.6)
[shift={(1.15cm,1.8cm)}] (1.15,3.6)--(1.5,3.6)
[xshift=1.1cm](1.15,3.6)--(1.5,3.6)
[yshift=-1.8cm] (1.15,3.6)--(1.5,3.6)
[shift={(3cm,1.8cm)}] (1.15,3.6)--(1.5,3.6)
[xshift=1.1cm] (1.15,3.6)--(1.5,3.6)
[yshift=-1.8cm] (1.15,3.6)--(1.5,3.6)
[shift={(1.15cm, 1.8cm)}](1.15,3.6)--(1.5,3.6)
[xshift=1.1cm] (1.15,3.6)--(1.5,3.6)
[yshift=-1.8cm] (1.15,3.6)--(1.5,3.6)
[shift={(1.15cm,1.8cm)}] (1.15,3.6)--(1.5,3.6)
[xshift=1.1cm](1.15,3.6)--(1.5,3.6)
[yshift=-1.8cm] (1.15,3.6)--(1.5,3.6)
[shift={(3cm,1.8cm)}] (1.15,3.6)--(1.5,3.6)
[xshift=1.1cm] (1.15,3.6)--(1.5,3.6)
[yshift=-1.8cm] (1.15,3.6)--(1.5,3.6)
[shift={(1.15cm, 1.8cm)}](1.15,3.6)--(1.5,3.6)
[xshift=1.1cm] (1.15,3.6)--(1.5,3.6)
[yshift=-1.8cm] (1.15,3.6)--(1.5,3.6)
[shift={(1.15cm,1.8cm)}] (1.15,3.6)--(1.5,3.6)
[xshift=1.1cm](1.15,3.6)--(1.5,3.6)
[yshift=-1.8cm] (1.15,3.6)--(1.5,3.6)
[shift={(3cm,1.8cm)}] (1.15,3.6)--(1.5,3.6)
[xshift=1.1cm] (1.15,3.6)--(1.5,3.6)
[yshift=-1.8cm] (1.15,3.6)--(1.5,3.6)
[shift={(1.15cm, 1.8cm)}](1.15,3.6)--(1.5,3.6)
[xshift=1.1cm] (1.15,3.6)--(1.5,3.6)
[yshift=-1.8cm] (1.15,3.6)--(1.5,3.6)
[shift={(1.15cm,1.8cm)}] (1.15,3.6)--(1.5,3.6)
[xshift=1.1cm](1.15,3.6)--(1.5,3.6)
[yshift=-1.8cm] (1.15,3.6)--(1.5,3.6)
[shift={(3cm,1.8cm)}] (1.15,3.6)--(1.5,3.6)
[xshift=1.1cm] (1.15,3.6)--(1.5,3.6)
[yshift=-1.8cm] (1.15,3.6)--(1.5,3.6)
[shift={(1.15cm, 1.8cm)}](1.15,3.6)--(1.5,3.6)
[xshift=1.1cm] (1.15,3.6)--(1.5,3.6)
[yshift=-1.8cm] (1.15,3.6)--(1.5,3.6)
[shift={(1.15cm,1.8cm)}] (1.15,3.6)--(1.5,3.6)
[xshift=1.1cm](1.15,3.6)--(1.5,3.6)
[yshift=-1.8cm] (1.15,3.6)--(1.5,3.6);
\fill[black] (1.35,0) circle (3pt)
[xshift=2.25cm](1.35,0) circle (3pt)
[xshift=2.25cm] (1.35,0) circle (3pt)
[xshift=4.1cm](1.35,0) circle (3pt)
[xshift=2.25cm] (1.35,0) circle (3pt)
[xshift=2.25cm](1.35,0) circle (3pt)
[xshift=4.1cm] (1.35,0) circle (3pt)
[xshift=2.25cm] (1.35,0) circle (3pt)
[xshift=2.25cm](1.35,0) circle (3pt)
[xshift=4.1cm](1.35,0) circle (3pt)
[xshift=2.25cm] (1.35,0) circle (3pt)
[xshift=2.25cm](1.35,0) circle (3pt)
[xshift=4.1cm](1.35,0) circle (3pt)
[xshift=2.25cm] (1.35,0) circle (3pt)
[xshift=2.25cm](1.35,0) circle (3pt);

\draw (1.35,0)--(1.9,1.8)
[xshift=2.25cm] (1.35,0)--(1.9,1.8)
[xshift=2.25cm] (1.35,0)--(1.9,1.8)
[xshift=4.1cm] (1.35,0)--(1.9,1.8)
[xshift=2.25cm] (1.35,0)--(1.9,1.8)
[xshift=2.25cm] (1.35,0)--(1.9,1.8)
[xshift=4.1cm] (1.35,0)--(1.9,1.8)
[xshift=2.25cm] (1.35,0)--(1.9,1.8)
[xshift=2.25cm] (1.35,0)--(1.9,1.8)
[xshift=4.1cm] (1.35,0)--(1.9,1.8)
[xshift=2.25cm] (1.35,0)--(1.9,1.8)
[xshift=2.25cm] (1.35,0)--(1.9,1.8)
[xshift=4.1cm] (1.35,0)--(1.9,1.8)
[xshift=2.25cm] (1.35,0)--(1.9,1.8)
[xshift=2.25cm] (1.35,0)--(1.9,1.8);

\draw (0.8,1.8)--(1.35,0)
[xshift=2.25cm]  (0.8,1.8)--(1.35,0)
[xshift=2.25cm]  (0.8,1.8)--(1.35,0)
[xshift=4.1cm]  (0.8,1.8)--(1.35,0)
[xshift=2.25cm]  (0.8,1.8)--(1.35,0)
[xshift=2.25cm]  (0.8,1.8)--(1.35,0)
[xshift=4.1cm]  (0.8,1.8)--(1.35,0)
[xshift=2.25cm]  (0.8,1.8)--(1.35,0)
[xshift=2.25cm]  (0.8,1.8)--(1.35,0)
[xshift=4.1cm]  (0.8,1.8)--(1.35,0)
[xshift=2.25cm]  (0.8,1.8)--(1.35,0)
[xshift=2.25cm]  (0.8,1.8)--(1.35,0)
[xshift=4.1cm]  (0.8,1.8)--(1.35,0)
[xshift=2.25cm]  (0.8,1.8)--(1.35,0)
[xshift=2.25cm]  (0.8,1.8)--(1.35,0);
\draw [densely dotted](7.3,3.6)--(8.1,3.6)
[yshift=-1.8cm] (7.3,3.6)--(8.1,3.6)
[yshift=-1.8cm] (7.3,3.6)--(8.1,3.6)
[xshift=8.6cm] (7.3,3.6)--(8.1,3.6)
[yshift=1.8cm] (7.3,3.6)--(8.1,3.6)
[yshift=1.8cm] (7.3,3.6)--(8.1,3.6)
[xshift=8.6cm] (7.3,3.6)--(8.1,3.6)
[yshift=-1.8cm] (7.3,3.6)--(8.1,3.6)
[yshift=-1.8cm] (7.3,3.6)--(8.1,3.6)
[xshift=8.6cm] (7.3,3.6)--(8.1,3.6)
[yshift=1.8cm] (7.3,3.6)--(8.1,3.6)
[yshift=1.8cm] (7.3,3.6)--(8.1,3.6)
[xshift=8.6cm] (7.3,3.6)--(8.1,3.6)
[yshift=-1.8cm] (7.3,3.6)--(8.1,3.6)
[yshift=-1.8cm] (7.3,3.6)--(8.1,3.6);
%
\fill[black] (2,9) circle (3pt)
[xshift=8.6cm] (2,9) circle (3pt)
[xshift=8.6cm] (2,9) circle (3pt)
[xshift=8.6cm] (2,9) circle (3pt)
[xshift=8.6cm] (2,9) circle (3pt);
\draw[densely dashed] (2,9)--(1.35,4)
[xshift=8.6cm] (2,9)--(1.35,4)
[xshift=8.6cm] (2,9)--(1.35,4)
[xshift=8.6cm] ((2,9)--(1.35,4)
[xshift=8.6cm] (2,9)--(1.35,4);
\draw[densely dashed] (0,4)--(2,9)
[xshift=8.6cm] (0,4)--(2,9)
[xshift=8.6cm] (0,4)--(2,9)
[xshift=8.6cm] (0,4)--(2,9)
[xshift=8.6cm] (0,4)--(2,9);
\draw[densely dashed] (2,9)--(6,4)
[xshift=8.6cm] (2,9)--(6,4)
[xshift=8.6cm] (2,9)--(6,4)
[xshift=8.6cm] (2,9)--(6,4)
[xshift=8.6cm] (2,9)--(6,4);

\draw[densely dashed] (3.6,4)--(2,9)
[xshift=8.6cm] (3.6,4)--(2,9)
[xshift=8.6cm] (3.6,4)--(2,9)
[xshift=8.6cm] (3.6,4)--(2,9)
[xshift=8.6cm] (3.6,4)--(2,9);
\draw [densely dotted](4.8,0)--(5.3,0)
[xshift=8.6cm] (4.8,0)--(5.3,0);
\draw [densely dotted](1.15,1.8)--(1.55,1.8);
\draw[red] (17.2,4)--(27.8,9);
\draw[red] (25.8,4)--(36.4,9);
\node (top1)[above] at(2,9){\footnotesize$(a,\top)$};
\node (top2)[above] at(10.6,9){\footnotesize$(b,\top)$};
\node (top3)[above] at(19.2,9){\footnotesize$(f_{a,b}(1),\top)$};
\node (top4)[above] at(27.8,9){\footnotesize$(f_{a,b}(1.x),\top)$};
\node (top5)[above] at(36.4,9){\footnotesize$(f_{a,b}(1.x.y),\top)$};
\node (d1) at(4,-1.3){\footnotesize$a$};
\node (d2) at(13,-1.3){\footnotesize$b$};
\node (d3) at(21.5,-1.3){\footnotesize$f_{a,b}(1)$};
\node (d3) at(30,-1.3){\footnotesize$f_{a,b}(1.x)$};
\node (d3) at(38.5,-1.3){\footnotesize$f_{a,b}(1.x.y)$};
\node (l2) [below] at(17.2,4){\tiny$(f_{a,b}(1),x)$};
\node (l3) [below] at(25.6,4){\tiny$(f_{a,b}(1.x),y)$};
\node (l4) [below] at(1.35,0){\tiny$(a,1)$};
\node (l5) [below] at(1.9,1.8){\tiny$(a,1.x)$};
\node (l6) [below] at(8.6,0){\tiny$(b,1)$};

\end{tikzpicture}

  \caption{ The strict order $<_{4}$}\label{F.5}
\end{figure}
The red lines in Fig.\ref{F.5} illustrate the cases: $(f_{a,b}(1),x)<_{4}(f_{a,b}(1.x),\top)$ and $(f_{a,b}(1.x),y)<_{4}(f_{a,b}(1.x.y),\top)$.

Based on the above observations, $<_{4}$ is defined after $<_{3}$ and is all linked together. Specifically, for given $(a,b)\in \mathbb{N}\times \mathbb{N}$ with $a<b$, we have:

$(b,y)<_{3}(f_{a,b}(y),\top)$ for any $y\in \mathbb{N}$;

$(f_{a,b}(y),z)<_{4}(f_{a,b}(y.z),\top)$ for any $z\in \mathbb{N}$;

$(f_{a,b}(y.z),u)<_{4}(f_{a,b}(y.z.u),\top)$ for any $u\in \mathbb{N}$.

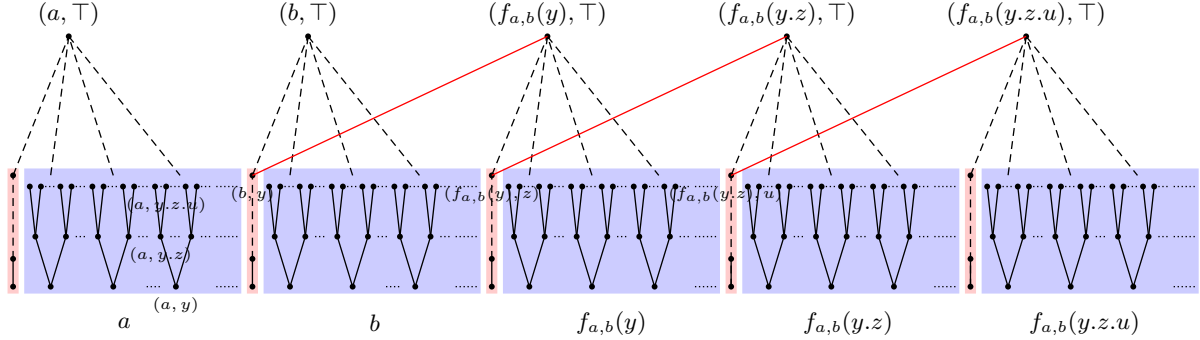
\begin{figure}[H]
  \centering
  \begin{tikzpicture} [line width=0.5pt,scale=0.37]
\fill[red,opacity=0.2] (-0.2,-0.25) rectangle (0.2,4.25)
[xshift=8.6cm] (-0.2,-0.25) rectangle (0.2,4.25)
[xshift=8.6cm] (-0.2,-0.25) rectangle (0.2,4.25)
[xshift=8.6cm] (-0.2,-0.25) rectangle (0.2,4.25)
[xshift=8.6cm] (-0.2,-0.25) rectangle (0.2,4.25);
\fill[blue,opacity=0.2] (0.4,-0.25) rectangle (8.2,4.25)
[xshift=8.6cm] (0.4,-0.25) rectangle (8.2,4.25)
[xshift=8.6cm] (0.4,-0.25) rectangle (8.2,4.25)
[xshift=8.6cm] (0.4,-0.25) rectangle (8.2,4.25)
[xshift=8.6cm] (0.4,-0.25) rectangle (8.2,4.25);
\fill[black] (0,0) circle (3pt)
[xshift=8.6cm] (0,0) circle (3pt)
[xshift=8.6cm] (0,0) circle (3pt)
[xshift=8.6cm] (0,0) circle (3pt)
[xshift=8.6cm] (0,0) circle (3pt);
\fill[black] (0,1) circle (3pt)
[xshift=8.6cm] (0,1) circle (3pt)
[xshift=8.6cm] (0,1) circle (3pt)
[xshift=8.6cm] (0,1) circle (3pt)
[xshift=8.6cm] (0,1) circle (3pt);
\fill[black] (0,4) circle (3pt)
[xshift=8.6cm] (0,4) circle (3pt)
[xshift=8.6cm] (0,4) circle (3pt)
[xshift=8.6cm] (0,4) circle (3pt)
[xshift=8.6cm] (0,4) circle (3pt);
\draw (0,0)--(0,1)
[xshift=8.6cm] (0,0)--(0,1)
[xshift=8.6cm] (0,0)--(0,1)
[xshift=8.6cm] (0,0)--(0,1)
[xshift=8.6cm] (0,0)--(0,1);
\draw [densely dashed](0,1)--(0,4)
[xshift=8.6cm] (0,1)--(0,4)
[xshift=8.6cm] (0,1)--(0,4)
[xshift=8.6cm] (0,0)--(0,4)
[xshift=8.6cm] (0,0)--(0,4);

\fill[black] (0.6,3.6) circle (3pt)
[xshift=1.1cm] (0.6,3.6) circle (3pt)
[xshift=1.15cm] (0.6,3.6) circle (3pt)
[xshift=1.1cm] (0.6,3.6) circle (3pt)
[xshift=1.15cm] (0.6,3.6) circle (3pt)
[xshift=1.1cm] (0.6,3.6) circle (3pt)
[xshift=3cm] (0.6,3.6) circle (3pt)
[xshift=1.1cm] (0.6,3.6) circle (3pt)
[xshift=1.15cm] (0.6,3.6) circle (3pt)
[xshift=1.1cm] (0.6,3.6) circle (3pt)
[xshift=1.15cm] (0.6,3.6) circle (3pt)
[xshift=1.1cm] (0.6,3.6) circle (3pt)
[xshift=3cm] (0.6,3.6) circle (3pt)
[xshift=1.1cm] (0.6,3.6) circle (3pt)
[xshift=1.15cm] (0.6,3.6) circle (3pt)
[xshift=1.1cm] (0.6,3.6) circle (3pt)
[xshift=1.15cm] (0.6,3.6) circle (3pt)
[xshift=1.1cm] (0.6,3.6) circle (3pt)
[xshift=3cm] (0.6,3.6) circle (3pt)
[xshift=1.1cm] (0.6,3.6) circle (3pt)
[xshift=1.15cm] (0.6,3.6) circle (3pt)
[xshift=1.1cm] (0.6,3.6) circle (3pt)
[xshift=1.15cm] (0.6,3.6) circle (3pt)
[xshift=1.1cm] (0.6,3.6) circle (3pt)
[xshift=3cm] (0.6,3.6) circle (3pt)
[xshift=1.1cm] (0.6,3.6) circle (3pt)
[xshift=1.15cm] (0.6,3.6) circle (3pt)
[xshift=1.1cm] (0.6,3.6) circle (3pt)
[xshift=1.15cm] (0.6,3.6) circle (3pt)
[xshift=1.1cm] (0.6,3.6) circle (3pt);

\fill[black] (0.8,1.8) circle (3pt)
[xshift=1.1cm](0.8,1.8) circle (3pt)
[xshift=1.15cm](0.8,1.8) circle (3pt)
[xshift=1.1cm](0.8,1.8) circle (3pt)
[xshift=1.15cm](0.8,1.8) circle (3pt)
[xshift=1.1cm](0.8,1.8) circle (3pt)
[xshift=3cm] (0.8,1.8) circle (3pt)
[xshift=1.1cm] (0.8,1.8) circle (3pt)
[xshift=1.15cm](0.8,1.8) circle (3pt)
[xshift=1.1cm](0.8,1.8) circle (3pt)
[xshift=1.15cm] (0.8,1.8) circle (3pt)
[xshift=1.1cm] (0.8,1.8) circle (3pt)
[xshift=3cm] (0.8,1.8) circle (3pt)
[xshift=1.1cm] (0.8,1.8) circle (3pt)
[xshift=1.15cm](0.8,1.8) circle (3pt)
[xshift=1.1cm](0.8,1.8) circle (3pt)
[xshift=1.15cm] (0.8,1.8) circle (3pt)
[xshift=1.1cm] (0.8,1.8) circle (3pt)
[xshift=3cm] (0.8,1.8) circle (3pt)
[xshift=1.1cm] (0.8,1.8) circle (3pt)
[xshift=1.15cm](0.8,1.8) circle (3pt)
[xshift=1.1cm](0.8,1.8) circle (3pt)
[xshift=1.15cm] (0.8,1.8) circle (3pt)
[xshift=1.1cm] (0.8,1.8) circle (3pt)
[xshift=3cm] (0.8,1.8) circle (3pt)
[xshift=1.1cm] (0.8,1.8) circle (3pt)
[xshift=1.15cm](0.8,1.8) circle (3pt)
[xshift=1.1cm](0.8,1.8) circle (3pt)
[xshift=1.15cm] (0.8,1.8) circle (3pt)
[xshift=1.1cm] (0.8,1.8) circle (3pt);

\fill[black] (1,3.6) circle (3pt)
[xshift=1.1cm](1,3.6) circle (3pt)
[xshift=1.15cm](1,3.6) circle (3pt)
[xshift=1.1cm](1,3.6) circle (3pt)
[xshift=1.15cm](1,3.6) circle (3pt)
[xshift=1.1cm] (1,3.6) circle (3pt)
[xshift=3cm](1,3.6) circle (3pt)
[xshift=1.1cm](1,3.6) circle (3pt)
[xshift=1.15cm](1,3.6) circle (3pt)
[xshift=1.1cm](1,3.6) circle (3pt)
[xshift=1.15cm] (1,3.6) circle (3pt)
[xshift=1.1cm] (1,3.6) circle (3pt)
[xshift=3cm](1,3.6) circle (3pt)
[xshift=1.1cm](1,3.6) circle (3pt)
[xshift=1.15cm](1,3.6) circle (3pt)
[xshift=1.1cm](1,3.6) circle (3pt)
[xshift=1.15cm] (1,3.6) circle (3pt)
[xshift=1.1cm] (1,3.6) circle (3pt)
[xshift=3cm](1,3.6) circle (3pt)
[xshift=1.1cm](1,3.6) circle (3pt)
[xshift=1.15cm](1,3.6) circle (3pt)
[xshift=1.1cm](1,3.6) circle (3pt)
[xshift=1.15cm] (1,3.6) circle (3pt)
[xshift=1.1cm] (1,3.6) circle (3pt)
[xshift=3cm](1,3.6) circle (3pt)
[xshift=1.1cm](1,3.6) circle (3pt)
[xshift=1.15cm](1,3.6) circle (3pt)
[xshift=1.1cm](1,3.6) circle (3pt)
[xshift=1.15cm] (1,3.6) circle (3pt)
[xshift=1.1cm] (1,3.6) circle (3pt);
\draw (0.6,3.6)--(0.8,1.8)
[xshift=1.1cm](0.6,3.6)--(0.8,1.8)
[xshift=1.15cm](0.6,3.6)--(0.8,1.8)
[xshift=1.1cm](0.6,3.6)--(0.8,1.8)
[xshift=1.15cm](0.6,3.6)--(0.8,1.8)
[xshift=1.1cm](0.6,3.6)--(0.8,1.8)
[xshift=3cm] (0.6,3.6)--(0.8,1.8)
[xshift=1.1cm](0.6,3.6)--(0.8,1.8)
[xshift=1.15cm](0.6,3.6)--(0.8,1.8)
[xshift=1.1cm](0.6,3.6)--(0.8,1.8)
[xshift=1.15cm] (0.6,3.6)--(0.8,1.8)
[xshift=1.1cm] (0.6,3.6)--(0.8,1.8)
[xshift=3cm] (0.6,3.6)--(0.8,1.8)
[xshift=1.1cm](0.6,3.6)--(0.8,1.8)
[xshift=1.15cm](0.6,3.6)--(0.8,1.8)
[xshift=1.1cm](0.6,3.6)--(0.8,1.8)
[xshift=1.15cm] (0.6,3.6)--(0.8,1.8)
[xshift=1.1cm] (0.6,3.6)--(0.8,1.8)
[xshift=3cm] (0.6,3.6)--(0.8,1.8)
[xshift=1.1cm](0.6,3.6)--(0.8,1.8)
[xshift=1.15cm](0.6,3.6)--(0.8,1.8)
[xshift=1.1cm](0.6,3.6)--(0.8,1.8)
[xshift=1.15cm] (0.6,3.6)--(0.8,1.8)
[xshift=1.1cm] (0.6,3.6)--(0.8,1.8)
[xshift=3cm] (0.6,3.6)--(0.8,1.8)
[xshift=1.1cm](0.6,3.6)--(0.8,1.8)
[xshift=1.15cm](0.6,3.6)--(0.8,1.8)
[xshift=1.1cm](0.6,3.6)--(0.8,1.8)
[xshift=1.15cm] (0.6,3.6)--(0.8,1.8)
[xshift=1.1cm] (0.6,3.6)--(0.8,1.8);

\draw (0.8,1.8)--(1,3.6)
[xshift=1.1cm](0.8,1.8)--(1,3.6)
[xshift=1.15cm](0.8,1.8)--(1,3.6)
[xshift=1.1cm](0.8,1.8)--(1,3.6)
[xshift=1.15cm](0.8,1.8)--(1,3.6)
[xshift=1.1cm](0.8,1.8)--(1,3.6)
[xshift=3cm] (0.8,1.8)--(1,3.6)
[xshift=1.1cm](0.8,1.8)--(1,3.6)
[xshift=1.15cm](0.8,1.8)--(1,3.6)
[xshift=1.1cm](0.8,1.8)--(1,3.6)
[xshift=1.15cm] (0.8,1.8)--(1,3.6)
[xshift=1.1cm] (0.8,1.8)--(1,3.6)
[xshift=3cm] (0.8,1.8)--(1,3.6)
[xshift=1.1cm](0.8,1.8)--(1,3.6)
[xshift=1.15cm](0.8,1.8)--(1,3.6)
[xshift=1.1cm](0.8,1.8)--(1,3.6)
[xshift=1.15cm] (0.8,1.8)--(1,3.6)
[xshift=1.1cm] (0.8,1.8)--(1,3.6)
[xshift=3cm] (0.8,1.8)--(1,3.6)
[xshift=1.1cm](0.8,1.8)--(1,3.6)
[xshift=1.15cm](0.8,1.8)--(1,3.6)
[xshift=1.1cm](0.8,1.8)--(1,3.6)
[xshift=1.15cm] (0.8,1.8)--(1,3.6)
[xshift=1.1cm] (0.8,1.8)--(1,3.6)
[xshift=3cm] (0.8,1.8)--(1,3.6)
[xshift=1.1cm](0.8,1.8)--(1,3.6)
[xshift=1.15cm](0.8,1.8)--(1,3.6)
[xshift=1.1cm](0.8,1.8)--(1,3.6)
[xshift=1.15cm] (0.8,1.8)--(1,3.6)
[xshift=1.1cm] (0.8,1.8)--(1,3.6);
\draw [densely dotted](1.15,3.6)--(1.5,3.6)
[xshift=1.1cm] (1.15,3.6)--(1.5,3.6)
[yshift=-1.8cm] (1.15,3.6)--(1.5,3.6)
[shift={(1.15cm, 1.8cm)}](1.15,3.6)--(1.5,3.6)
[xshift=1.1cm] (1.15,3.6)--(1.5,3.6)
[yshift=-1.8cm] (1.15,3.6)--(1.5,3.6)
[shift={(1.15cm,1.8cm)}] (1.15,3.6)--(1.5,3.6)
[xshift=1.1cm](1.15,3.6)--(1.5,3.6)
[yshift=-1.8cm] (1.15,3.6)--(1.5,3.6)
[shift={(3cm,1.8cm)}] (1.15,3.6)--(1.5,3.6)
[xshift=1.1cm] (1.15,3.6)--(1.5,3.6)
[yshift=-1.8cm] (1.15,3.6)--(1.5,3.6)
[shift={(1.15cm, 1.8cm)}](1.15,3.6)--(1.5,3.6)
[xshift=1.1cm] (1.15,3.6)--(1.5,3.6)
[yshift=-1.8cm] (1.15,3.6)--(1.5,3.6)
[shift={(1.15cm,1.8cm)}] (1.15,3.6)--(1.5,3.6)
[xshift=1.1cm](1.15,3.6)--(1.5,3.6)
[yshift=-1.8cm] (1.15,3.6)--(1.5,3.6)
[shift={(3cm,1.8cm)}] (1.15,3.6)--(1.5,3.6)
[xshift=1.1cm] (1.15,3.6)--(1.5,3.6)
[yshift=-1.8cm] (1.15,3.6)--(1.5,3.6)
[shift={(1.15cm, 1.8cm)}](1.15,3.6)--(1.5,3.6)
[xshift=1.1cm] (1.15,3.6)--(1.5,3.6)
[yshift=-1.8cm] (1.15,3.6)--(1.5,3.6)
[shift={(1.15cm,1.8cm)}] (1.15,3.6)--(1.5,3.6)
[xshift=1.1cm](1.15,3.6)--(1.5,3.6)
[yshift=-1.8cm] (1.15,3.6)--(1.5,3.6)
[shift={(3cm,1.8cm)}] (1.15,3.6)--(1.5,3.6)
[xshift=1.1cm] (1.15,3.6)--(1.5,3.6)
[yshift=-1.8cm] (1.15,3.6)--(1.5,3.6)
[shift={(1.15cm, 1.8cm)}](1.15,3.6)--(1.5,3.6)
[xshift=1.1cm] (1.15,3.6)--(1.5,3.6)
[yshift=-1.8cm] (1.15,3.6)--(1.5,3.6)
[shift={(1.15cm,1.8cm)}] (1.15,3.6)--(1.5,3.6)
[xshift=1.1cm](1.15,3.6)--(1.5,3.6)
[yshift=-1.8cm] (1.15,3.6)--(1.5,3.6)
[shift={(3cm,1.8cm)}] (1.15,3.6)--(1.5,3.6)
[xshift=1.1cm] (1.15,3.6)--(1.5,3.6)
[yshift=-1.8cm] (1.15,3.6)--(1.5,3.6)
[shift={(1.15cm, 1.8cm)}](1.15,3.6)--(1.5,3.6)
[xshift=1.1cm] (1.15,3.6)--(1.5,3.6)
[yshift=-1.8cm] (1.15,3.6)--(1.5,3.6)
[shift={(1.15cm,1.8cm)}] (1.15,3.6)--(1.5,3.6)
[xshift=1.1cm](1.15,3.6)--(1.5,3.6)
[yshift=-1.8cm] (1.15,3.6)--(1.5,3.6);
\fill[black] (1.35,0) circle (3pt)
[xshift=2.25cm](1.35,0) circle (3pt)
[xshift=2.25cm] (1.35,0) circle (3pt)
[xshift=4.1cm](1.35,0) circle (3pt)
[xshift=2.25cm] (1.35,0) circle (3pt)
[xshift=2.25cm](1.35,0) circle (3pt)
[xshift=4.1cm] (1.35,0) circle (3pt)
[xshift=2.25cm] (1.35,0) circle (3pt)
[xshift=2.25cm](1.35,0) circle (3pt)
[xshift=4.1cm](1.35,0) circle (3pt)
[xshift=2.25cm] (1.35,0) circle (3pt)
[xshift=2.25cm](1.35,0) circle (3pt)
[xshift=4.1cm](1.35,0) circle (3pt)
[xshift=2.25cm] (1.35,0) circle (3pt)
[xshift=2.25cm](1.35,0) circle (3pt);

\draw (1.35,0)--(1.9,1.8)
[xshift=2.25cm] (1.35,0)--(1.9,1.8)
[xshift=2.25cm] (1.35,0)--(1.9,1.8)
[xshift=4.1cm] (1.35,0)--(1.9,1.8)
[xshift=2.25cm] (1.35,0)--(1.9,1.8)
[xshift=2.25cm] (1.35,0)--(1.9,1.8)
[xshift=4.1cm] (1.35,0)--(1.9,1.8)
[xshift=2.25cm] (1.35,0)--(1.9,1.8)
[xshift=2.25cm] (1.35,0)--(1.9,1.8)
[xshift=4.1cm] (1.35,0)--(1.9,1.8)
[xshift=2.25cm] (1.35,0)--(1.9,1.8)
[xshift=2.25cm] (1.35,0)--(1.9,1.8)
[xshift=4.1cm] (1.35,0)--(1.9,1.8)
[xshift=2.25cm] (1.35,0)--(1.9,1.8)
[xshift=2.25cm] (1.35,0)--(1.9,1.8);

\draw (0.8,1.8)--(1.35,0)
[xshift=2.25cm]  (0.8,1.8)--(1.35,0)
[xshift=2.25cm]  (0.8,1.8)--(1.35,0)
[xshift=4.1cm]  (0.8,1.8)--(1.35,0)
[xshift=2.25cm]  (0.8,1.8)--(1.35,0)
[xshift=2.25cm]  (0.8,1.8)--(1.35,0)
[xshift=4.1cm]  (0.8,1.8)--(1.35,0)
[xshift=2.25cm]  (0.8,1.8)--(1.35,0)
[xshift=2.25cm]  (0.8,1.8)--(1.35,0)
[xshift=4.1cm]  (0.8,1.8)--(1.35,0)
[xshift=2.25cm]  (0.8,1.8)--(1.35,0)
[xshift=2.25cm]  (0.8,1.8)--(1.35,0)
[xshift=4.1cm]  (0.8,1.8)--(1.35,0)
[xshift=2.25cm]  (0.8,1.8)--(1.35,0)
[xshift=2.25cm]  (0.8,1.8)--(1.35,0);
\draw [densely dotted](7.3,3.6)--(8.1,3.6)
[yshift=-1.8cm] (7.3,3.6)--(8.1,3.6)
[yshift=-1.8cm] (7.3,3.6)--(8.1,3.6)
[xshift=8.6cm] (7.3,3.6)--(8.1,3.6)
[yshift=1.8cm] (7.3,3.6)--(8.1,3.6)
[yshift=1.8cm] (7.3,3.6)--(8.1,3.6)
[xshift=8.6cm] (7.3,3.6)--(8.1,3.6)
[yshift=-1.8cm] (7.3,3.6)--(8.1,3.6)
[yshift=-1.8cm] (7.3,3.6)--(8.1,3.6)
[xshift=8.6cm] (7.3,3.6)--(8.1,3.6)
[yshift=1.8cm] (7.3,3.6)--(8.1,3.6)
[yshift=1.8cm] (7.3,3.6)--(8.1,3.6)
[xshift=8.6cm] (7.3,3.6)--(8.1,3.6)
[yshift=-1.8cm] (7.3,3.6)--(8.1,3.6)
[yshift=-1.8cm] (7.3,3.6)--(8.1,3.6);
%
\fill[black] (2,9) circle (3pt)
[xshift=8.6cm] (2,9) circle (3pt)
[xshift=8.6cm] (2,9) circle (3pt)
[xshift=8.6cm] (2,9) circle (3pt)
[xshift=8.6cm] (2,9) circle (3pt);
\draw[densely dashed] (2,9)--(1.35,4)
[xshift=8.6cm] (2,9)--(1.35,4)
[xshift=8.6cm] (2,9)--(1.35,4)
[xshift=8.6cm] ((2,9)--(1.35,4)
[xshift=8.6cm] (2,9)--(1.35,4);
\draw[densely dashed] (0,4)--(2,9)
[xshift=8.6cm] (0,4)--(2,9)
[xshift=8.6cm] (0,4)--(2,9)
[xshift=8.6cm] (0,4)--(2,9)
[xshift=8.6cm] (0,4)--(2,9);
\draw[densely dashed] (2,9)--(6,4)
[xshift=8.6cm] (2,9)--(6,4)
[xshift=8.6cm] (2,9)--(6,4)
[xshift=8.6cm] (2,9)--(6,4)
[xshift=8.6cm] (2,9)--(6,4);

\draw[densely dashed] (3.6,4)--(2,9)
[xshift=8.6cm] (3.6,4)--(2,9)
[xshift=8.6cm] (3.6,4)--(2,9)
[xshift=8.6cm] (3.6,4)--(2,9)
[xshift=8.6cm] (3.6,4)--(2,9);
\draw [densely dotted](4.8,0)--(5.3,0)
[xshift=8.6cm] (4.8,0)--(5.3,0);
\draw [densely dotted](5.65,1.8)--(6.05,1.8);
\draw [densely dotted](5.15,3.6)--(5.35,3.6);
\draw[red] (8.6,4)--(19.2,9);
\draw[red] (17.2,4)--(27.8,9);
\draw[red] (25.8,4)--(36.4,9);
\node (top1)[above] at(2,9){\footnotesize$(a,\top)$};
\node (top2)[above] at(10.6,9){\footnotesize$(b,\top)$};
\node (top3)[above] at(19.2,9){\footnotesize$(f_{a,b}(y),\top)$};
\node (top4)[above] at(27.8,9){\footnotesize$(f_{a,b}(y.z),\top)$};
\node (top5)[above] at(36.4,9){\footnotesize$(f_{a,b}(y.z.u),\top)$};
\node (d1) at(4,-1.3){\footnotesize$a$};
\node (d2) at(13,-1.3){\footnotesize$b$};
\node (d3) at(21.5,-1.3){\footnotesize$f_{a,b}(y)$};
\node (d3) at(30,-1.3){\footnotesize$f_{a,b}(y.z)$};
\node (d3) at(38.5,-1.3){\footnotesize$f_{a,b}(y.z.u)$};
\node (l1) [below] at(8.6,4){\tiny$(b,y)$};
\node (l2) [below] at(17.2,4){\tiny$(f_{a,b}(y),z)$};
\node (l3) [below] at(25.6,4){\tiny$(f_{a,b}(y.z),u)$};
\node (l4) [below] at(5.85,0){\tiny$(a,y)$};
\node (l5) [below] at(5.3,1.8){\tiny$(a,y.z)$};
\node (l6) [below] at(5.5,3.6){\tiny$(a,y.z.u)$};
\end{tikzpicture}

  \caption{ Assembling the strict order $<_{3}$ and $<_{4}$}\label{F.6}
\end{figure}

And so it goes on and on. The process can be depicted as Fig.\ref{F.6}.
\begin{figure}[H]
  \centering
  \begin{tikzpicture} [line width=0.5pt,scale=0.37]
\fill[red,opacity=0.2] (-0.2,-0.25) rectangle (0.2,4.25)
[xshift=8.6cm] (-0.2,-0.25) rectangle (0.2,4.25)
[xshift=8.6cm] (-0.2,-0.25) rectangle (0.2,4.25)
[xshift=8.6cm] (-0.2,-0.25) rectangle (0.2,4.25)
[xshift=8.6cm] (-0.2,-0.25) rectangle (0.2,4.25);
\fill[blue,opacity=0.2] (0.4,-0.25) rectangle (8.2,4.25)
[xshift=8.6cm] (0.4,-0.25) rectangle (8.2,4.25)
[xshift=8.6cm] (0.4,-0.25) rectangle (8.2,4.25)
[xshift=8.6cm] (0.4,-0.25) rectangle (8.2,4.25)
[xshift=8.6cm] (0.4,-0.25) rectangle (8.2,4.25);
\fill[black] (0,0) circle (3pt)
[xshift=8.6cm] (0,0) circle (3pt)
[xshift=8.6cm] (0,0) circle (3pt)
[xshift=8.6cm] (0,0) circle (3pt)
[xshift=8.6cm] (0,0) circle (3pt);
\fill[black] (0,1) circle (3pt)
[xshift=8.6cm] (0,1) circle (3pt)
[xshift=8.6cm] (0,1) circle (3pt)
[xshift=8.6cm] (0,1) circle (3pt)
[xshift=8.6cm] (0,1) circle (3pt);
\fill[black] (0,4) circle (3pt)
[xshift=8.6cm] (0,4) circle (3pt)
[xshift=8.6cm] (0,4) circle (3pt)
[xshift=8.6cm] (0,4) circle (3pt)
[xshift=8.6cm] (0,4) circle (3pt);
\draw (0,0)--(0,1)
[xshift=8.6cm] (0,0)--(0,1)
[xshift=8.6cm] (0,0)--(0,1)
[xshift=8.6cm] (0,0)--(0,1)
[xshift=8.6cm] (0,0)--(0,1);
\draw [densely dashed](0,1)--(0,4)
[xshift=8.6cm] (0,1)--(0,4)
[xshift=8.6cm] (0,1)--(0,4)
[xshift=8.6cm] (0,0)--(0,4)
[xshift=8.6cm] (0,0)--(0,4);

\fill[black] (0.6,3.6) circle (3pt)
[xshift=1.1cm] (0.6,3.6) circle (3pt)
[xshift=1.15cm] (0.6,3.6) circle (3pt)
[xshift=1.1cm] (0.6,3.6) circle (3pt)
[xshift=1.15cm] (0.6,3.6) circle (3pt)
[xshift=1.1cm] (0.6,3.6) circle (3pt)
[xshift=3cm] (0.6,3.6) circle (3pt)
[xshift=1.1cm] (0.6,3.6) circle (3pt)
[xshift=1.15cm] (0.6,3.6) circle (3pt)
[xshift=1.1cm] (0.6,3.6) circle (3pt)
[xshift=1.15cm] (0.6,3.6) circle (3pt)
[xshift=1.1cm] (0.6,3.6) circle (3pt)
[xshift=3cm] (0.6,3.6) circle (3pt)
[xshift=1.1cm] (0.6,3.6) circle (3pt)
[xshift=1.15cm] (0.6,3.6) circle (3pt)
[xshift=1.1cm] (0.6,3.6) circle (3pt)
[xshift=1.15cm] (0.6,3.6) circle (3pt)
[xshift=1.1cm] (0.6,3.6) circle (3pt)
[xshift=3cm] (0.6,3.6) circle (3pt)
[xshift=1.1cm] (0.6,3.6) circle (3pt)
[xshift=1.15cm] (0.6,3.6) circle (3pt)
[xshift=1.1cm] (0.6,3.6) circle (3pt)
[xshift=1.15cm] (0.6,3.6) circle (3pt)
[xshift=1.1cm] (0.6,3.6) circle (3pt)
[xshift=3cm] (0.6,3.6) circle (3pt)
[xshift=1.1cm] (0.6,3.6) circle (3pt)
[xshift=1.15cm] (0.6,3.6) circle (3pt)
[xshift=1.1cm] (0.6,3.6) circle (3pt)
[xshift=1.15cm] (0.6,3.6) circle (3pt)
[xshift=1.1cm] (0.6,3.6) circle (3pt);

\fill[black] (0.8,1.8) circle (3pt)
[xshift=1.1cm](0.8,1.8) circle (3pt)
[xshift=1.15cm](0.8,1.8) circle (3pt)
[xshift=1.1cm](0.8,1.8) circle (3pt)
[xshift=1.15cm](0.8,1.8) circle (3pt)
[xshift=1.1cm](0.8,1.8) circle (3pt)
[xshift=3cm] (0.8,1.8) circle (3pt)
[xshift=1.1cm] (0.8,1.8) circle (3pt)
[xshift=1.15cm](0.8,1.8) circle (3pt)
[xshift=1.1cm](0.8,1.8) circle (3pt)
[xshift=1.15cm] (0.8,1.8) circle (3pt)
[xshift=1.1cm] (0.8,1.8) circle (3pt)
[xshift=3cm] (0.8,1.8) circle (3pt)
[xshift=1.1cm] (0.8,1.8) circle (3pt)
[xshift=1.15cm](0.8,1.8) circle (3pt)
[xshift=1.1cm](0.8,1.8) circle (3pt)
[xshift=1.15cm] (0.8,1.8) circle (3pt)
[xshift=1.1cm] (0.8,1.8) circle (3pt)
[xshift=3cm] (0.8,1.8) circle (3pt)
[xshift=1.1cm] (0.8,1.8) circle (3pt)
[xshift=1.15cm](0.8,1.8) circle (3pt)
[xshift=1.1cm](0.8,1.8) circle (3pt)
[xshift=1.15cm] (0.8,1.8) circle (3pt)
[xshift=1.1cm] (0.8,1.8) circle (3pt)
[xshift=3cm] (0.8,1.8) circle (3pt)
[xshift=1.1cm] (0.8,1.8) circle (3pt)
[xshift=1.15cm](0.8,1.8) circle (3pt)
[xshift=1.1cm](0.8,1.8) circle (3pt)
[xshift=1.15cm] (0.8,1.8) circle (3pt)
[xshift=1.1cm] (0.8,1.8) circle (3pt);

\fill[black] (1,3.6) circle (3pt)
[xshift=1.1cm](1,3.6) circle (3pt)
[xshift=1.15cm](1,3.6) circle (3pt)
[xshift=1.1cm](1,3.6) circle (3pt)
[xshift=1.15cm](1,3.6) circle (3pt)
[xshift=1.1cm] (1,3.6) circle (3pt)
[xshift=3cm](1,3.6) circle (3pt)
[xshift=1.1cm](1,3.6) circle (3pt)
[xshift=1.15cm](1,3.6) circle (3pt)
[xshift=1.1cm](1,3.6) circle (3pt)
[xshift=1.15cm] (1,3.6) circle (3pt)
[xshift=1.1cm] (1,3.6) circle (3pt)
[xshift=3cm](1,3.6) circle (3pt)
[xshift=1.1cm](1,3.6) circle (3pt)
[xshift=1.15cm](1,3.6) circle (3pt)
[xshift=1.1cm](1,3.6) circle (3pt)
[xshift=1.15cm] (1,3.6) circle (3pt)
[xshift=1.1cm] (1,3.6) circle (3pt)
[xshift=3cm](1,3.6) circle (3pt)
[xshift=1.1cm](1,3.6) circle (3pt)
[xshift=1.15cm](1,3.6) circle (3pt)
[xshift=1.1cm](1,3.6) circle (3pt)
[xshift=1.15cm] (1,3.6) circle (3pt)
[xshift=1.1cm] (1,3.6) circle (3pt)
[xshift=3cm](1,3.6) circle (3pt)
[xshift=1.1cm](1,3.6) circle (3pt)
[xshift=1.15cm](1,3.6) circle (3pt)
[xshift=1.1cm](1,3.6) circle (3pt)
[xshift=1.15cm] (1,3.6) circle (3pt)
[xshift=1.1cm] (1,3.6) circle (3pt);
\draw (0.6,3.6)--(0.8,1.8)
[xshift=1.1cm](0.6,3.6)--(0.8,1.8)
[xshift=1.15cm](0.6,3.6)--(0.8,1.8)
[xshift=1.1cm](0.6,3.6)--(0.8,1.8)
[xshift=1.15cm](0.6,3.6)--(0.8,1.8)
[xshift=1.1cm](0.6,3.6)--(0.8,1.8)
[xshift=3cm] (0.6,3.6)--(0.8,1.8)
[xshift=1.1cm](0.6,3.6)--(0.8,1.8)
[xshift=1.15cm](0.6,3.6)--(0.8,1.8)
[xshift=1.1cm](0.6,3.6)--(0.8,1.8)
[xshift=1.15cm] (0.6,3.6)--(0.8,1.8)
[xshift=1.1cm] (0.6,3.6)--(0.8,1.8)
[xshift=3cm] (0.6,3.6)--(0.8,1.8)
[xshift=1.1cm](0.6,3.6)--(0.8,1.8)
[xshift=1.15cm](0.6,3.6)--(0.8,1.8)
[xshift=1.1cm](0.6,3.6)--(0.8,1.8)
[xshift=1.15cm] (0.6,3.6)--(0.8,1.8)
[xshift=1.1cm] (0.6,3.6)--(0.8,1.8)
[xshift=3cm] (0.6,3.6)--(0.8,1.8)
[xshift=1.1cm](0.6,3.6)--(0.8,1.8)
[xshift=1.15cm](0.6,3.6)--(0.8,1.8)
[xshift=1.1cm](0.6,3.6)--(0.8,1.8)
[xshift=1.15cm] (0.6,3.6)--(0.8,1.8)
[xshift=1.1cm] (0.6,3.6)--(0.8,1.8)
[xshift=3cm] (0.6,3.6)--(0.8,1.8)
[xshift=1.1cm](0.6,3.6)--(0.8,1.8)
[xshift=1.15cm](0.6,3.6)--(0.8,1.8)
[xshift=1.1cm](0.6,3.6)--(0.8,1.8)
[xshift=1.15cm] (0.6,3.6)--(0.8,1.8)
[xshift=1.1cm] (0.6,3.6)--(0.8,1.8);

\draw (0.8,1.8)--(1,3.6)
[xshift=1.1cm](0.8,1.8)--(1,3.6)
[xshift=1.15cm](0.8,1.8)--(1,3.6)
[xshift=1.1cm](0.8,1.8)--(1,3.6)
[xshift=1.15cm](0.8,1.8)--(1,3.6)
[xshift=1.1cm](0.8,1.8)--(1,3.6)
[xshift=3cm] (0.8,1.8)--(1,3.6)
[xshift=1.1cm](0.8,1.8)--(1,3.6)
[xshift=1.15cm](0.8,1.8)--(1,3.6)
[xshift=1.1cm](0.8,1.8)--(1,3.6)
[xshift=1.15cm] (0.8,1.8)--(1,3.6)
[xshift=1.1cm] (0.8,1.8)--(1,3.6)
[xshift=3cm] (0.8,1.8)--(1,3.6)
[xshift=1.1cm](0.8,1.8)--(1,3.6)
[xshift=1.15cm](0.8,1.8)--(1,3.6)
[xshift=1.1cm](0.8,1.8)--(1,3.6)
[xshift=1.15cm] (0.8,1.8)--(1,3.6)
[xshift=1.1cm] (0.8,1.8)--(1,3.6)
[xshift=3cm] (0.8,1.8)--(1,3.6)
[xshift=1.1cm](0.8,1.8)--(1,3.6)
[xshift=1.15cm](0.8,1.8)--(1,3.6)
[xshift=1.1cm](0.8,1.8)--(1,3.6)
[xshift=1.15cm] (0.8,1.8)--(1,3.6)
[xshift=1.1cm] (0.8,1.8)--(1,3.6)
[xshift=3cm] (0.8,1.8)--(1,3.6)
[xshift=1.1cm](0.8,1.8)--(1,3.6)
[xshift=1.15cm](0.8,1.8)--(1,3.6)
[xshift=1.1cm](0.8,1.8)--(1,3.6)
[xshift=1.15cm] (0.8,1.8)--(1,3.6)
[xshift=1.1cm] (0.8,1.8)--(1,3.6);
\draw [densely dotted](1.15,3.6)--(1.5,3.6)
[xshift=1.1cm] (1.15,3.6)--(1.5,3.6)
[yshift=-1.8cm] (1.15,3.6)--(1.5,3.6)
[shift={(1.15cm, 1.8cm)}](1.15,3.6)--(1.5,3.6)
[xshift=1.1cm] (1.15,3.6)--(1.5,3.6)
[yshift=-1.8cm] (1.15,3.6)--(1.5,3.6)
[shift={(1.15cm,1.8cm)}] (1.15,3.6)--(1.5,3.6)
[xshift=1.1cm](1.15,3.6)--(1.5,3.6)
[yshift=-1.8cm] (1.15,3.6)--(1.5,3.6)
[shift={(3cm,1.8cm)}] (1.15,3.6)--(1.5,3.6)
[xshift=1.1cm] (1.15,3.6)--(1.5,3.6)
[yshift=-1.8cm] (1.15,3.6)--(1.5,3.6)
[shift={(1.15cm, 1.8cm)}](1.15,3.6)--(1.5,3.6)
[xshift=1.1cm] (1.15,3.6)--(1.5,3.6)
[yshift=-1.8cm] (1.15,3.6)--(1.5,3.6)
[shift={(1.15cm,1.8cm)}] (1.15,3.6)--(1.5,3.6)
[xshift=1.1cm](1.15,3.6)--(1.5,3.6)
[yshift=-1.8cm] (1.15,3.6)--(1.5,3.6)
[shift={(3cm,1.8cm)}] (1.15,3.6)--(1.5,3.6)
[xshift=1.1cm] (1.15,3.6)--(1.5,3.6)
[yshift=-1.8cm] (1.15,3.6)--(1.5,3.6)
[shift={(1.15cm, 1.8cm)}](1.15,3.6)--(1.5,3.6)
[xshift=1.1cm] (1.15,3.6)--(1.5,3.6)
[yshift=-1.8cm] (1.15,3.6)--(1.5,3.6)
[shift={(1.15cm,1.8cm)}] (1.15,3.6)--(1.5,3.6)
[xshift=1.1cm](1.15,3.6)--(1.5,3.6)
[yshift=-1.8cm] (1.15,3.6)--(1.5,3.6)
[shift={(3cm,1.8cm)}] (1.15,3.6)--(1.5,3.6)
[xshift=1.1cm] (1.15,3.6)--(1.5,3.6)
[yshift=-1.8cm] (1.15,3.6)--(1.5,3.6)
[shift={(1.15cm, 1.8cm)}](1.15,3.6)--(1.5,3.6)
[xshift=1.1cm] (1.15,3.6)--(1.5,3.6)
[yshift=-1.8cm] (1.15,3.6)--(1.5,3.6)
[shift={(1.15cm,1.8cm)}] (1.15,3.6)--(1.5,3.6)
[xshift=1.1cm](1.15,3.6)--(1.5,3.6)
[yshift=-1.8cm] (1.15,3.6)--(1.5,3.6)
[shift={(3cm,1.8cm)}] (1.15,3.6)--(1.5,3.6)
[xshift=1.1cm] (1.15,3.6)--(1.5,3.6)
[yshift=-1.8cm] (1.15,3.6)--(1.5,3.6)
[shift={(1.15cm, 1.8cm)}](1.15,3.6)--(1.5,3.6)
[xshift=1.1cm] (1.15,3.6)--(1.5,3.6)
[yshift=-1.8cm] (1.15,3.6)--(1.5,3.6)
[shift={(1.15cm,1.8cm)}] (1.15,3.6)--(1.5,3.6)
[xshift=1.1cm](1.15,3.6)--(1.5,3.6)
[yshift=-1.8cm] (1.15,3.6)--(1.5,3.6);
\fill[black] (1.35,0) circle (3pt)
[xshift=2.25cm](1.35,0) circle (3pt)
[xshift=2.25cm] (1.35,0) circle (3pt)
[xshift=4.1cm](1.35,0) circle (3pt)
[xshift=2.25cm] (1.35,0) circle (3pt)
[xshift=2.25cm](1.35,0) circle (3pt)
[xshift=4.1cm] (1.35,0) circle (3pt)
[xshift=2.25cm] (1.35,0) circle (3pt)
[xshift=2.25cm](1.35,0) circle (3pt)
[xshift=4.1cm](1.35,0) circle (3pt)
[xshift=2.25cm] (1.35,0) circle (3pt)
[xshift=2.25cm](1.35,0) circle (3pt)
[xshift=4.1cm](1.35,0) circle (3pt)
[xshift=2.25cm] (1.35,0) circle (3pt)
[xshift=2.25cm](1.35,0) circle (3pt);

\draw (1.35,0)--(1.9,1.8)
[xshift=2.25cm] (1.35,0)--(1.9,1.8)
[xshift=2.25cm] (1.35,0)--(1.9,1.8)
[xshift=4.1cm] (1.35,0)--(1.9,1.8)
[xshift=2.25cm] (1.35,0)--(1.9,1.8)
[xshift=2.25cm] (1.35,0)--(1.9,1.8)
[xshift=4.1cm] (1.35,0)--(1.9,1.8)
[xshift=2.25cm] (1.35,0)--(1.9,1.8)
[xshift=2.25cm] (1.35,0)--(1.9,1.8)
[xshift=4.1cm] (1.35,0)--(1.9,1.8)
[xshift=2.25cm] (1.35,0)--(1.9,1.8)
[xshift=2.25cm] (1.35,0)--(1.9,1.8)
[xshift=4.1cm] (1.35,0)--(1.9,1.8)
[xshift=2.25cm] (1.35,0)--(1.9,1.8)
[xshift=2.25cm] (1.35,0)--(1.9,1.8);

\draw (0.8,1.8)--(1.35,0)
[xshift=2.25cm]  (0.8,1.8)--(1.35,0)
[xshift=2.25cm]  (0.8,1.8)--(1.35,0)
[xshift=4.1cm]  (0.8,1.8)--(1.35,0)
[xshift=2.25cm]  (0.8,1.8)--(1.35,0)
[xshift=2.25cm]  (0.8,1.8)--(1.35,0)
[xshift=4.1cm]  (0.8,1.8)--(1.35,0)
[xshift=2.25cm]  (0.8,1.8)--(1.35,0)
[xshift=2.25cm]  (0.8,1.8)--(1.35,0)
[xshift=4.1cm]  (0.8,1.8)--(1.35,0)
[xshift=2.25cm]  (0.8,1.8)--(1.35,0)
[xshift=2.25cm]  (0.8,1.8)--(1.35,0)
[xshift=4.1cm]  (0.8,1.8)--(1.35,0)
[xshift=2.25cm]  (0.8,1.8)--(1.35,0)
[xshift=2.25cm]  (0.8,1.8)--(1.35,0);
\draw [densely dotted](7.3,3.6)--(8.1,3.6)
[yshift=-1.8cm] (7.3,3.6)--(8.1,3.6)
[yshift=-1.8cm] (7.3,3.6)--(8.1,3.6)
[xshift=8.6cm] (7.3,3.6)--(8.1,3.6)
[yshift=1.8cm] (7.3,3.6)--(8.1,3.6)
[yshift=1.8cm] (7.3,3.6)--(8.1,3.6)
[xshift=8.6cm] (7.3,3.6)--(8.1,3.6)
[yshift=-1.8cm] (7.3,3.6)--(8.1,3.6)
[yshift=-1.8cm] (7.3,3.6)--(8.1,3.6)
[xshift=8.6cm] (7.3,3.6)--(8.1,3.6)
[yshift=1.8cm] (7.3,3.6)--(8.1,3.6)
[yshift=1.8cm] (7.3,3.6)--(8.1,3.6)
[xshift=8.6cm] (7.3,3.6)--(8.1,3.6)
[yshift=-1.8cm] (7.3,3.6)--(8.1,3.6)
[yshift=-1.8cm] (7.3,3.6)--(8.1,3.6);
%
\fill[black] (2,9) circle (3pt)
[xshift=8.6cm] (2,9) circle (3pt)
[xshift=8.6cm] (2,9) circle (3pt)
[xshift=8.6cm] (2,9) circle (3pt)
[xshift=8.6cm] (2,9) circle (3pt);
\draw[densely dashed] (2,9)--(1.35,4)
[xshift=8.6cm] (2,9)--(1.35,4)
[xshift=8.6cm] (2,9)--(1.35,4)
[xshift=8.6cm] ((2,9)--(1.35,4)
[xshift=8.6cm] (2,9)--(1.35,4);
\draw[densely dashed] (0,4)--(2,9)
[xshift=8.6cm] (0,4)--(2,9)
[xshift=8.6cm] (0,4)--(2,9)
[xshift=8.6cm] (0,4)--(2,9)
[xshift=8.6cm] (0,4)--(2,9);
\draw[densely dashed] (2,9)--(6,4)
[xshift=8.6cm] (2,9)--(6,4)
[xshift=8.6cm] (2,9)--(6,4)
[xshift=8.6cm] (2,9)--(6,4)
[xshift=8.6cm] (2,9)--(6,4);

\draw[densely dashed] (3.6,4)--(2,9)
[xshift=8.6cm] (3.6,4)--(2,9)
[xshift=8.6cm] (3.6,4)--(2,9)
[xshift=8.6cm] (3.6,4)--(2,9)
[xshift=8.6cm] (3.6,4)--(2,9);
\draw [densely dotted](4.8,0)--(5.3,0)
[xshift=8.6cm] (4.8,0)--(5.3,0);
\draw [densely dotted](5.65,1.8)--(6.05,1.8);
\draw [densely dotted](5.15,3.6)--(5.35,3.6);
\draw[red] (8.6,4)--(19.2,9);
\draw[red] (17.2,4)--(27.8,9);
\draw[red] (25.8,4)--(36.4,9);
\draw[blue] (5.85,0)--(19.2,9);
\draw[blue] (5.3,1.8)--(27.8,9);
\draw[blue] (5.5,3.6)--(36.4,9);
\node (top1)[above] at(2,9){\footnotesize$(a,\top)$};
\node (top2)[above] at(10.6,9){\footnotesize$(b,\top)$};
\node (top3)[above] at(19.2,9){\footnotesize$(f_{a,b}(y),\top)$};
\node (top4)[above] at(27.8,9){\footnotesize$(f_{a,b}(y.z),\top)$};
\node (top5)[above] at(36.4,9){\footnotesize$(f_{a,b}(y.z.u),\top)$};
\node (d1) at(4,-1.3){\footnotesize$a$};
\node (d2) at(13,-1.3){\footnotesize$b$};
\node (d3) at(21.5,-1.3){\footnotesize$f_{a,b}(y)$};
\node (d3) at(30,-1.3){\footnotesize$f_{a,b}(y.z)$};
\node (d3) at(38.5,-1.3){\footnotesize$f_{a,b}(y.z.u)$};
\node (l1) [below] at(8.6,4){\tiny$(b,y)$};
\node (l2) [below] at(17.2,4){\tiny$(f_{a,b}(y),z)$};
\node (l3) [below] at(25.6,4){\tiny$(f_{a,b}(y.z),u)$};
\node (l4) [below] at(5.85,0){\tiny$(a,y)$};
\node (l5) [below] at(5.3,1.8){\tiny$(a,y.z)$};
\node (l6) [below] at(5.5,3.6){\tiny$(a,y.z.u)$};
\end{tikzpicture}

  \caption{ Assembling the strict order $<_{2}$, $<_{3}$ and $<_{4}$}\label{F.7}
\end{figure}
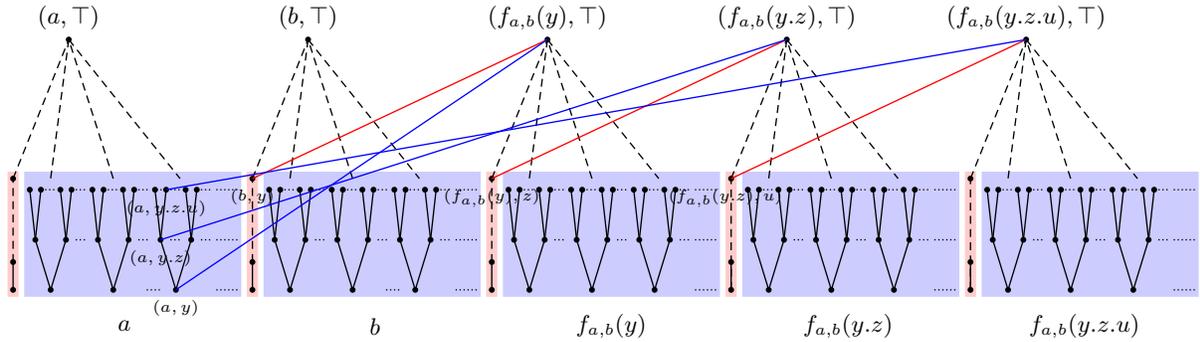

 In Fig.\ref{F.7}, the red lines are the same as Fig\ref{F.6} and the blue lines add the cases of $<_{2}$: $(a,y)<_{2}(f_{a,b}(y),\top)$, $(a,y.z)<_{2}(f_{a,b}(y.z),\top)$ and $(a,y.z.u)<_{2}(f_{a,b}(y.z.u),\top)$.
\begin{lemma}\label{b}
  Let $P$ be equipped with the order $\leq$. Then $P$ is an irreducible subset of $\Sigma P$.
  \begin{proof}
    From the definition of irreducibility, it suffices to prove that $U\cap V\neq \emptyset$ for any non-empty Scott open sets $U,V$ of $P$. Now we choose $(n_{0},x)\in U$, $(m_{0},y)\in V$. If $n_{0}=m_{0}$, then $(n_{0},\top)\in U\cap V$ through the Scott-openness of $U,V$. Otherwise, $n_{0}\neq m_{0}$. Without loss of generality, we can assume that $n_{0}<m_{0}$.

     Using again the fact that $V$ is Scott open, it is straight forward to show that there exists $a_{1}\in \mathbb{N}$ such that $(m_{0}, a_{1})\in V$. From the definition of $<_{3}$, we can see that $(m_{0},a_{1})<_{3}(f_{n_{0},m_{0}}(a_{1}),\top)$. Whence, $(f_{n_{0},m_{0}}(a_{1}),\top)\in V$. The Scott openness of $V$ implies that $(f_{n_{0},m_{0}}(a_{1}),a_{2})\in V$ for some $a_{2}\in \mathbb{N}$. Due to the definition of $<_{4}$, we conclude that $(f_{n_{0},m_{0}}(a_{1}),a_{2})<_{4}(f_{n_{0},m_{0}}(a_{1}.a_{2}),\top)$. It follows that $(f_{n_{0},m_{0}}(a_{1}.a_{2}),\top)\in V$.

    By induction on $\mathbb{N}$, for any $n\in \mathbb{N}$, there exists $(f_{n_{0},m_{0}}(a_{1}.a_{2}.\cdots .a_{n}),\top)\in V$. It is worth noting that $\sup_{n\in \mathbb{N}}(n_{0},a_{1}.\cdots.a_{n})=(n_{0},\top)\in U$. This indicates that there exists $k\in \mathbb{N}$ such that $(n_{0},a_{1}.\cdots.a_{k})\in U$ by the Scott openness of $U$. Through the definition of $<_{2}$, we can deduce that $(f_{n_{0},m_{0}}(a_{1}.\cdots.a_{k}),\top)\in U$. This means that $(f_{n_{0},m_{0}}(a_{1}.\cdots.a_{k}),\top)\in U\cap V$.
  \end{proof}
\end{lemma}
\begin{figure}[H]
  \centering
  \includegraphics[height=5cm, width=16cm]{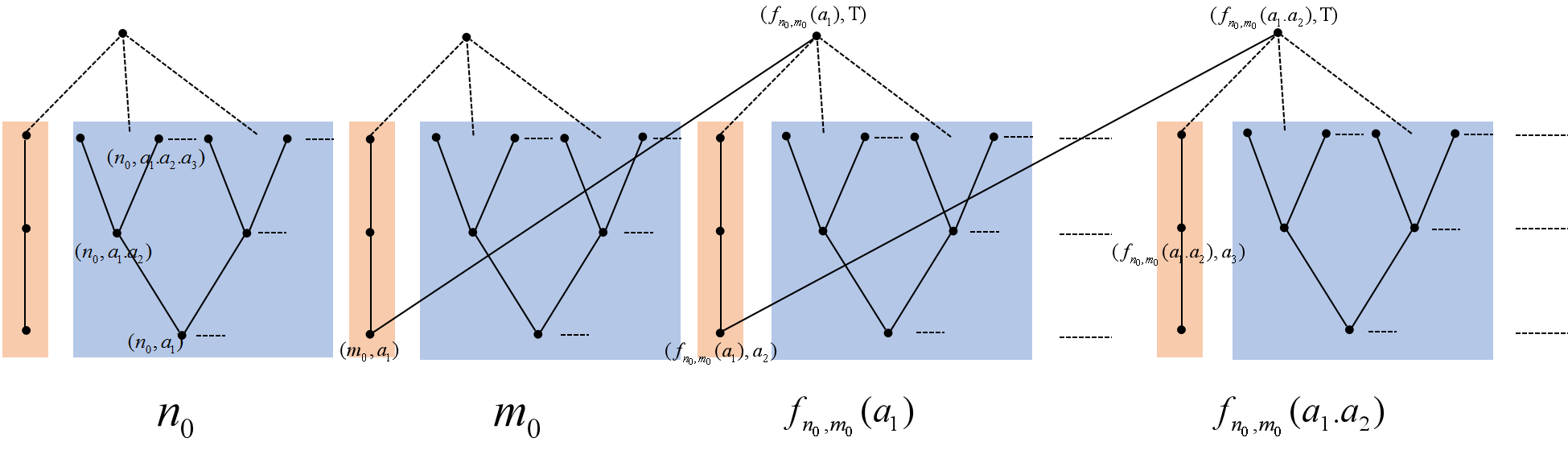}\\
  \caption{ The proof of Lemma \ref{b} ($n<m$)}\label{F.8}
\end{figure}
\begin{lemma}\label{a}
  Let $M=\{\bigcap_{x\in E}\da x\mid E\subseteq P\}$. Then $(M,\subseteq)$ is a bounded complete dcpo.
  \begin{proof}
    Obviously, it remains to testify that $M$ is a dcpo. Let $B=\{(n,m)\in \mathbb{N}\times \mathbb{N}\mid n<m\}$. In order to determine what the intersections of two principal ideals of $P$ are, we first classify the principal ideals $\da x$ of $P$. Let us classify them into five types.

    Type I is the case $\{(n_{0},s)\mid s\leq s_{0}\}$ for some $n_{0}\in \mathbb{N},s_{0}\in \mathbb{N}^{<\mathbb{N}}$.

    Type II is the case $\{(m_{0},n)\mid n\leq n_{0}\}$ for some $m_{0},n_{0}\in \mathbb{N}$.

    Type III is the case $L_{n_{0}}$ for some $n_{0}\in \mathbb{N}\backslash \bigcup_{(n,m)\in B}i(n,m)$.

    Type IV is the case $\da(f_{m_{0},n_{0}}(s_{0}),\top)=L_{f_{m_{0},n_{0}}(s_{0})}\cup \{(m_{0},s_{0})\}\cup\{(n_{0},n)\mid n\leq s_{0}\}$ for some $(m_{0},n_{0})\in B, s_{0}\in \mathbb{N}^{<\mathbb{N}}$ with $ {\mid}s_{0}{\mid}=1$.

    Type V is the case $\da(f_{m_{0},n_{0}}(s_{0}),\top)=L_{f_{m_{0},n_{0}}(s_{0})}\cup \{(m_{0},s)\mid s\leq s_{0}\}\cup\{(f_{m_{0},n_{0}}(s_{0}^{*}),n)\mid n\leq n_{0}^{*}\}$ for some $(m_{0},n_{0})\in B, s_{0}=s_{0}^{*}.n_{0}^{*}\in \mathbb{N}^{<\mathbb{N}}$ with $s_{0}^{*}\in \mathbb{N}^{<\mathbb{N}},n_{0}^{*}\in \mathbb{N}$.

    All those cases are depicted as below, as blue regions.
\begin{figure}[H]
  \centering
  \includegraphics[height=5cm]{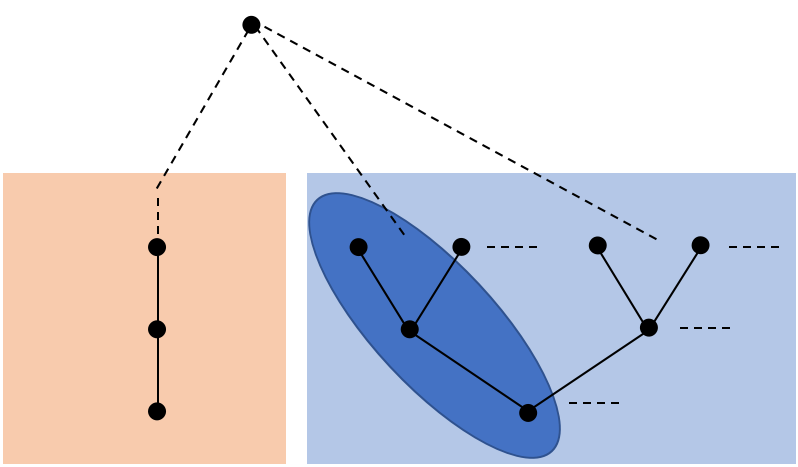}\\
  \caption{ The Type I ideals}
\end{figure}
\begin{figure}[H]
  \centering
  \includegraphics[height=5cm]{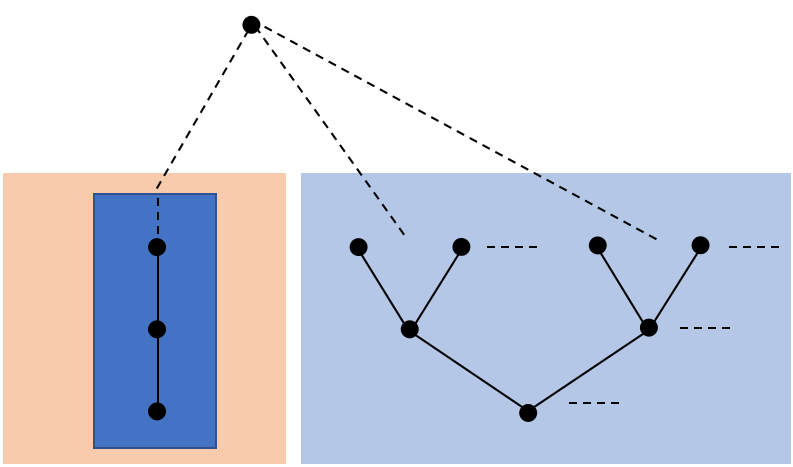}\\
  \caption{ The Type II ideals}
\end{figure}
\begin{figure}[H]
  \centering
  \includegraphics[height=5cm]{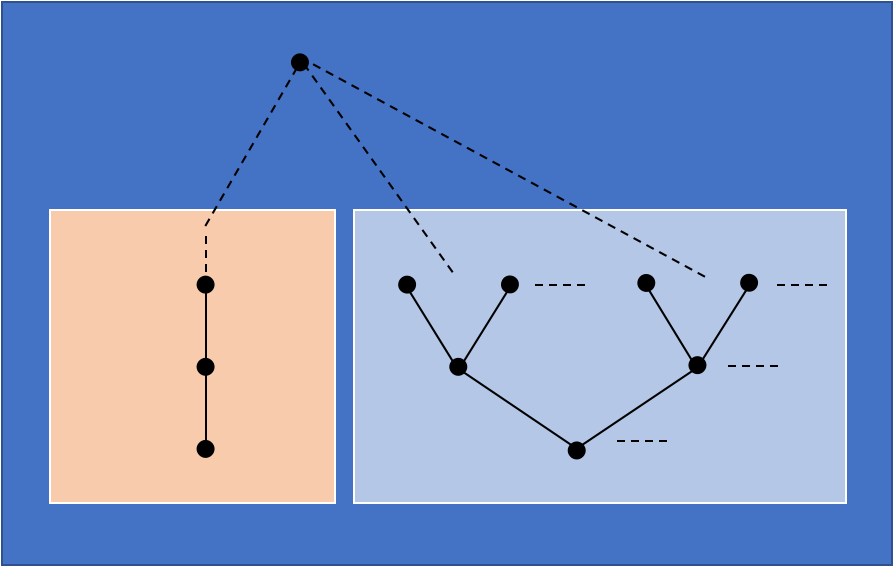}\\
  \caption{ The Type III ideals}
\end{figure}
\begin{figure}[H]
  \centering
  \includegraphics[height=5cm, width=16cm]{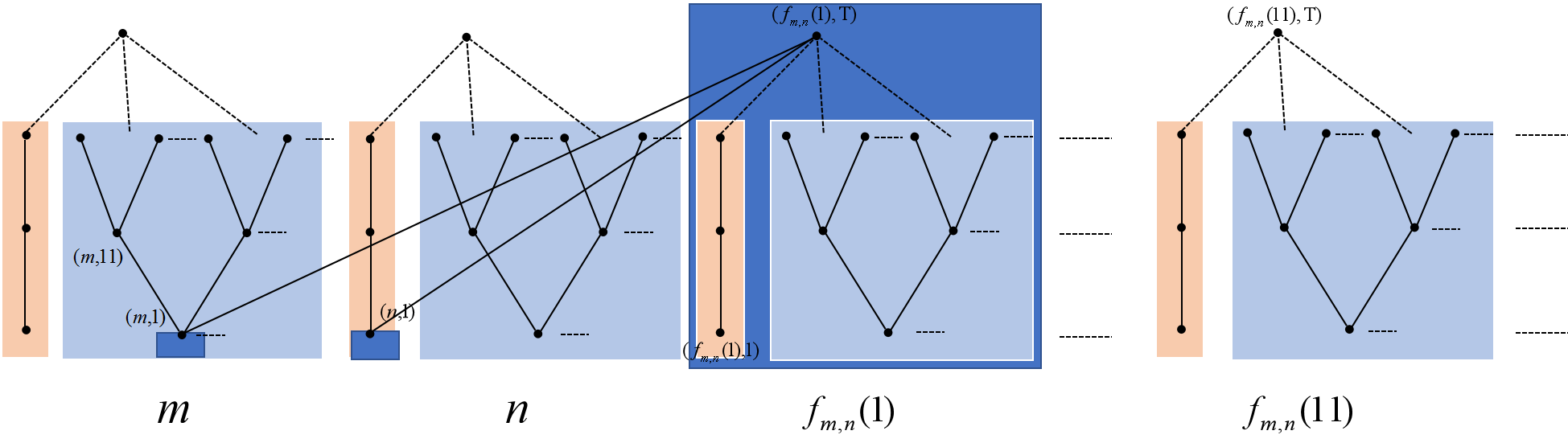}\\
  \caption{ The Type IV ideals}
\end{figure}
\begin{figure}[H]
  \centering
  \includegraphics[height=5cm, width=16cm]{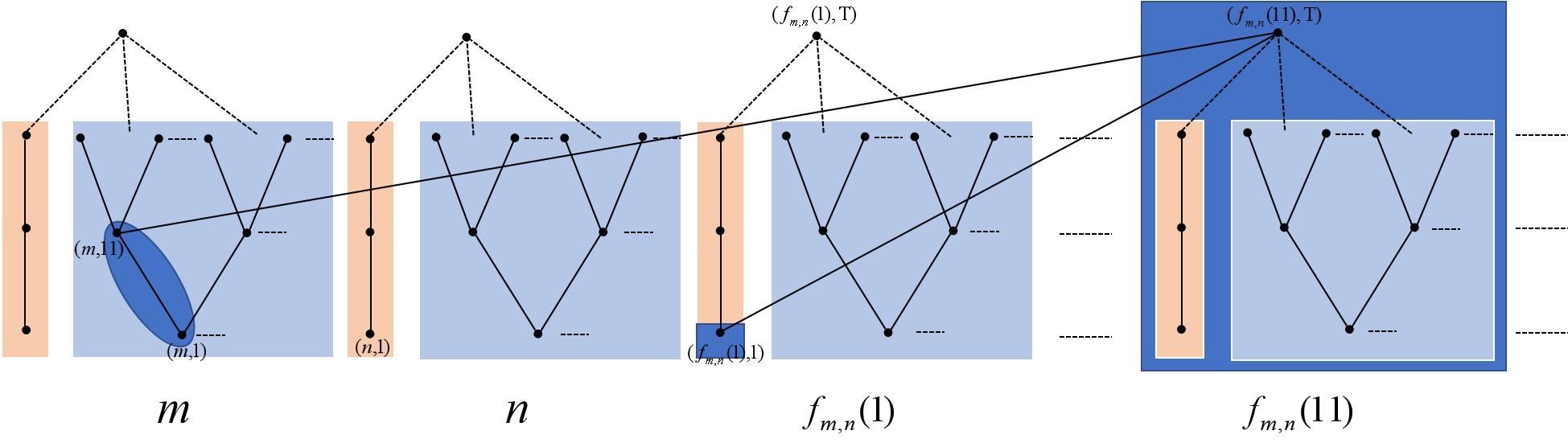}\\
  \caption{ The Type V ideals}
\end{figure}
    I have summarized all the subsets of $P$ that we can obtain by intersecting two principal ideals in the following table.
      \begin{table}[H]
  \centering

  \label{T.smallcontext}
  \vspace{0.2cm}
  \centering
  \scalebox{1}{
  \begin{tabular}{lccccc}
  \hline
           &  Type I & Type II  & Type III  &  Type IV  &  Type V\\
  \hline
  Type I & I/$\emptyset$      & $\emptyset$      & I/$\emptyset$      & I/$\emptyset$ & I/$\emptyset$    \\
  \hline
  Type II &      & II/$\emptyset$      & II/$\emptyset$     & II/$\emptyset$   & II/$\emptyset$  \\
  \hline
  Type III &      &      & III/$\emptyset$      & I/II/$\emptyset$  & I/II/$\emptyset$   \\
  \hline
  Type IV &       &       &       & I/II/IV/I$\cup$$\mathrm{II}^{1}$/$\emptyset$  & I/II/I$\cup$$\mathrm{II}^{1}$/I$\cup$$\mathrm{II}^{2}$/$\emptyset$\\
  \hline
  Type V &       &       &       &  & I/II/V/I$\cup$$\mathrm{II}^{2}$/$\emptyset$ \\
  \hline
  \end{tabular}}
\end{table}

 In above table, Type I$\cup \mathrm{II}^{1} =\{(m_{0},s_{0})\}\cup\{(n_{0},n)\mid n\leq k_{0}\}$ for some $(m_{0},n_{0})\in B, s_{0}\in \mathbb{N}^{<\mathbb{N}},k_{0}\in \mathbb{N}$ with $ {\mid}s_{0}{\mid}=1$, $k_{0}\leq s_{0}$.

Type I$\cup \mathrm{II}^{2} =\{(m_{0},s)\mid s\leq s_{0}\}\cup\{(f_{m_{0},n_{0}}(s_{0}),n)\mid n\leq k_{0}\}$ for some $(m_{0},n_{0})\in B, s_{0}\in \mathbb{N}^{<\mathbb{N}},k_{0}\in \mathbb{N}$.

The two cases are depicted as below, as blue regions.
\begin{figure}[H]
  \centering
  \includegraphics[height=5cm, width=16cm]{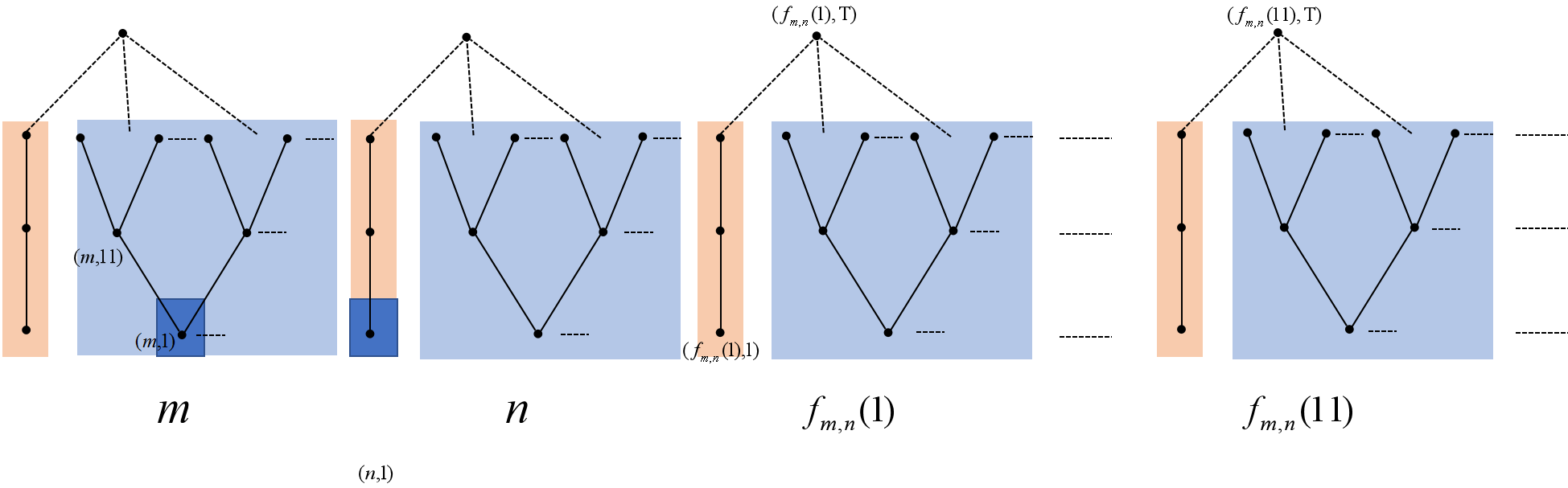}\\
  \caption{ The Type $\mathrm{I}\cup \mathrm{II}^{1}$ ideals}
\end{figure}
\begin{figure}[H]
  \centering
  \includegraphics[height=5cm, width=16cm]{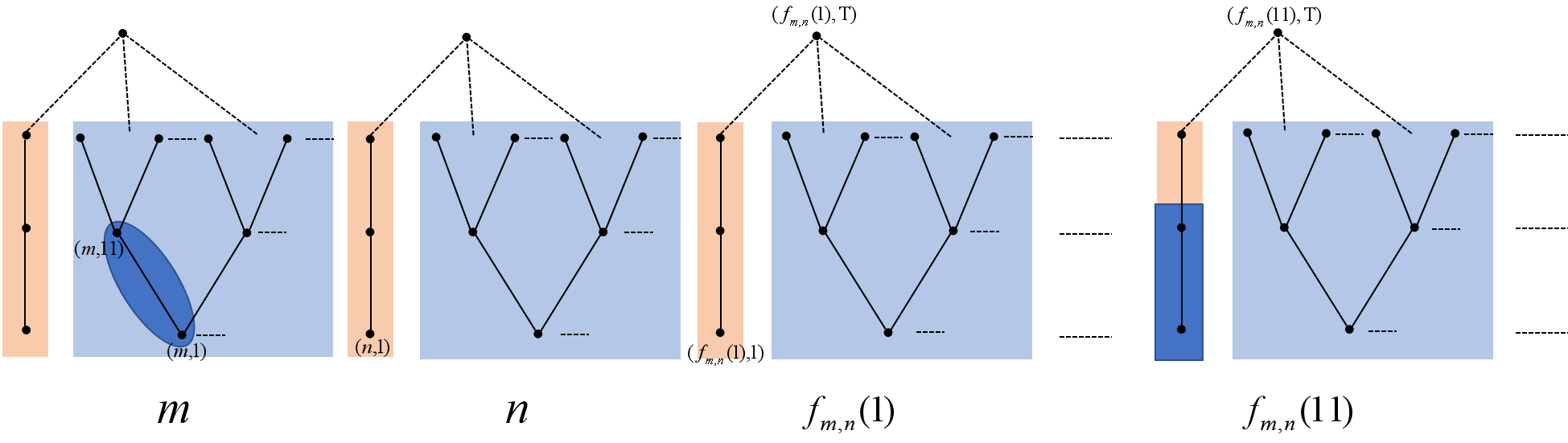}\\
  \caption{ The Type $\mathrm{I}\cup \mathrm{II}^{2}$ ideals}
\end{figure}

    The interesting cases are what happens when you intersect two Type IV ideals, two Type V ideals or Type IV ideal and Type V ideal.

    The corresponding cell of two Type IV ideals says ¡°I/II/IV/I$\cup$$ \mathrm{II}^{1}$/$\emptyset$¡±, and that means that the intersection can be a Type I ideal, a Type II ideal, a type IV ideal, the Type $\mathrm{I}\cup \mathrm{II}^{1}$ ideal or the empty set.

    Suppose $I_{1}=\da(f_{m_{1},n_{1}}(s_{1}),\top)=L_{f_{m_{1},n_{1}}(s_{1})}\cup \{(m_{1},s_{1})\}\cup\{(n_{1},n)\mid n\leq s_{1}\}$ for some $(m_{1},n_{1})\in B, s_{1}\in \mathbb{N}^{<\mathbb{N}}$ with $ {\mid}s_{1}{\mid}=1$,

    $I_{2}=\da(f_{m_{2},n_{2}}(s_{2}),\top)=L_{f_{m_{2},n_{2}}(s_{2})}\cup \{(m_{2},s_{2})\}\cup\{(n_{2},n)\mid n\leq s_{2}\}$ for some $(m_{2},n_{2})\in B, s_{2}\in \mathbb{N}^{<\mathbb{N}}$ with $ {\mid}s_{2}{\mid}=1$. Then $I_{1},I_{2}$ are two Type IV ideals.

    We now distinguish the following cases for $f_{m_{2},n_{2}}(s_{2})$:

    Case 1, $f_{m_{2},n_{2}}(s_{2})<m_{1}$: Then $I_{1}\cap I_{2}=\emptyset$.

    Case 2, $f_{m_{2},n_{2}}(s_{2})=m_{1}$: Then $I_{1}\cap I_{2}=\{(m_{1},s_{1})\}$, which is a Type I ideal.

    Case 3, $m_{1}<f_{m_{2},n_{2}}(s_{2})<n_{1}$: Then $I_{1}\cap I_{2}=\emptyset$ in case $m_{2}\neq m_{1}$. In case $m_{2}=m_{1}$, if $s_{1}\neq s_{2}$, then $I_{1}\cap I_{2}=\emptyset$. Otherwise, $I_{1}\cap I_{2}=\{(m_{1},s_{1})\}$, which is a Type I ideal.

    Case 4, $f_{m_{2},n_{2}}(s_{2})=n_{1}$: In case $m_{2}\neq m_{1}$, we conclude that $I_{1}\cap I_{2}=\{(n_{1},n)\mid n\leq s_{1}\}$, which is a Type II ideal. In case $m_{2}=m_{1}$, if $s_{1}\neq  s_{2}$, then we have the same result as in case $m_{2}\neq m_{1}$. Otherwise, $s_{1}=  s_{2}$. Then $I_{1}\cap I_{2}=\{(m_{1},s_{1})\}\cup \{(n_{1},s)\mid s\leq s_{1}\}$, which is a Type I$\cup$$\mathrm{II}^{1}$ ideal.

    Case 5, $n_{1}<f_{m_{2},n_{2}}(s_{2})<f_{m_{1},n_{1}}(s_{1})$: Then the case $m_{1}\neq m_{2},n_{1}\neq n_{2}$ implies that $I_{1}\cap I_{2}=\emptyset$.

    In case $m_{1}=m_{2},n_{1}\neq n_{2}$, suppose $s_{1}\neq s_{2}$. Then $I_{1}\cap I_{2}=\emptyset$. Otherwise, $s_{1}=s_{2}$. This implies that $ I_{1}\cap I_{2}=\{(m_{1},s_{1})\}$, which is a Type I ideal.

     In case, $m_{1}=m_{2},n_{1}= n_{2}$. As a result, $s_{1}< s_{2}$ follows immediately due to the fact that $f_{m_{1},n_{1}}$ is an monotone injection, which contradicts the assumption that ${\mid}s_{1}{\mid}={\mid}s_{2}{\mid}=1$.

     The case $m_{1}\neq m_{2}$, $n_{1}=n_{2}$ implies that $I_{1}\cap I_{2}=\{(n_{1},n)\mid n\leq \min\{s_{1},s_{2}\}\}$, which is a Type II ideal.

    Case 6, $f_{m_{2},n_{2}}(s_{2})=f_{m_{1},n_{1}}(s_{1})$: Then $(m_{1},n_{1})=(m_{2},n_{2}),s_{1}=s_{2}$ by the property of $i$ and $f_{m_{1},n_{1}}$. This reveals that $I_{1}\cap I_{2}=I_{1}$, which is a type IV ideal.

    Note that the remain case for $f_{m_{2},n_{2}}(s_{2})$ is symmetric with the above cases. This covers all possible cases and we have confirmed that the intersection of two type IV can be a Type I ideal, a Type II ideal, a type IV ideal, the type $\mathrm{I}\cup \mathrm{II}^{1}$ ideal or the empty set.

  The corresponding cell of Type IV ideal and Type V ideal says ¡°I/II//I$\cup$$ \mathrm{II}^{1}$/I$\cup$$\mathrm{II}^{2}$/$\emptyset$¡±, and that means that the intersection can be a Type I ideal, a Type II ideal, the Type $\mathrm{I}\cup \mathrm{II}^{1}$ ideal, the Type I$\cup$$\mathrm{II}^{2}$ or the empty set.

   Suppose $I_{1}=\da(f_{m_{1},n_{1}}(s_{1}),\top)=L_{f_{m_{1},n_{1}}(s_{1})}\cup \{(m_{1},s_{1})\}\cup\{(n_{1},n)\mid n\leq s_{1}\}$ for some $(m_{1},n_{1})\in B, s_{1}\in \mathbb{N}^{<\mathbb{N}}$ with $ {\mid}s_{1}{\mid}=1$,

    $I_{2}=\da(f_{m_{2},n_{2}}(s_{2}),\top)=L_{f_{m_{2},n_{2}}(s_{2})}\cup \{(m_{2},s)\mid s\leq s_{2}\}\cup\{(f_{m_{2},n_{2}}(s_{2}^{*}),n)\mid n\leq n_{2}^{*}\}$ for some $(m_{2},n_{2})\in B, s_{2}=s_{2}^{*}.n_{2}^{*}\in \mathbb{N}^{<\mathbb{N}}$ with $s_{2}^{*}\in \mathbb{N}^{<\mathbb{N}},n_{2}^{*}\in \mathbb{N}$. Then $I_{1}$ is a Type IV ideal, $I_{2}$ a Type V ideal.

    We now distinguish the following cases for $f_{m_{2},n_{2}}(s_{2})$:

    Case 1, $f_{m_{2},n_{2}}(s_{2})<m_{1}$: Then $I_{1}\cap I_{2}=\emptyset$.

    Case 2, $f_{m_{2},n_{2}}(s_{2})=m_{1}$: Then $I_{1}\cap I_{2}=\{(m_{1},s_{1})\}$, which is a Type I ideal.

    Case 3, $m_{1}<f_{m_{2},n_{2}}(s_{2})<n_{1}$: Then $I_{1}\cap I_{2}=\emptyset$ in case $m_{2}\neq m_{1}$. In case $m_{2}=m_{1}$, if $s_{1}\nleq s_{2}$, then $I_{1}\cap I_{2}=\emptyset$. Otherwise, $I_{1}\cap I_{2}=\{(m_{1},s_{1})\}$, which is a Type I ideal.

    Case 4, $f_{m_{2},n_{2}}(s_{2})=n_{1}$: In case $m_{2}\neq m_{1}$, we conclude that $I_{1}\cap I_{2}=\{(n_{1},n)\mid n\leq s_{1}\}$, which is a Type II ideal. In case $m_{2}=m_{1}$, if $s_{1}\nleq  s_{2}$, then we have the same result as in case $m_{2}\neq m_{1}$. Otherwise, $s_{1}\leq s_{2}$. Then $I_{1}\cap I_{2}=\{(m_{1},s_{1})\}\cup \{(n_{1},s)\mid s\leq s_{1}\}$, which is a Type I$\cup$$\mathrm{II}^{1}$ ideal.

    Case 5, $n_{1}<f_{m_{2},n_{2}}(s_{2})<f_{m_{1},n_{1}}(s_{1})$: Then the case $m_{1}\neq m_{2},n_{1}\neq f_{m_{2},n_{2}}(s_{2}^{*})$ implies that $I_{1}\cap I_{2}=\emptyset$.

    In case $m_{1}=m_{2},n_{1}\neq f_{m_{2},n_{2}}(s_{2}^{*})$, suppose $s_{1}\nleq s_{2}$. Then $I_{1}\cap I_{2}=\emptyset$. Otherwise, $s_{1}\leq s_{2}$. This implies that $ I_{1}\cap I_{2}=\{(m_{1},s_{1})\}$, which is a Type I ideal.

     In case, $m_{1}=m_{2},n_{1}= f_{m_{2},n_{2}}(s_{2}^{*})$. If $s_{1}\nleq s_{2}$, then $I_{1}\cap I_{2}=\{(n_{1},n)\mid n\leq \min\{s_{1},n_{2}^{*}\}\}$, which is a Type II ideal. Otherwise $s_{1}\leq s_{2}$, then $I_{1}\cap I_{2}=\{(m_{1},s_{1})\}\cup \{(n_{1},n)\mid n\leq \min\{s_{1},n_{2}^{*}\}\}$, which is a Type $\mathrm{I}\cup \mathrm{II}^{1}$ ideal.

     In case $m_{1}\neq m_{2}$, $n_{1}=f_{m_{2},n_{2}}(s_{2}^{*})$, then $I_{1}\cap I_{2}=\{(n_{1},n)\mid n\leq \min\{s_{1},n_{2}^{*}\}\}$, which is a Type II ideal.

    Case 6, $f_{m_{2},n_{2}}(s_{2})=f_{m_{1},n_{1}}(s_{1})$: Then $(m_{1},n_{1})=(m_{2},n_{2}),s_{1}=s_{2}$ by the property of $i$ and $f_{m_{1},n_{1}}$, which contradicts the assumption that $s_{1}\neq s_{2}$.

    Case 7, $f_{m_{2},n_{2}}(s_{2})>f_{m_{1},n_{1}}(s_{1})$: Then we need to distinguish $f_{m_{1},n_{1}}(s_{1})$ in this case.

    Case 7.1, $f_{m_{1},n_{1}}(s_{1})<m_{2}$: Then $I_{1}\cap I_{2}=\emptyset$.

    Case 7.2, $f_{m_{1},n_{1}}(s_{1})=m_{2}$: Then $I_{1}\cap I_{2}=\{(m_{2},s)\mid s\leq s_{2}\}$, which is a Type I ideal.

    Case 7.3, $m_{2}<f_{m_{1},n_{1}}(s_{1})<f_{m_{2},n_{2}}(s_{2}^{*})$: Then $I_{1}\cap I_{2}=\emptyset$ in case $m_{2}\neq m_{1}$. In case $m_{2}=m_{1}$, if $s_{1}\nleq s_{2}$, then $I_{1}\cap I_{2}=\emptyset$. Otherwise, $I_{1}\cap I_{2}=\{(m_{1},s_{1})\}$, which is a Type I ideal.

    Case 7.4, $f_{m_{1},n_{1}}(s_{1})=f_{m_{2},n_{2}}(s_{2}^{*})$: Then $(m_{1},n_{1})=(m_{2},n_{2}), s_{2}^{*}=s_{1}$ by the property of $i$ and $f_{m_{1},n_{1}}$. This implies that $I_{1}\cap I_{2}=\{(m_{1},s_{1})\}\cup\{(f_{m_{1},n_{1}}(s_{1}),n)\mid n\leq n_{2}^{*}\}$ which is a Type $\mathrm{I}\cup \mathrm{II}^{2}$ ideal.

    Case 7.5, $f_{m_{2},n_{2}}(s_{2}^{*})<f_{m_{1},n_{1}}(s_{1})<f_{m_{2},n_{2}}(s_{2})$: Then the case $m_{1}\neq m_{2},n_{1}\neq f_{m_{2},n_{2}}(s_{2}^{*})$ implies that $I_{1}\cap I_{2}=\emptyset$.

    In case $m_{1}=m_{2},n_{1}\neq f_{m_{2},n_{2}}(s_{2}^{*})$, suppose $s_{1}\nleq s_{2}$. Then $I_{1}\cap I_{2}=\emptyset$. Otherwise, $s_{1}\leq s_{2}$. This implies that $ I_{1}\cap I_{2}=\{(m_{1},s_{1})\}$, which is a Type I ideal.

     In case, $m_{1}=m_{2},n_{1}= f_{m_{2},n_{2}}(s_{2}^{*})$. If $s_{1}\nleq s_{2}$, then $I_{1}\cap I_{2}=\{(n_{1},n)\mid n\leq \min\{s_{1},n_{2}^{*}\}\}$, which is a Type II ideal. Otherwise $s_{1}\leq s_{2}$, then $I_{1}\cap I_{2}=\{(m_{1},s_{1})\}\cup \{(n_{1},n)\mid n\leq \min\{s_{1},n_{2}^{*}\}\}$, which is a Type $\mathrm{I}\cup \mathrm{II}^{1}$ ideal.

     In case $m_{1}\neq m_{2}$, $n_{1}=f_{m_{2},n_{2}}(s_{2}^{*})$, then $I_{1}\cap I_{2}=\{(n_{1},n)\mid n\leq \min\{s_{1},n_{2}^{*}\}\}$, which is a Type II ideal.

This covers all possible cases and we have confirmed that the intersection of Type IV and Type V can be a Type I ideal, a Type II ideal, the type $\mathrm{I}\cup \mathrm{II}^{1}$ ideal, $\mathrm{I}\cup \mathrm{II}^{2}$ ideal or the empty set.

The corresponding cell of two Type V ideals says ¡°I/II/V/I$\cup$$ \mathrm{II}^{2}$/$\emptyset$¡±, and that means that the intersection can be a Type I ideal, a Type II ideal, a type V ideal, the Type $\mathrm{I}\cup \mathrm{II}^{2}$ ideal or the empty set.

Suppose $I_{1}=\da(f_{m_{1},n_{1}}(s_{1}),\top)=L_{f_{m_{1},n_{1}}(s_{1})}\cup \{(m_{1},s)\mid s\leq s_{1}\}\cup\{(f_{m_{1},n_{1}}(s_{1}^{*}),n)\mid n\leq n_{1}^{*}\}$ for some $(m_{1},n_{1})\in B, s_{1}=s_{1}^{*}.n_{1}^{*}\in \mathbb{N}^{<\mathbb{N}}$ with $s_{1}^{*}\in \mathbb{N}^{<\mathbb{N}},n_{1}^{*}\in \mathbb{N}$,

    $I_{2}=\da(f_{m_{2},n_{2}}(s_{2}),\top)=L_{f_{m_{2},n_{2}}(s_{2})}\cup \{(m_{2},s)\mid s\leq s_{2}\}\cup\{(f_{m_{2},n_{2}}(s_{2}^{*}),n)\mid n\leq n_{2}^{*}\}$ for some $(m_{2},n_{2})\in B, s_{2}=s_{2}^{*}.n_{2}^{*}\in \mathbb{N}^{<\mathbb{N}}$ with $s_{2}^{*}\in \mathbb{N}^{<\mathbb{N}},n_{2}^{*}\in \mathbb{N}$. Then $I_{1},I_{2}$ are two Type V ideals.

We now distinguish the following cases for $f_{m_{2},n_{2}}(s_{2})$:

    Case 1, $f_{m_{2},n_{2}}(s_{2})<m_{1}$: Then $I_{1}\cap I_{2}=\emptyset$.

    Case 2, $f_{m_{2},n_{2}}(s_{2})=m_{1}$: Then $I_{1}\cap I_{2}=\{(m_{1},s)\mid s\leq s_{1}\}$, which is a Type I ideal.

    Case 3, $m_{1}<f_{m_{2},n_{2}}(s_{2})<f_{m_{1},n_{1}}(s_{1}^{*})$: Then $I_{1}\cap I_{2}=\emptyset$ in case $m_{2}\neq m_{1}$. In case $m_{2}=m_{1}$, if $\da s_{1}\cap \da s_{2}=\emptyset$, then $I_{1}\cap I_{2}=\emptyset$. Otherwise, $I_{1}\cap I_{2}=\{(m_{1},s)\mid s\leq\inf\{s_{1},s_{2}\}$, which is a Type I ideal.

    Case 4, $f_{m_{2},n_{2}}(s_{2})=f_{m_{1},n_{1}}(s_{1}^{*})$: Then $(m_{1},n_{1})=(m_{2},n_{2}), s_{1}^{*}=s_{2}$ by the property of $i$ and $f_{m_{1},n_{1}}$. This implies that $I_{1}\cap I_{2}=\{(m_{1},s)\mid s\leq s_{2})\}\cup\{(f_{m_{2},n_{2}}(s_{2}),n)\mid n\leq n_{1}^{*}\}$ which is a Type $\mathrm{I}\cup \mathrm{II}^{2}$ ideal.

    Case 5, $f_{m_{1},n_{1}}(s_{1}^{*})<f_{m_{2},n_{2}}(s_{2})<f_{m_{1},n_{1}}(s_{1})$: Then the case $m_{1}\neq m_{2},f_{m_{1},n_{1}}(s_{1}^{*})\neq f_{m_{2},n_{2}}(s_{2}^{*})$ implies that $I_{1}\cap I_{2}=\emptyset$.

    In case $m_{1}=m_{2},f_{m_{1},n_{1}}(s_{1}^{*})\neq f_{m_{2},n_{2}}(s_{2}^{*})$, suppose $\da s_{1}\cap \da s_{2}=\emptyset$. Then $I_{1}\cap I_{2}=\emptyset$. Otherwise, $\da s_{1}\cap \da s_{2}\neq \emptyset$. This implies that $ I_{1}\cap I_{2}=\{(m_{1},s)\mid s\leq \inf\{s_{1},s_{2}\}\}$, which is a Type I ideal.

     In case, $m_{1}=m_{2},f_{m_{1},n_{1}}(s_{1}^{*})= f_{m_{2},n_{2}}(s_{2}^{*})$.  Then $(m_{1},n_{1})=(m_{2},n_{2}), s_{1}^{*}=s_{2}^{*}$ by the property of $i$ and $f_{m_{1},n_{1}}$. This implies that $I_{1}\cap I_{2}=\{(m_{1},s)\mid s\leq s_{1}^{*})\}\cup\{(f_{m_{1},n_{1}}(s_{1}^{*}),n)\mid n\leq \min\{n_{1}^{*},n_{2}^{*}\}\}$, which is a Type $\mathrm{I}\cup \mathrm{II}^{2}$ ideal.

     The case $m_{1}\neq m_{2}$, $f_{m_{1},n_{1}}(s_{1}^{*})= f_{m_{2},n_{2}}(s_{2}^{*})$ indicates that $m_{1}=m_{2}$, which contradicts the assumption that $m_{1}\neq m_{2}$.

    Case 6, $f_{m_{2},n_{2}}(s_{2})=f_{m_{1},n_{1}}(s_{1})$: Then $(m_{1},n_{1})=(m_{2},n_{2}),s_{1}=s_{2}$ by the property of $i$ and $f_{m_{1},n_{1}}$. This reveals that $I_{1}\cap I_{2}=I_{1}$, which is a type V ideal.

    Note that the remain case for $f_{m_{2},n_{2}}(s_{2})$ is symmetric with the above cases. This covers all possible cases and we have confirmed that the intersection of two Type V ideals can be a Type I ideal, a Type II ideal, a type V ideal, the type $\mathrm{I}\cup \mathrm{II}^{2}$ ideal or the empty set.

    The Type I$\cup$$ \mathrm{II}^{1}$ ideals, I$\cup$$ \mathrm{II}^{2}$ ideals themselves intersect with sets of Type I$\cup$$ \mathrm{II}^{1}$, I$\cup$$ \mathrm{II}^{2}$, I, II, III, or IV, V as follows:

    \begin{table}[H]
    \footnotesize{
  \centering

  \label{T.smallcontext}
  \vspace{0.2cm}
  \centering
  \scalebox{1}{
  \begin{tabular}{lccccccc}
  \hline
           &  Type I & Type II  & Type III  &  Type IV  &  Type V &  Type I$\cup \mathrm{II}^{1}$ &  Type I$\cup \mathrm{II}^{2}$\\
  \hline
  Type I$\cup \mathrm{II}^{1}$ & I/$\emptyset$      & II/$\emptyset$      & I/II/$\emptyset$      & I/II/I$\cup \mathrm{II}^{1}$/$\emptyset$   & I/II/I$\cup \mathrm{II}^{1}$/$\emptyset$  & I/II/I$\cup \mathrm{II}^{1}$/$\emptyset$ & I/II/I$\cup \mathrm{II}^{1}$/$\emptyset$\\
  \hline
\hline
  Type I$\cup \mathrm{II}^{2}$ & I/$\emptyset$      & II/$\emptyset$      & I/II/$\emptyset$      & I/II/I$\cup \mathrm{II}^{1}$/$\emptyset$   & I/II/I$\cup \mathrm{II}^{2}$/$\emptyset$ &  & I/II/I$\cup \mathrm{II}^{2}$/$\emptyset$\\
  \end{tabular}}}
\end{table}

    The interesting cases are these cases which we discussed.

    The corresponding cell of Type I$\cup \mathrm{II}^{1}$ ideals and Type IV ideals says ¡°I/II/I$\cup$$ \mathrm{II}^{1}$/$\emptyset$¡±, and that means that the intersection can be a Type I ideal, a Type II ideal, the Type $\mathrm{I}\cup \mathrm{II}^{1}$ ideal or the empty set.

    Assume $I_{1}=\{(m_{1},s_{1})\}\cup\{(n_{1},n)\mid n\leq k_{1}\}$ for some $(m_{1},n_{1})\in B, s_{1}\in \mathbb{N}^{<\mathbb{N}},k_{1}\in \mathbb{N}$ with $ {\mid}s_{1}{\mid}=1$, $k_{1}\leq s_{1}$,

    $I_{2}=\da(f_{m_{2},n_{2}}(s_{2}),\top)=L_{f_{m_{2},n_{2}}(s_{2})}\cup \{(m_{2},s_{2})\}\cup\{(n_{2},n)\mid n\leq s_{2}\}$ for some $(m_{2},n_{2})\in B, s_{2}\in \mathbb{N}^{<\mathbb{N}}$ with $ {\mid}s_{2}{\mid}=1$. Then $I_{1}$ is a Type I$\cup \mathrm{II}^{1}$ ideal, $I_{2}$ a type IV ideal.

    The only interesting cases are Type $\mathrm{I}\cup \mathrm{II}^{1}$ ideals, $\mathrm{I}\cup \mathrm{II}^{2}$ ideals.

    We now distinguish two cases:

    Case 1, $m_{1}=f_{m_{2},n_{2}}(s_{2})$: Then $n_{1}=n_{2}$. Note that $m_{1}<n_{1}$. This means that $f_{m_{2},n_{2}}(s_{2})<n_{2}$, which contradicts that $f_{m_{2},n_{2}}(s_{2})>n_{2}$.

    Case 2, $m_{1}= m_{2}$: If $n_{1}=f_{m_{2},n_{2}}(s_{2})$, then $I_{1}\cap I_{2}=I_{1}$. Otherwise, $n_{1}=n_{2}$. Then $I_{1}\cap I_{2}=\{(m_{1},s_{1})\}\cup \{(n_{1},n)\mid n\leq \min\{k_{1},s_{2}\}\}$, which is a Type $\mathrm{I}\cup \mathrm{II}^{1}$ ideal.

The corresponding cell of Type I$\cup \mathrm{II}^{1}$ ideals and Type V ideals says ¡°I/II/I$\cup$$ \mathrm{II}^{1}$/$\emptyset$¡±, and that means that the intersection can be a Type I ideal, a Type II ideal, the Type $\mathrm{I}\cup \mathrm{II}^{1}$ ideal or the empty set.

Assume $I_{1}=\{(m_{1},s_{1})\}\cup\{(n_{1},n)\mid n\leq k_{1}\}$ for some $(m_{1},n_{1})\in B, s_{1}\in \mathbb{N}^{<\mathbb{N}},k_{1}\in \mathbb{N}$ with $ {\mid}s_{1}{\mid}=1$, $k_{1}\leq s_{1}$,

    $I_{2}=\da(f_{m_{2},n_{2}}(s_{2}),\top)=L_{f_{m_{2},n_{2}}(s_{2})}\cup \{(m_{2},s)\mid s\leq s_{2}\}\cup\{(f_{m_{2},n_{2}}(s_{2}^{*}),n)\mid n\leq n_{2}^{*}\}$ for some $(m_{2},n_{2})\in B, s_{2}=s_{2}^{*}.n_{2}^{*}\in \mathbb{N}^{<\mathbb{N}}$ with $s_{2}^{*}\in \mathbb{N}^{<\mathbb{N}},n_{2}^{*}\in \mathbb{N}$. Then $I_{1}$ is a Type $\mathrm{I}\cup \mathrm{II}^{1}$ ideal, $I_{2}$ a Type V ideal.

   The only interesting cases are Type $\mathrm{I}\cup \mathrm{II}^{1}$ ideals, $\mathrm{I}\cup \mathrm{II}^{2}$ ideals.

    We now distinguish two cases:

    Case 1, $m_{1}=f_{m_{2},n_{2}}(s_{2})$: Then $n_{1}=f_{m_{2},n_{2}}(s_{2}^{*})$. Note that $m_{1}<n_{1}$. This means that $f_{m_{2},n_{2}}(s_{2})<f_{m_{2},n_{2}}(s_{2}^{*})$, which contradicts that $f_{m_{2},n_{2}}$ is a monotone injection.

    Case 2, $m_{1}= m_{2}$: If $n_{1}=f_{m_{2},n_{2}}(s_{2})$, then $I_{1}\cap I_{2}=I_{1}$. Otherwise, $n_{1}=f_{m_{2},n_{2}}(s_{2}^{*})$. Then $I_{1}\cap I_{2}=\{(m_{1},s_{1})\}\cup \{(n_{1},n)\mid n\leq \min\{k_{1},n_{2}^{*}\}\}$, which is a Type $\mathrm{I}\cup \mathrm{II}^{1}$ ideal.

    The corresponding cell of two Type I$\cup \mathrm{II}^{1}$ ideals says ¡°I/II/I$\cup$$ \mathrm{II}^{1}$/$\emptyset$¡±, and that means that the intersection can be a Type I ideal, a Type II ideal, the Type $\mathrm{I}\cup \mathrm{II}^{1}$ ideal or the empty set.

Assume $I_{1}=\{(m_{1},s_{1})\}\cup\{(n_{1},n)\mid n\leq k_{1}\}$ for some $(m_{1},n_{1})\in B, s_{1}\in \mathbb{N}^{<\mathbb{N}},k_{1}\in \mathbb{N}$ with $ {\mid}s_{1}{\mid}=1$, $k_{1}\leq s_{1}$,

 $I_{2}=\{(m_{2},s_{2})\}\cup\{(n_{2},n)\mid n\leq k_{2}\}$ for some $(m_{2},n_{2})\in B, s_{2}\in \mathbb{N}^{<\mathbb{N}},k_{2}\in \mathbb{N}$ with $ {\mid}s_{2}{\mid}=1$, $k_{2}\leq s_{2}$. Then $I_{1}$ and $I_{2}$ are two Type $\mathrm{I}\cup \mathrm{II}^{1}$ ideals.

The only interesting cases are Type $\mathrm{I}\cup \mathrm{II}^{1}$ ideals, $\mathrm{I}\cup \mathrm{II}^{2}$ ideals. Then $m_{1}=m_{2}, n_{1}=n_{2}$. This means that $I_{1}\cap I_{2}=\{(m_{1},s_{1})\}\cup\{(n_{1},n)\mid n\leq \min\{k_{1},k_{2}\}\}$, which is a Type $\mathrm{I}\cup \mathrm{II}^{1}$ ideal.

 The corresponding cell of Type I$\cup \mathrm{II}^{1}$ ideals and Type I$\cup \mathrm{II}^{2}$ ideals says ¡°I/II/I$\cup$$ \mathrm{II}^{1}$/$\emptyset$¡±, and that means that the intersection can be a Type I ideal, a Type II ideal, the Type $\mathrm{I}\cup \mathrm{II}^{1}$ ideal or the empty set.

 Assume $I_{1}=\{(m_{1},s_{1})\}\cup\{(n_{1},n)\mid n\leq k_{1}\}$ for some $(m_{1},n_{1})\in B, s_{1}\in \mathbb{N}^{<\mathbb{N}},k_{1}\in \mathbb{N}$ with $ {\mid}s_{1}{\mid}=1$, $k_{1}\leq s_{1}$,

 $I_{2}=\{(m_{2},s)\mid s\leq s_{2}\}\cup\{(f_{m_{2},n_{2}}(s_{2}),n)\mid n\leq k_{2}\}$ for some $(m_{2},n_{2})\in B, s_{2}\in \mathbb{N}^{<\mathbb{N}},k_{2}\in \mathbb{N}$. Then $I_{1}$ is a Type $\mathrm{I}\cup \mathrm{II}^{1}$ ideal, $I_{2}$ is a Type $\mathrm{I}\cup \mathrm{II}^{2}$ ideal.

 The only interesting cases are Type $\mathrm{I}\cup \mathrm{II}^{1}$ ideals, $\mathrm{I}\cup \mathrm{II}^{2}$ ideals. Then $m_{1}=m_{2}$, $n_{1}=f_{m_{2},n_{2}}(s_{2})$. It follows that $I_{1}\cap I_{2}=\{(m_{1},s_{1})\}\cup\{(n_{1},n)\mid n\leq \min\{k_{1},k_{2}\}\}$, which is a Type $\mathrm{I}\cup \mathrm{II}^{1}$ ideal.

 The corresponding cell of Type I$\cup \mathrm{II}^{2}$ ideals and Type IV ideals says ¡°I/II/I$\cup$$ \mathrm{II}^{1}$/$\emptyset$¡±, and that means that the intersection can be a Type I ideal, a Type II ideal, the Type $\mathrm{I}\cup \mathrm{II}^{1}$ ideal or the empty set.

 Assume $I_{1}=\{(m_{1},s)\mid s\leq s_{1}\}\cup\{(f_{m_{1},n_{1}}(s_{1}),n)\mid n\leq k_{1}\}$ for some $(m_{1},n_{1})\in B, s_{1}\in \mathbb{N}^{<\mathbb{N}},k_{1}\in \mathbb{N}$,

    $I_{2}=\da(f_{m_{2},n_{2}}(s_{2}),\top)=L_{f_{m_{2},n_{2}}(s_{2})}\cup \{(m_{2},s_{2})\}\cup\{(n_{2},n)\mid n\leq s_{2}\}$ for some $(m_{2},n_{2})\in B, s_{2}\in \mathbb{N}^{<\mathbb{N}}$ with $ {\mid}s_{2}{\mid}=1$. Then $I_{1}$ is a Type I$\cup \mathrm{II}^{2}$ ideal, $I_{2}$ a type IV ideal.

    The only interesting cases are Type $\mathrm{I}\cup \mathrm{II}^{1}$ ideals, $\mathrm{I}\cup \mathrm{II}^{2}$ ideals.

    We now distinguish two cases:

    Case 1, $m_{1}=f_{m_{2},n_{2}}(s_{2})$: Then $f_{m_{1},n_{1}}(s_{1})=n_{2}$. Note that $m_{1}<f_{m_{1},n_{1}}(s_{1})$. This means that $f_{m_{2},n_{2}}(s_{2})<n_{2}$, which contradicts that $f_{m_{2},n_{2}}(s_{2})>n_{2}$.

    Case 2, $m_{1}= m_{2}$: If $f_{m_{1},n_{1}}(s_{1})=f_{m_{2},n_{2}}(s_{2})$, then $I_{1}\cap I_{2}=I_{1}$. Otherwise, $f_{m_{1},n_{1}}(s_{1})=n_{2}$. Then $I_{1}\cap I_{2}=\{(m_{2},s_{2})\}\cup \{(n_{2},n)\mid n\leq \min\{k_{1},s_{2}\}\}$, which is a Type $\mathrm{I}\cup \mathrm{II}^{1}$ ideal.

  The corresponding cell of Type I$\cup \mathrm{II}^{2}$ ideals and Type V ideals says ¡°I/II/I$\cup$$ \mathrm{II}^{2}$/$\emptyset$¡±, and that means that the intersection can be a Type I ideal, a Type II ideal, the Type $\mathrm{I}\cup \mathrm{II}^{2}$ ideal or the empty set.

Assume $I_{1}=\{(m_{1},s)\mid s\leq s_{1}\}\cup\{(f_{m_{1},n_{1}}(s_{1}),n)\mid n\leq k_{1}\}$ for some $(m_{1},n_{1})\in B, s_{1}\in \mathbb{N}^{<\mathbb{N}},k_{1}\in \mathbb{N}$,

    $I_{2}=\da(f_{m_{2},n_{2}}(s_{2}),\top)=L_{f_{m_{2},n_{2}}(s_{2})}\cup \{(m_{2},s)\mid s\leq s_{2}\}\cup\{(f_{m_{2},n_{2}}(s_{2}^{*}),n)\mid n\leq n_{2}^{*}\}$ for some $(m_{2},n_{2})\in B, s_{2}=s_{2}^{*}.n_{2}^{*}\in \mathbb{N}^{<\mathbb{N}}$ with $s_{2}^{*}\in \mathbb{N}^{<\mathbb{N}},n_{2}^{*}\in \mathbb{N}$. Then $I_{1}$ is a Type $\mathrm{I}\cup \mathrm{II}^{2}$ ideal, $I_{2}$ a Type V ideal.

   The only interesting cases are Type $\mathrm{I}\cup \mathrm{II}^{1}$ ideals, $\mathrm{I}\cup \mathrm{II}^{2}$ ideals.

    We now distinguish two cases:

    Case 1, $m_{1}=f_{m_{2},n_{2}}(s_{2})$: Then $f_{m_{1},n_{1}}(s_{1})=f_{m_{2},n_{2}}(s_{2}^{*})$. Note that $m_{1}<f_{m_{1},n_{1}}(s_{1})$. This means that $f_{m_{2},n_{2}}(s_{2})<f_{m_{2},n_{2}}(s_{2}^{*})$, which contradicts that $f_{m_{2},n_{2}}$ is a monotone injection.

    Case 2, $m_{1}= m_{2}$: If $f_{m_{1},n_{1}}(s_{1})=f_{m_{2},n_{2}}(s_{2})$, then $s_{1}=s_{2}$ from the property of $i$ and $f_{m_{1},n_{1}}$. This means that $I_{1}\cap I_{2}=I_{1}$. Otherwise, $f_{m_{1},n_{1}}(s_{1})=f_{m_{2},n_{2}}(s_{2}^{*})$. Then $(m_{1},n_{1})=(m_{2},n_{2})$, $s_{1}=s_{2}^{*}$. So $I_{1}\cap I_{2}=\{(m_{1},s)\mid s\leq s_{1}\}\cup \{(f_{m_{1},n_{1}}(s_{1}),n)\mid n\leq \min\{k_{1},n_{2}^{*}\}\}$, which is a Type $\mathrm{I}\cup \mathrm{II}^{2}$ ideal.

    The corresponding cell of two Type I$\cup \mathrm{II}^{2}$ ideals says ¡°I/II/I$\cup$$ \mathrm{II}^{2}$/$\emptyset$¡±, and that means that the intersection can be a Type I ideal, a Type II ideal, the Type $\mathrm{I}\cup \mathrm{II}^{2}$ ideal or the empty set.

Assume $I_{1}=\{(m_{1},s)\mid s\leq s_{1}\}\cup\{(f_{m_{1},n_{1}}(s_{1}),n)\mid n\leq k_{1}\}$ for some $(m_{1},n_{1})\in B, s_{1}\in \mathbb{N}^{<\mathbb{N}},k_{1}\in \mathbb{N}$,

 $I_{2}=\{(m_{2},s)\mid s\leq s_{2}\}\cup\{(f_{m_{2},n_{2}}(s_{2}),n)\mid n\leq k_{2}\}$ for some $(m_{2},n_{2})\in B, s_{2}\in \mathbb{N}^{<\mathbb{N}},k_{2}\in \mathbb{N}$. Then $I_{1}$ and $I_{2}$ are two Type $\mathrm{I}\cup \mathrm{II}^{2}$ ideals.

The only interesting cases are Type $\mathrm{I}\cup \mathrm{II}^{1}$ ideals, $\mathrm{I}\cup \mathrm{II}^{2}$ ideals. Then $m_{1}=m_{2}, f_{m_{1},n_{1}}(s_{1})=f_{m_{2},n_{2}}(s_{2})$. It follows that $(m_{1},n_{1})=(m_{2},n_{2}), s_{1}=s_{2}$. This means that $I_{1}\cap I_{2}=\{(m_{1},s)\mid s\leq s_{1}\}\cup\{(f_{m_{1},n_{1}}(s_{1}),n)\mid n\leq \min\{k_{1},k_{2}\}\}$, which is a Type $\mathrm{I}\cup \mathrm{II}^{2}$ ideal.

  This covers all cases to be considered and we conclude that the finite intersection of principle ideals of $P$ consists of Type I through Type V ideals and Type I$\cup \mathrm{II}^{1}$ ideals, Type I$\cup \mathrm{II}^{2}$ ideals or $\emptyset$. It is apparent that $M$ also consists of Type I through Type IV ideals and Type I$\cup$II ideals, $\emptyset$.

 Let $(I_{i})_{i\in I}$ be a directed subset of $M$ without maximum element. Since Type III ideals and Type IV ideals, Type V ideals are maximal element of $M$, we deduce that none of the elements in $(I_{i})_{i\in I}$ belong to Type III ideals or Type IV ideals, Type V ideals. Now we need to distinguish the following two cases:

 Case 1, $\{I_{i}\mid i\in I\}$ belong to Type I ideals or Type II ideals: Then it is pretty easy to see that $(I_{i})_{i\in I}$ is made up of either Type I ideals or Type II ideals. Hence, $\sup_{i\in I}I_{i}$ exists since $P$ is a dcpo.

Case 2, there exists $i_{0}\in I$ such that $I_{i_{0}}$ is a Type $\mathrm{I}\cup \mathrm{II}^{1}$ ideal or a Type $\mathrm{I}\cup \mathrm{II}^{2}$ ideal: Then $\{I_{i}\mid i\geq i_{0}\}$ are Type $\mathrm{I}\cup \mathrm{II}^{1}$ ideals or a Type $\mathrm{I}\cup \mathrm{II}^{2}$ ideals. Note that Type $\mathrm{I}\cup \mathrm{II}^{1}$ ideals can only strictly increase finite times. This means that there exists $i_{1}\geq i_{0}$ such that $(I_{i})_{i\geq i_{1}}$ are Type $\mathrm{I}\cup \mathrm{II}^{2}$ ideals.

For any $i\geq i_{1}$, provide $I_{i}=\{(m_{i},s)\mid s\leq s_{i}\}\cup\{(f_{m_{i},n_{i}}(s_{i}),n)\mid n\leq k_{i}\}$ for some $(m_{i},n_{i})\in B, s_{i}\in \mathbb{N}^{<\mathbb{N}},k_{i}\in \mathbb{N}$. Then $f_{m_{i},n_{i}}(s_{i})=f_{m_{j},n_{j}}(s_{j})$ for any $i,j\geq i_{1}$. This implies that $(m_{i},n_{i})=(m_{j},n_{j}), s_{i}=s_{j}$. Let $m_{i}=m_{0}$, $n_{i}=n_{0}$, $s_{i}=s_{0}$. Then $I_{i}=\{(m_{0},s)\mid s\leq s_{0}\}\cup\{(f_{m_{0},n_{0}}(s_{0}),n)\mid n\leq k_{i}\}$.

We claim that $\sup_{i\in I}I_{i}=\da (f_{m_{0},n_{0}}(s_{0}),\top)$. Obviously, $\da (f_{m_{0},n_{0}}(s_{0}),\top)$ is an upper bound of $(I_{i})_{i\in I}$. For any $I$ is an upper bound of $(I_{i})_{i\in I}$ in $M$, we have that $\bigcup_{i\geq i_{1}}\{(f_{m_{0},n_{0}}(s_{0}),n)\mid n\leq k_{i}\}\subseteq I$. Note that $I$ is a Scott closed set of $P$. Then $\sup_{i\geq i_{1}}(f_{m_{0},n_{0}}(s_{0}),k_{i})=(f_{m_{0},n_{0}}(s_{0}),\top)\in I$. Thus, $M$ is a dcpo.

  \end{proof}
\end{lemma}
\begin{theorem}
  $\Sigma(M)$ is a non-sober countable bounded complete dcpo.
  \begin{proof}
    Clearly, $M$ is countable. From Lemma \ref{a}, it suffices to prove that $\Sigma M$ is not sober. We define $g: P\rightarrow M$ by $g(x)=\da x$ for any $x\in P$. By the proof of Lemma \ref{a}, we can conclude that $g$ is Scott continuous. The Scott irreducibility of $g(P)$ follows immediately from Lemma \ref{b}. Note that $\sup f(P)$ does not exist. As a result, $\Sigma M$ is not a sober space.

  \end{proof}
\end{theorem}
\begin{remark}\label{nonsober}
  Let $R=M\cup\{T\}$. Then $R$ is a countable non-sober complete lattice.
\end{remark}

Note that the non-sober complete lattice $R$ constructed above is not distributive. Thus it remains to know whether there a distributive countable non-sober complete lattice. We now  answer  this problem.
\begin{theorem}
  Let $\mathcal{F}=\{\da F\mid F\subseteq_{fin} R\}$. Then $(\mathcal{F},\subseteq)$ is a countable non-sober distributive complete lattice.
  \begin{proof}
    It is easy to see that $\mathcal{F}$ is a distributive lattice. It remains to prove that $\mathcal{F}$ is a non-sober complete lattice.

    Claim 1: $\mathcal{F}$ is a complete lattice.

    Let $(\da F_{i})_{i\in I}$ be a filter family of subsets of $R$. It suffices to prove that $\bigcap_{i\in I}\da F_{i}\in \mathcal{F}$. Note that $\da F_{i}=\ua^{op}F_{i}$ is a finitely generated upper set in $(R,\geq)$ for any $i\in I$. This means that $(\ua^{op} F_{i})_{i\in I}$ is a filtered family of $(R,\geq)$. Because $R^{op}$ is a dcpo and every element of $R^{op}$ is compact. We conclude that $\bigcap_{i\in F}\da F_{i}=\bigcap_{i\in I}\ua^{op}F_{i}\in \mathcal{F}$ with the help of Rudin's Lemma.

    Claim 2: $R$ is non-sober.

     Define $f : R \rightarrow \mathcal{F}$ by $f(x) = \da x$ and $g : \mathcal{F} \rightarrow R$ by $g(A) = \sup A$ for any $A \in \mathcal{F}$.

     It is evident to confirm that $(f,g)$ is a pair of adjoint. This implies that $g$ is Scott continuous. Note that $g\circ f=id_{R}$ and $f$ is Scott continuous. Therefore, $R$ is a Scott retract of $\mathcal{F}$. Suppose $\mathcal{F}$ is sober. Then $R$ is sober, which contradicts Remark \ref{nonsober}.
  \end{proof}
\end{theorem}

\section{Conclusions}

In this paper we constructed a countable complete lattice whose Scott space is non-sober, thus answered a problem posed by Achim Jung.
Based this complete lattice, we further came up with a countable distributive complete lattice whose Scott space is non-sober. One of the
the useful results we obtained is that if $P$ and $Q$ are dcpos such that $Idl(P)$ and $Idl(Q)$ are countable, then the topology of
$\Sigma P\times \Sigma Q$ coincides with the Scott topology of $P\times Q$.

\section{Reference}
\bibliographystyle{plain}

\end{document}